\documentclass[a4paper,reqno]{amsart}
\usepackage{amsfonts}
\theoremstyle{plain}
\begingroup
\newtheorem{theorem}{Theorem}[section]
\newtheorem{lemma}[theorem]{Lemma}
\newtheorem{proposition}[theorem]{Proposition}

\endgroup

\theoremstyle{definition}
\begingroup
\newtheorem{definition}[theorem]{Definition}
\newtheorem{remark}[theorem]{Remark}

\endgroup

\theoremstyle{remark}
\begingroup

\endgroup

\mathchardef\emptyset="001F

\numberwithin{equation}{section}

\newcommand{\e}{\varepsilon}
\newcommand{\Om}{\Omega}
\newcommand{\om}{\omega}
\newcommand{\Ga}{{\mathit\Gamma}}

\newcommand{\R}{{\mathbb R}}
\newcommand{\Rn}{{\R}^n}
\newcommand{\Rm}{{\R}^m}
\newcommand{\Rp}{{\mathcal R}}

\newcommand{\Mmn}{{\mathbb M}^{m{\times}n}}
\newcommand{\Mnn}{{\mathbb M}^{n{\times}n}}
\newcommand{\wto}{\rightharpoonup}

\newcommand{\setmeno}{\!\setminus\!}
\newcommand{\into}{{\int_{\Omega}}}
\newcommand{\hn}{{\mathcal H}^{n-1}}
\newcommand{\Ln}{{\mathcal L}^n}
\newcommand{\E}{{\mathcal E}}
\newcommand{\Ec}{\E_{}^{el}}
\newcommand{\W}{{\mathcal W}}
\newcommand{\F}{{\mathcal F}}
\newcommand{\G}{{\mathcal G}}

\newcommand{\ki}{_k^i}
\newcommand{\kj}{_k^j}
\newcommand{\kim}{_k^{i-1}}
\newcommand{\kjm}{_k^{j-1}}
\newcommand{\K}{{\mathcal K}}

\newcommand{\subsethn}{\mathrel{\mathop{\smash\subset
\vphantom{=}}\limits^\sim} }
\newcommand{\eqhn}{\mathrel{\mathop{\smash=
\vphantom{\scriptscriptstyle a}}\limits^\sim}}

\newcommand{\aplim}{\mathop{\rm ap\,lim}}
\newcommand{\sigmap}{{\sigma^p}}

\newcommand{\ski}{s\ki}
\newcommand{\skim}{s\kim}
\newcommand{\tki}{t\ki}
\newcommand{\tkim}{t\kim}
\newcommand{\tkj}{t\kj}
\newcommand{\tkjm}{t\kjm}
\newcommand{\etak}{\eta_k}
\newcommand{\zetak}{\zeta_k}

\newcommand{\psiki}{\psi\ki}

\newcommand{\psikj}{\psi\kj}
\newcommand{\psikjm}{\psi\kjm}
\newcommand{\Psik}{\Psi_k}
\newcommand{\Xk}{X_k}
\newcommand{\uk}{u_k}
\newcommand{\uki}{u\ki}

\newcommand{\ukj}{u\kj}
\newcommand{\ukjm}{u\kjm}
\newcommand{\vk}{v_k}
\newcommand{\gki}{\Ga\ki}
\newcommand{\gkim}{\Ga\kim}
\newcommand{\gkj}{\Ga\kj}
\newcommand{\gkjm}{\Ga\kjm}

\newcommand{\phik}{\varphi_k}
\newcommand{\chik}{\chi_k}
\newcommand{\rhokj}{\rho\kj}
\newcommand{\sigmakj}{\sigma\kj}

\newcommand{\taumi}{\tau_m^i}
\newcommand{\taumim}{\tau_m^{i-1}}
\newcommand{\sk}{\sigma_k}

\title[Quasistatic crack growth in finite elasticity]
{Quasistatic crack growth in finite elasticity}
\author[Gianni Dal Maso]{Gianni Dal Maso}
\address[Gianni Dal Maso]{SISSA, Via Beirut 2-4, 34014 Trieste,
Italy}
\email[Gianni Dal Maso]{dalmaso@sissa.it}
\author[Gilles A.~Francfort]{Gilles A.~Francfort}
\address[Gilles A.~Francfort]{LPMTM, Universit\'e Paris 13, Av.~J.B.~Cl\'ement, 93430 Villetaneuse,
France}
\email[Gilles A.~Francfort]{francfor@galilee.univ-paris13.fr}
\author[Rodica Toader]{Rodica Toader}
\address[Rodica Toader]{Dipartimento di Ingegneria Civile, Via delle 
Scienze 208, 33100
 Udine, Italy}
\email[Rodica Toader]{toader@dic.uniud.it}

\begin{document}
\begin{abstract}
In this paper, we prove a new existence result for a variational model of crack growth in brittle materials proposed in \cite{F-M}. We consider the case of 
$n$-dimensional finite elasticity, for an arbitrary $n\ge1$, with a quasiconvex bulk energy and with prescribed boundary deformations and applied loads, both depending on time. 
\end{abstract}
\maketitle
{\small

\bigskip
\keywords{\noindent {\bf Keywords:} variational models,
energy minimization, free-discontinuity problems, quasiconvexity,
crack propagation, quasistatic evolution, brittle fracture, 
Griffith's criterion.}

\bigskip
\subjclass{\noindent {\bf 2000 Mathematics Subject Classification:}
35R35, 74R10, 49Q10, 35A35, 35B30, 35J25.}
}
\bigskip
\bigskip

\begin{section}{INTRODUCTION}

In this paper we present a new existence result for a 
variational model of quasistatic growth for brittle cracks 
introduced in \cite{F-M} and based on Griffith's idea (see 
\cite{Gri}) that the crack growth is determined by the competition 
between the elastic energy of the body and the work needed to 
produce a new crack, or extend an existing one. The main feature of this model is 
that the crack path is not prescribed, but is a result of energy 
balance. In order to obtain our existence theorems in any space 
dimension and for a general bulk energy, we introduce a 
mathematical formulation of the problem in a suitable space of 
functions which may exhibit jump discontinuities on sets of codimension one. 

We now describe the model in more detail.
The reference configuration is a bounded open set $\Om$ of $\Rn$ 
with Lipschitz boundary $\partial\Om=\partial_D\Om\cup\partial_N\Om$, with $\partial_D\Om\cap\partial_N\Om=\emptyset$. On the Dirichlet part $\partial_D\Om$ of the boundary we prescribe the boundary deformation, while on the Neumann part  $\partial_N\Om$ we apply the surface forces. In our formulation, a crack is 
any rectifiable set $\Ga$ contained in $\overline\Om$ and with 
finite $n-1$ dimensional Hausdorff measure. We assume that the work 
done to produce the crack $\Ga$ can be written as
$$
\K(\Ga):=\int_{\Ga\setminus \partial_N\Om}\kappa(x,\nu_\Ga(x))
\,d\hn(x)\,,
$$
where $\nu_\Ga$ is a unit normal vector field on $\Ga$ and $\hn$
is the $n-1$ dimensional Hausdorff measure. The function $\kappa(x,\nu)$ depends on the material and satisfies the standard hypotheses which guarantee the lower semicontinuity of $\K$. Since $\kappa(x,\nu)$ depends on the position $x$ and on the orientation $\nu$, we are able to deal with 
heterogeneous and anisotropic materials (see Subsection~\ref{cracks}).

We adopt the framework of hyperelasticity and assume that the bulk energy of the uncracked part of the body is given by
$$
\W(\nabla 
u):=\int_{\Om\setminus\Ga}W(x,\nabla u(x))\,dx\,,
$$ 
where 
$u\colon\Om\setmeno\Ga\to\Rn$ is the unknown deformation of the body, 
and $W(x,\xi)$ is a given function depending on the material. We only
suppose that $W(x,\xi)$ is quasiconvex with respect to $\xi$ and 
satisfies suitable growth and regularity conditions (see 
Subsection~\ref{bulk}). It is convenient to consider the deformation 
of the uncracked part $\Om\setmeno\Ga$ of the body as a function $u$ 
defined almost everywhere on $\Om$, whose discontinuity set $S(u)$ is 
contained in $\Ga$. An adequate functional setting for these 
deformations is a suitable subspace of the space $GSBV(\Om;\Rn)$ 
introduced in \cite{DG-A} and studied in~\cite{A} (see 
Section~\ref{spaces}).

For every time $t \in [0,T]$, the applied load is given by a system of $t$-dependent body and 
surface forces. We assume that these forces are 
conservative and that their work is given by 

$$
\F(t)(u):=\int_{\Om\setminus\Ga}F(t,x,u(x))\,dx\,, 
\qquad
\G(t)(u):=\int_{\partial_S\Om}G(t,x,u(x))\,d\hn(x)\,,
$$
where 
$\partial_S\Om$ is a subset of $\partial_N\Om$, and $F$ and $G$ satisfy suitable 
regularity and growth conditions (see Subsections~\ref{body} 
and~\ref{surfacef} and Section~\ref{ipgen}). To avoid interactions between cracks and surface 
forces, we impose that all cracks remain at a positive distance from 
$\partial_S\Om$ (see~(\ref{separation}) and 
Remark~\ref{boundforces}).

We adopt the following terminology: an admissible configuration is a pair 
$(u,\Ga)$, where $\Ga$ is an admissible crack and $u$ is an 
admissible deformation with jump set $S(u)$ contained in $\Ga$. The 
total energy of $(u,\Ga)$ at time $t$ is defined 
by
$$
\E(t)(u,\Ga):=\W(\nabla 
u)+\K(\Ga)-\F(t)(u)-\G(t)(u)\,.
$$

For every time $t\in[0,T]$ we prescribe a ``continuous" boundary 
displacement $\psi(t)$ on $\partial_D\Om\setmeno\Ga(t)$, where 
$\Ga(t)$ is the unknown crack at time $t$. We thus assume that $\psi(t)$ 
is the trace on $\partial_D\Om$ of a function in a suitable Sobolev 
space on $\Om$, so that we cannot a priori impose a ``strong discontinuity"
at the boundary, like a jump in the prescribed displacement $\psi(t)$. Moreover, we also assume 
that $t\mapsto\psi(t)$ is sufficiently regular 
(see Subsection~\ref{bounddef}).
The set $AD(\psi(t),\Ga(t))$ of 
admissible deformations with crack $\Ga(t)$ and boundary displacement 
$\psi(t)$ is then defined as the set of deformations $u$ in a 
suitable subspace of $GSBV(\Om;\Rn)$, whose jump set $S(u)$ is 
contained in $\Ga(t)$ and whose trace agrees with $\psi(t)$ on 
$\partial_D\Om\setmeno\Ga(t)$.

In the spirit of Griffith's original theory, a 
minimum energy configuration at time $t$ is an admissible
configuration $(u(t),\Ga(t))$, with $u(t)\in AD(\psi(t),\Ga(t))$, 
such that
$$
\E(t)(u(t),\Ga(t))\le \E(t)(u,\Ga)
$$
for every 
admissible crack $\Ga$ containing $\Ga(t)$ and for every 
$u\in 
AD(\psi(t),\Ga)$. In other words, the energy of 
$(u(t),\Ga(t))$ can not be reduced by choosing a larger crack and, 
possibly, a new deformation with the same boundary condition (see Subsection~\ref{minen}). 

An irreversible quasistatic evolution of minimum energy 
configurations is a function $t\mapsto\!(u(t),\Ga(t))$ which 
satisfies the following conditions:
\begin{itemize}
\smallskip
\item[(a)] static 
equilibrium: for every $t\in[0,T]$ the pair $(u(t),\Ga(t))$ is a 
minimum energy configuration at time $t$;
\smallskip
\item[(b)] 
irreversibility: $\Ga(s)$ is contained in $\Ga(t)$ for 
$0\le s<t\le T$;
\smallskip
\item[(c)] nondissipativity: the derivative of the internal 
energy equals the power of the applied forces.
\smallskip
\end{itemize}
In condition (c) the (loosely named) internal energy is defined by 
$$
\E^{in}(t)(u(t),\Ga(t)):=\W(\nabla u(t))+\K(\Ga(t))\,,
$$
while the power of the external forces is given by
\begin{equation}\label{power}
\begin{array}{c}
\displaystyle
\int_{\partial_D\Om\setminus\Ga(t)}\partial_\xi 
W(x,\nabla u(t))\nu\,\dot\psi(t)\,d\hn+
\int_\Om 
\partial_zF(t,x,u(t))\dot 
u(t)\,dx+{}\\
\displaystyle
{}+\int_{\partial_S\Om} 
\partial_zG(t,x,u(t))\dot u(t)\,d\hn\,,
\end{array}
\end{equation}
where $\nu$ is the outer unit normal to 
$\partial\Om$, $ \dot\psi(t)$ and $\dot u(t)$ denote the time 
derivatives of $\psi(t)$ and $u(t)$, while $\partial_\xi W$, $\partial_z F$, and $\partial_zG$ are the partial derivatives of $W(x,\xi)$, $F(t,x,z)$, and $G(t,x,z)$ with respect to $\xi$ and $z$. 
As $\partial_\xi W(x,\nabla u(t))\nu$ is the boundary 
traction corresponding to the deformation $u(t)$, the first term in 
(\ref{power}) is the power of the surface force which produces the 
boundary displacement $\psi(t)$ on ${\partial_D\Om\setmeno\Ga(t)}$. 

Unfortunately, formula 
(\ref{power}) makes sense only if $\psi(t)$ and $u(t)$ are 
sufficiently regular with respect to $t$ (see Remark~\ref{conservation}), while there are quasistatic 
evolutions that are discontinuous with 
respect to $t$. Therefore we prefer to express the conservation of 
energy in a weaker form, which makes sense even if $u$ is not regular 
(see Subsection~\ref{quasistatic}). 

The main result of this paper 
is the following existence theorem: if $(u_0,\Ga_0)$ is a minimum 
energy configuration at time $t=0$, then there exists an 
irreversible quasistatic evolution
$t\mapsto(u(t),\Ga(t))$ with $(u(0),\Ga(0))=(u_0,\Ga_0)$ (see 
Theorem~\ref{main}). 

As for the hypothesis of the previous 
theorem, we remark that for every initial boundary displacement 
$\psi(0)$ and for every crack $\Ga$ there exists a minimum energy 
 configuration $(u_0,\Ga_0)$ at time $t=0$ with $\Ga_0$ 
containing $\Ga$ (see Theorem~\ref{esistmin}). For special 
initial displacements $u_0$ it is easy to determine the cracks 
$\Ga_0$ such that $(u_0,\Ga_0)$ can be used as initial condition in 
the existence theorem.
For instance, if $\psi(0)$ and $u_0$ coincide with the identity map 
$u_{id}$, then $(u_{id},\Ga_0)$ 
is a minimum energy 
configuration at time $t=0$ for every crack $\Ga_0$, under very 
natural assumptions on $W$, $F$, and $G$ (see 
Remark~\ref{identical}). 

Previous results on this subject have been 
obtained in \cite{DM-T} in the case $n=2$ for a scalar-valued $u$ and 
for $W(\xi)=|\xi|^2$, which corresponds to the antiplane case in 
linear elasticity. In that paper the admissible cracks are assumed to 
be connected, or with a uniform bound on the number of connected 
components. This restriction allows to simplify the mathematical 
formulation of the problem. Indeed, in \cite{DM-T} the cracks are 
assumed to be closed and, consequently the deformations belong to a 
suitable Sobolev space.

These results were extended to the case of 
planar linear elasticity by Chambolle in \cite{Ch}. In both papers 
the existence of a solution is obtained by an approximation argument, 
where the approximating cracks converge in the sense of the Hausdorff 
metric, while the approximating deformation gradients converge strongly 
in $L^2$.

The paper \cite{Fra-Lar} removes the restriction on the 
connected components of $\Ga$ and on the dimension of the space, 
and introduces a weak formulation
in the space $SBV(\Om)$. The function $u$ is still scalar-valued and this 
hypothesis is used to obtain some compactness results which need a 
uniform $L^\infty$-bound that, in the scalar case, can be easily 
obtained by truncation. It also provides a jump transfer theorem that is
instrumental in the present analysis (see Subsection~\ref{jtr}).

In the present paper, we deal 
with the vector case, where the deformation $u$ maps a subset $\Om$ of 
$\Rn$ into $\Rn$ (or, more  generally, into $\Rm$, so as to include the antiplane case when $m=1$). This forces us to introduce a weaker formulation in 
the larger space $GSBV(\Om;\Rm)$, where a compactness theorem holds 
under more general hypotheses. Another new feature of this paper is 
that we consider the case of finite elasticity, with an arbitrary 
quasiconvex bulk energy with polynomial growth, and allow for a large 
class of body and surface forces. In truth however, our formulation is not all
encompassing; it does not allow for constant body loads like gravity, or conservative
surface loads like pressure.

As in prior works \cite{DM-T}, \cite{Ch}, \cite{Fra-Lar}, our result is obtained by time discretization. We fix a sequence of subdivisions $(\tki)_{0\le i\le k}$ of the interval $[0,T]$, with
$ 0=t_k^0<t_k^1<\cdots<t_k^{k-1}<t_k^{k}=T$ and $\lim_{k}\max_{i} (\tki-\tkim)= 0$,
and define by induction an approximate solution $(\uki,\gki)$ at time $\tki$:
let $(u_k^0, \Ga_k^0):=(u_0,\Ga_0)$, and, for $i=1,\ldots,k$, let $(\uki,\gki)$ be a solution of the minimum problem
$$
\min\, \{\E(\tki)(u,\Ga):\gkim\subset\Ga,\ u\in AD(\psi(\tki),\Ga)\}\,,
$$
whose existence can be deduced from the $GSBV$ compactness theorem 
of~\cite{A2}.
For every $t\in[0,T]$ we consider the piecewise constant interpolations
$$
\begin{array}{ccc}
\vphantom{\displaystyle\min_A }
\tau_k(t):=\tki\,,
&\uk(t):=\uki\,,
&\Ga_k(t):=\gki\,,
\end{array}
$$
where $i$ is the largest integer such that $\tki\le t$. The solution $(u(t),\Ga(t))$ of the continuous-time evolution problem will be obtained by passing to the limit in the sequence $(u_k(t),\Ga_k(t))$ as $k\to\infty$.

To this aim we introduce a new notion of convergence of sets, called 
$\sigmap$-convergence, related to the notion of jump sets of $SBV$ functions. We study the main properties of this convergence and prove, in particular, a compactness theorem (see Subsections~\ref{sets} and~\ref{compa}). Using these results we show that there exists a subsequence, still denoted $\Ga_k(t)$, such that $\Ga_k(t)\cup\partial_N\Om$ $\sigmap$-converges to a crack $\Ga(t)$ for every $t\in[0,T]$. Since this subsequence does not depend on~$t$, the cracks $\Ga(t)$ are easily shown to satisfy the irreversibility condition~(b).

We now fix $t\in[0,T]$ and, using the $GSBV$ compactness theorem, we extract a further subsequence of $u_k(t)$, depending on $t$, which converges to some function $u(t)$. The very definition of $\sigmap$-convergence implies that the jump set $S(u(t))$ is contained in $\Ga(t)$, so that $u(t)\in AD(\psi(t),\Ga(t))$. Then we show that $(u(t),\Ga(t))$ is a minimum energy configuration at time $t$ (condition~(a)), using the fact that $(u_k(t),\Ga_k(t))$ satisfies the same property at time $\tau_k(t)$. A crucial tool in the proof of this stability result for minimizers is the jump transfer theorem, established in~\cite{Fra-Lar} in the case of $SBV$ functions, and extended here to the $GSBV$ setting (see Subsections~\ref{jtr} and~\ref{cm}). 

It remains to prove condition (c) on the conservation of energy, in the weak form given in Subsection~\ref{quasistatic}. To this aim
we introduce the functions 
$$
\begin{array}{c} \vphantom{\displaystyle \min_A}
\theta_k(t):= \langle \partial\W(\nabla\uk(t)), \nabla \dot\psi(t) \rangle
- \langle \partial\F(\tau_k(t))(\uk(t)), \dot\psi(t) \rangle - {}\\
 {} - \dot \F(t)(\uk(t))- \langle \partial\G(\tau_k(t))(\uk(t)), 
\dot\psi(t) \rangle -
\dot \G(t)(\uk(t)) \,,
\end{array}
$$
where $ \langle\cdot,\cdot \rangle$ denotes the duality pairing in suitable $L^p$ spaces, $ \partial\F$ and $ \partial\G$ are the differentials of the functionals $\F$ and $\G$ in these function spaces, while $ \dot \F$ and $\dot \G$ are the time derivatives of $\F$ and $\G$. By using the minimality property which defines $(\uki,\gki)$, we prove the fundamental estimate
$$
\E(\tau_k(t))(u_k(t),\Ga_k(t))\le \E(0)(u_0,\Ga_0)+\int_0^{\tau_k(t)}\theta_k(s)\,ds+R_k\,,
$$
with $R_k\to0$ as $k\to\infty$ (see Section~\ref{discrete}).

The main difficulty is to pass to the limit in the first term in the definition of $\theta_k(t)$, since $\nabla\uk(t)$ converges only weakly. We overcome this difficulty by proving a technical result (Lemma~\ref{stressconv}), which shows that convergence of the energies implies convergence of the stresses (see Subsection~\ref{mth}). Since $\E$ is lower semicontinuous, we can pass to the limit in the previous estimate, obtaining the energy inequality
$$
\E(t)(u(t),\Ga(t))\le \E(0)(u_0,\Ga_0)+\int_0^t\theta(s)\,ds\,,
$$
where 
$$
\begin{array}{c} \vphantom{\displaystyle \min_A}
\theta(t):= \langle \partial\W(\nabla u(t)), \nabla \dot\psi(t) \rangle
- \langle \partial\F(t)(u(t)), \dot\psi(t) \rangle - {}\\
 {} - \dot \F(t)(u(t))- \langle \partial\G(t)(u(t)), 
\dot\psi(t) \rangle -
\dot \G(t)(u(t)) \,.
\end{array}
$$

Recalling the weak formulation of condition (c) given in Subsection~\ref{quasistatic}, it remains to show that
$$
\E(t)(u(t),\Ga(t))\ge \E(0)(u_0,\Ga_0)+\int_0^t\theta(s)\,ds\,.
$$
To prove this inequality we consider a sequence of subdivisions $(\ski)_{0\le i\le i_k}$ 
of the interval $[0,t]$, with
$0=s_k^0<s_k^1<\cdots<s_k^{i_k-1}<s_k^{i_k}=t$ and 
$\lim_k\max_i(\ski-\skim)= 0$, and compare $\E(\skim)(u(\skim),\Ga(\skim))$ with $\E(\ski)(u(\ski),\Ga(\ski))$, thanks to the minimality property of $(u(\skim),\Ga(\skim))$ given by condition (a). In this way we obtain an estimate of the form
$$
\E(t)(u(t),\Ga(t))\ge \E(0)(u_0,\Ga_0)+\Theta_k(t)\,,
$$
where $\Theta_k(t)$ is an intricate expression that can be written in terms of  Riemann sums of the functions $\theta(s)$, $\nabla\dot\psi(s)$, $\dot\psi(s)$, 
$\dot \F(s)(u(s))$, and $\dot \G(s)(u(s))$ (see Section~\ref{proof}). Although these functions are only Lebesgue integrable, $\Theta_k(t)$ converges to the integral of $\theta$ on $[0,t]$ for a suitable choice of the subdivisions $(\ski)$ (see Subsections~\ref{aRiem} and~\ref{cRiem}).

Finally we prove a result which can be used to justify the numerical approximation of the quasistatic evolution based on time discretization. Even if the deformation $u(t)$ is not uniquely determined by the crack $\Ga(t)$, for every $t\in[0,T]$ the elastic energies and the crack energies of the discrete-time problems converge to the corresponding energies for the continuous-time problems (see Section~\ref{convergence}). 

\end{section}

\begin{section}{SPACES OF FUNCTIONS WITH BOUNDED VARIATION}\label{spaces}

Throughout the paper $\Ln$ and $\hn$ denote the Lebesgue measure in 
$\Rn$ and the $n-1$ dimensional Hausdorff measure, respectively. 
Unless otherwise specified, the expression {\em almost everywhere\/} 
(abbreviated as {\em a.e.\/}) always refers to~$\Ln$.
If $1\le r\le\infty$ and $E$ is a set, we use the notation 
$\|\cdot\|_r$ or$\|\cdot\|_{r,E}$ for the $L^r$ norm on $E$ with 
respect to $\Ln$ or $\hn$ (or to some other measure as dictated by
 the context).

Given two sets $A,\, B$ in $\Rn$ we write $A\subsethn B$ if
$\hn(A\setmeno B)=0$ and we write $A\eqhn B$ if $\hn(A\triangle B)=0$,
where $A\triangle B:= (A\setmeno B)\cup (B\setmeno A)$ denotes the 
symmetric difference of $A$ and~$B$.

We say that a set $\Ga\subset\Rn$ is {\em rectifiable\/} if there exists a
sequence $\Ga_i$ of $C^1$-manifolds of dimension ${n-1}$ such that
$\Ga\eqhn\bigcup_i\Ga_i$ (these sets are called $(\hn, n-1)$ 
rectifiable in~\cite{F}). A {\em unit normal vector field\/} $\nu$ on 
$\Ga$ is an $\hn$-measurable function $\nu\colon\Ga\to \Rn$, with 
$|\nu(x)|=1$ for $\hn$-a.e.\ $x\in\Ga$, such that $\nu(x)$ is normal 
to $\Ga_i$ for $\hn$-a.e.\ $x\in\Ga_i$ and for every $i$. It is 
well-known that every rectifiable set has a unit normal vector field 
(indeed, infinitely many, since there is no continuity assumption) 
and that the definition does not depend on the decomposition of $\Ga$ 
(see, e.g., \cite[Remark 2.87]{A-F-P}).

Let $U$ be a bounded open set in $\Rn$ and let $u\colon U\to\Rm$ be 
a measurable function. Given $x\in U$ we say that $\tilde 
u(x)\in\Rm$ is the {\em approximate limit\/} of $u$ at $x$, and write 
$\tilde u(x)=\aplim\limits_{y\to x}u(y)$, if for every $\e>0$ we have
\begin{equation}\label{aplim}
\lim_{\rho\to0+}\rho^{-n}{\mathcal L}^n(\{y\in B_\rho(x):|u(y)-\tilde 
u(x)|>\e \})=0\,,
\end{equation}
where $B_\rho(x)$ is the open ball with centre $x$ and radius $\rho$. 
We define the {\em jump set\/} $S(u)$ of $u$ as 
the set of points $x\in U$ where the approximate limit of $u$ does 
not exist.
Given $x\in U$ such that $\tilde u(x)$ exists, we say that the 
${m{\times}n}$ matrix $\nabla u(x)$ is the {\em approximate 
differential\/} of $u$ at $x$ if
$$
\aplim\limits_{y\to x}\frac{u(y)-\tilde u(x)-\nabla u(x)(y-x)}{|y-x|}=0\,.
$$

The space $BV(U;\Rm)$ of {\em functions of bounded variation\/} is 
defined as the set of all $u\in L^1(U;\Rm)$ such that the 
distributional gradient $Du$ is a bounded Radon measure on $U$ with 
values in the space $\Mmn$ of ${m{\times}n}$ matrices. If $u\in 
BV(U;\Rm)$, we can consider the Lebesgue decomposition 
$Du=D^au+D^su$, where $D^au$ is absolutely continuous with respect to 
$\Ln$ and $D^su$ is singular with respect to~$\Ln$. In this case the 
approximate differential $\nabla u(x)$ exists for a.e.\ $x\in U$ and 
the function $\nabla u$ belongs to $L^1(U;\Mmn)$ and coincides a.e.\ 
with the
density of $D^au$ with respect to $\Ln$ (Calder\'on-Zygmund Theorem, 
see, e.g., \cite[Theorem 3.83]{A-F-P}). Note that $S(u)$ coincides
with the complement of the set of Lebesgue points for $u$, up to
a set of $\hn$-measure $0$.

The space $SBV(U;\Rm)$ of {\em special functions of bounded 
variation\/} is defined as the set of all $u\in BV(U;\Rm)$ such that 
$D^su$ is concentrated on $S(u)$, i.e., $|D^su|({U\setmeno S(u)})=0$. 
As usual, $SBV_{loc}(U;\Rm)$ denotes the space of functions which 
belong to $SBV(U';\Rm)$ for every open set $U'\subset\subset U$.

Let us fix an exponent $p$, with $1<p<+\infty$. The space 
$SBV^p(U;\Rm)$ is defined as the set of functions $u\in SBV(U;\Rm)$ 
with $\nabla u\in L^p(U;\Mmn)$ and $\hn(S(u))<+\infty$. 

\begin{definition}\label{defconv}
A sequence  $u_k$ converges to  $u$  weakly in 
$SBV^p(U;\Rm)$ if and only if
  $u_k,\,u\in SBV^p(U;\Rm)\cap L^\infty(U;\Rm)$, $u_k\to u$ a.e.\ in $U$, 
$\nabla u_k\wto\nabla u$ weakly in $L^p(U;\Mmn)$, and 
$\|u_k\|_\infty$ and $\hn(S(u_k))$ are bounded uniformly with respect 
to~$k$.
\end{definition}

 If $u\in W^{1,p}(U;\Rm)$, then $u\in SBV^p(U;\Rm)$ and 
$S(u)\eqhn\emptyset$.
The following compactness theorem is proved in \cite{A} (see also 
\cite[Section~4.2]{A-F-P}).

\begin{theorem}\label{compsbv}
Let $u_k$ be a sequence in $SBV^p(U;\Rm)$ such that $\|u_k\|_{\infty}$,
$\|\nabla u_k\|_p$, and
$\hn(S(u_k))$ are bounded uniformly with
respect to $k$. Then there exists a subsequence which converges 
weakly in $SBV^p(U;\Rm)$.
\end{theorem}

This result is not enough for the study of the fracture problem in 
dimension $n$, because we have no a priori bound on the $L^\infty$ 
norm of the solutions. To overcome this difficulty we have to use the 
wider space $GSBV(U;\Rm)$ of {\em generalized special functions of 
bounded variation\/}, defined as the set of all functions $u\colon 
U\to\Rm$ such that $\varphi(u)\in SBV_{loc}(U,\Rm)$ for every 
$\varphi\in C^1(\Rm;\Rm)$ with ${\rm supp}(\nabla 
\varphi)\subset\subset\Rm$. It is easy to see that $SBV(U;\Rm)\subset 
GSBV(U;\Rm)$ and
$GSBV(U;\Rm)\cap L^\infty(U;\Rm)=SBV_{loc}(U;\Rm)\cap L^\infty(U;\Rm)$.
If $u\in GSBV(U;\Rm)$, then the approximate differential $\nabla 
u(x)$ exists for a.e.\ $x\in U$ (see \cite[Propositions 1.3 
and~1.4]{A3}).

We define $GSBV^p(U;\Rm)$ as the set of functions $u\in GSBV(U;\Rm)$ 
such that $\nabla u\in L^p(U;\Mmn)$ and $\hn(S(u))<+\infty$. If $u\in 
GSBV^p(U;\Rm)$ and $\varphi\in C^1(\Rm;\Rm)$ with ${\rm supp}(\nabla 
\varphi)\subset\subset\Rm$, then the function $v:=\varphi(u)$ 
belongs to $SBV_{loc}(U;\Rm)$ and $S(v)\subsethn S(u)$. As $\nabla 
v=\nabla\varphi(u)\nabla u$
a.e.\ in $U$, we have $\nabla v\in L^p(U;\Mmn)$. Since by 
\cite[Section~3.9]{A-F-P}
$$
|Dv|(U)\le \int_U |\nabla v(x)|\,dx+2\|v\|_{\infty}\hn(S(v))\,,
$$
we deduce that $v\in BV(U;\Rm)$, and recalling that $v\in SBV_{loc}(U;\Rm)$
we conclude that $v\in SBV^p(U;\Rm)$. The previous discussion shows that  $SBV^p(U;\Rm)\subset 
GSBV^p(U;\Rm)$ and
$GSBV^p(U;\Rm)\cap L^\infty(U;\Rm)=SBV^p(U;\Rm)\cap L^\infty(U;\Rm)$.
As usual we set $SBV^p(U):=SBV^p(U;\R)$ and $GSBV^p(U):=GSBV^p(U;\R)$.

The following proposition proves some basic properties of the space 
$GSBV^p(U;\Rm)$.
Note that the same properties do not hold for $GSBV(U;\Rm)$ (see 
\cite[Remark 4.27]{A-F-P}).

\begin{proposition}\label{vectorspace}
$GSBV^p(U;\Rm)$ is a vector space. A function 
$u:=(u^1,\ldots,u^m)\colon U\to\Rm$ belongs to $GSBV^p(U;\Rm)$ if and 
only if each component $u^i$ belongs to $GSBV^p(U)$.
\end{proposition}
\begin{proof}
Let $u,\, v\in GSBV^p(U;\Rm)$ and let $\varphi\in C^1(\Rm;\Rm)$ with 
${\rm supp}(\nabla \varphi)\subset\subset\Rm$. We have to prove that 
the function $w:=\varphi(u+v)$ belongs to $SBV_{loc}(U;\Rm)$. To this 
aim we consider a function
$\varphi_1\in C^1_c(\Rm;\Rm)$ such that $\varphi_1(z)=z$
for $|z|\le 1$, and we define $\varphi_k(z):=k\varphi_1(z/k)$. Then
$\varphi_k\in C^1_c(\Rm;\Rm)$,
$\varphi_k(z)=z$ for $|z|\le k$, and
$|\nabla \varphi_k|\le C$ for some constant $C$ independent of $k$.
Let $w_k:=\varphi(\varphi_k(u)+\varphi_k(v))$. Since $\varphi_k(u)$ 
and $\varphi_k(v)$ belong to $SBV^p(U;\Rm)$, the functions $w_k$ belong to 
$SBV^p(U;\Rm)$. As
$$
\nabla w_k=\nabla \varphi(\varphi_k(u)+\varphi_k(v))[\nabla 
\varphi_k(u)\nabla u+\nabla \varphi_k(v)\nabla v]\,,
$$
the sequence $\nabla w_k$ is bounded in $L^p(U;\Mmn)$. Moreover 
$S(w_k)\subsethn S(u)\cup S(v)$ and 
$\|w_k\|_{\infty}\le\|\varphi\|_{\infty}<+\infty$. Since $w_k$ 
converges to $w$ a.e.\ in $U$, from Theorem~\ref{compsbv} we deduce 
that $w\in SBV^p(U;\Rm)$.

Let $u:=(u^1,\ldots,u^m)$ be a function in $GSBV^p(U;\Rm)$,
let $i=1,\ldots,m$, and let $\psi\in C^1(\R)$ with
${\rm supp}\,\psi'\subset\subset\R$. In order to prove that $u^i\in 
GSBV^p(U)$ it is enough to show that $\psi(u^i)$ belongs to $SBV(U)$. 
Let $\psi_k$ be a sequence in $C^1_c(\Rm)$ such that 
$\|\psi_k\|_{\infty}\le 1$, $\|\nabla \psi_k\|_{\infty}\le 1$, and 
$\psi_k(z)=1$ for $|z|\le k$, and let $v_k:=\psi(u^i)\psi_k(u)$. Then 
$v_k\in SBV^p(U;\Rm)$, $\|v_k\|_{\infty}\le\|\psi\|_{\infty}<+\infty$,
$S(v_k)\subsethn S(u)$, and
$$
\nabla v_k=\psi'(u^i)\psi_k(u)\nabla u^i+\psi(u^i)\nabla \psi_k(u)\nabla u\,,
$$
so that $\nabla v_k$ is bounded in $L^p(U;\Rm)$. Since $v_k\to 
\psi(u^i)$ a.e.\ in $U$, from Theorem~\ref{compsbv} we deduce that 
$\psi(u^i)\in SBV^p(U)$.

Conversely, if all components $u^i$ of $u$ belong to $GSBV^p(U)$, 
then it is easy to see that $u^ie_i$ belongs to $GSBV^p(U;\Rm)$, 
$e_i$ being the $i^{th}$ vector of the canonical basis of $\Rm$, so 
that $u\in GSBV^p(U;\Rm)$ by the vector space property.
\end{proof}

{}From Proposition~\ref{vectorspace} and \cite[Theorems 4.34 and 
4.40]{A-F-P} we obtain that $S(u)$ is rectifiable for every $u\in 
GSBV^p(U;\Rm)$, and, if $\nu_u$ is a unit normal vector field on 
$S(u)$, then for $\hn$-a.e.\ $x\in S(u)$ there exist two distinct 
vectors $u^+(x),\, u^-(x)\in\Rm$ such that \begin{equation}\label{upm}
u^+(x)=
\aplim\limits_{y\to x,\,y\in H^+(x)}u(y)\,,
\qquad
u^-(x)=
\aplim\limits_{y\to x,\,y\in H^-(x)}u(y)\,,
\end{equation}
where $H^+(x):=\{y\in U:(y-x)\cdot\nu_u(x)>0\}$ and
$H^-(x):=\{y\in U:(y-x)\cdot\nu_u(x)<0\}$.

We introduce now the notion of trace on the boundary in the $GSBV^p$ setting.

\begin{proposition}\label{traces}
If $U$ has a Lipschitz boundary and $u\in GSBV^p(U;\Rm)$, then there 
exists a function $\tilde u\colon\partial U\to \Rm$ such that
\begin{equation}\label{tracex}
\aplim\limits_{y\to x,\,y\in U}u(y)=\tilde u(x)
\end{equation}
for $\hn$-a.e.\ $x\in \partial U$.
\end{proposition}

The $\hn$-a.e.\ defined function $\tilde u\colon\partial U\to \Rm$ is 
called the trace of $u$ on $\partial U$ and in the rest of the paper 
will be denoted simply by $u$.

\begin{proof}[Proof of Proposition~\ref{traces}]
Let us fix a bounded open set $U_0$ containing $\overline U$. Given 
$u\in GSBV^p(U;\Rm)$, let $u_0$ be the function defined by $u_0:=u$ on 
$U$ and $u_0:=0$ on ${U_0\setminus U}$. For every $\varphi\in 
C^1(\Rm;\Rm)$, with ${\rm supp}(\nabla \varphi)\subset\subset\Rm$, 
we have $\varphi(u)\in SBV^p(U;\Rm)$. Therefore $\varphi(u_0)\in 
SBV^p(U_0;\Rm)$ and we conclude that $u_0\in GSBV^p(U_0;\Rm)$.
By the definition of $S(u_0)$, for every $x\in \partial U\setmeno S(u_0)$
we have
$$
\aplim\limits_{y\to x}u_0(y)=\tilde u_0(x)\,,
$$
which implies (\ref{tracex}) with $\tilde u(x):=\tilde u_0(x)$.
We can choose a unit normal vector field $\nu_{u_0}$ on $S(u_0)$ 
which coincides with the outward unit normal to $\partial U$
$\hn$-a.e.\ on $S(u_0)\cap\partial U$ (see, e.g., \cite[Proposition 
2.85]{A-F-P}). Since in this case
$$
\lim_{\rho\to0+}\rho^{-n}\Ln(B_{\rho}(x)\cap(H^-(x)\triangle U))=0
$$
for $\hn$-a.e.\ $x\in S(u_0)\cap\partial U$,
the second equality in (\ref{upm}) gives (\ref{tracex}) with $\tilde 
u(x):=u_0^-(x)$.
\end{proof}

For every $q\ge 1$ we set $GSBV^p_q(U;\Rm):=GSBV^p(U;\Rm)\cap L^q(U;\Rm)$.

\begin{lemma}\label{rem3} Assume that $U$ has a Lipschitz boundary and that
$u\in GSBV^p_q(U;\Rm)$ for some $q\ge 1$. If
$S(u)\eqhn\emptyset$, then $u$ belongs to
$W^{1,p}(U;\Rm)\cap L^q(U;\Rm)$.
\end{lemma}

\begin{proof}
Let
$\varphi\in C^1_c(\Rm;\Rm)$ be a function such that $\varphi(z)=z$
for $|z|\le 1$, and let $\varphi_k(z):=k\varphi(z/k)$. Then
$\varphi_k\in C^1_c(\Rm;\Rm)$,
$\varphi_k(z)=z$ for $|z|\le k$, and
$|\nabla \varphi_k|\le C$ for some constant $C$ independent of $k$. Under
our assumptions on $u$, the functions $v_k:=\varphi_k(u)$ belong to
$SBV^p(U;\Rm)\cap L^\infty(U;\Rm)$ and
$S(v_k)\eqhn\emptyset$. This implies that $v_k\in
W^{1,p}(U,\Rm)\cap L^\infty(U;\Rm)$. Since $U$ has a Lipschitz 
boundary there exists a constant $\gamma>0$, depending only on $p$, 
$q$, and $U$, such that
$\|v_k\|_p \le \gamma(\|\nabla v_k\|_p + \|v_k\|_q)$ for every~$k$.
As
$|\varphi_k(z)|\le C |z|$ for every $z\in \Rm$ and every $k$, we have
$\|v_k\|_q\le C\|u\|_q$.
Since $\nabla v_k=\nabla\varphi_k(u)\nabla u$, we have also $\|\nabla 
v_k\|_p\le C \|\nabla u\|_p$. Therefore $v_k$ is bounded in
$W^{1,p}(U;\Rm)$. As $v_k$ converges to $u$ pointwise a.e.\ on
$U$, we conclude that $u\in W^{1,p}(U;\Rm)$.
\end{proof}

In the spirit of Definition \ref{defconv}, we introduce the following notion of convergence.

\begin{definition}
A sequence $u_k$ converges to $u$ weakly in 
$GSBV^p(U;\Rm)$ if and only if $u_k,\, u$ belong to $GSBV^p(U;\Rm)$, $u_k\to u$ 
a.e.\ in $U$, $\nabla u_k\wto\nabla u$ weakly in $L^p(U;\Mmn)$, and 
$\hn(S(u_k))$ is bounded uniformly with respect to~$k$.
 \end{definition}
  It is immediate 
that weak convergence in $SBV^p(U;\Rm)$ implies weak convergence in 
$GSBV^p(U;\Rm)$.
The following compactness theorem for $GSBV^p(U;\Rm)$ is proved in
\cite[Theorem 2.2]{A2} (see also \cite[Section~4.5]{A-F-P}).

\begin{theorem}\label{compgsbv} Let $u_k$ be a sequence in 
$GSBV^p(U;\Rm)$ such that $\|u_k\|_1$,
$\|\nabla u_k\|_p$, and
$\hn(S(u_k))$ are bounded uniformly with
respect to~$k$. Then there exists a subsequence which converges 
weakly in $GSBV^p(U;\Rm)$.
\end{theorem}

We recall that a function $W\colon \Mmn\to\R$ is said to be {\em 
quasiconvex\/} if
$$
\int_U W(\xi+\nabla \varphi(x))\,dx\ge W(\xi)\Ln(U)
$$
for every $\xi \in \Mmn$ and every $\varphi\in C^1_c(U;\Rm)$.
The following theorem collects the lower semicontinuity results with 
respect to weak convergence in $GSBV^p(U;\Rm)$ that we shall use in 
the rest of the paper.

\begin{theorem}\label{semigsbv}
Let $W\colon U\times \Mmn\to\R$ be a Carath\'eodory function satisfying
\begin{eqnarray}
&W(x,\cdot) \hbox{ is quasiconvex on }\Mmn \hbox{ for every }x\in 
U\,,\label{quasic}\\
&a_0|\xi|^p-b_0(x)\le W(x,\xi)\le a_1|\xi|^p+b_1(x)\quad\hbox{for 
every }(x,\xi)\in U\times \Mmn\label{growth}
\end{eqnarray}
for some constants $a_0>0$, $a_1>0$, and some nonnegative functions 
$b_0,\, b_1\in L^1(U)$. Let $\kappa\colon U{\times}\Rn\to\R$ be a 
lower semicontinuous function such that
\begin{eqnarray}
&\kappa(x,\cdot) \hbox{ is a norm on }\Rn \hbox{ for every }x\in U\,,
\label{norm}\\
&
\kappa_1|\nu|\le \kappa(x,\nu )\le \kappa_2|\nu| \quad\hbox{for every 
}(x,\nu)\in U{\times}\Rn
\label{kappa}
\end{eqnarray}
for some constants $\kappa_1>0$ and $\kappa_2>0$.
If $u_k$ converges to $u$ weakly in $GSBV^p(U;\Rm)$, then
for every $\hn$-measurable set $E$, with $\hn(E)<+\infty$, we have
\begin{eqnarray}
& \displaystyle\int_UW(x,\nabla u(x))\,dx\le
\liminf_{k\to\infty}\int_UW(x,\nabla u_k(x))\,dx\,,\label{W}\\
&\displaystyle\int_{S(u)\setminus E}\kappa(x,\nu_u(x) ) \,d\hn(x)\le
\liminf_{k\to\infty}
\int_{S(u_k)\setminus E} \kappa(x,\nu_{u_k}(x) ) \,d\hn(x)\,, \label{setE}
\end{eqnarray}
where $\nu_{u_k}$ and $\nu_u$ are unit normal vector fields on 
$S(u_k)$ and $S(u)$, respectively.
\end{theorem}

\begin{proof} If $u_k$ converges to $u$ weakly in $SBV^p(U;\Rm)$, 
inequality (\ref{W}) is proved in \cite{A3} (see also \cite[Theorem 
5.29]{A-F-P}). The general case is proved in~\cite{Kri}.

If $E=\emptyset$, the proof of (\ref{setE}) can be found in 
\cite[Theorem 3.7]{A2} when $\kappa$ does not depend on~$x$. The 
extension to the case of a general $\kappa$ can be obtained by 
standard localization techniques. When $E\neq\emptyset$ is compact 
it is enough to replace $U$ by $U\setmeno E$. To
prove (\ref{setE}) in the general case let $\e>0$ and let $K\subset 
E$ be a compact set such that $\hn(E\setmeno K)<\e$.
Since $S(u)\setmeno E\subset S(u)\setmeno K$ and $S(u_k)\setmeno 
K\subset (S(u_k)\setmeno E)\cup (E\setmeno K)$, we have
\begin{eqnarray*}
&\displaystyle\int_{S(u)\setminus E}\kappa(x,\nu_u (x)) \,d\hn(x)\le
\int_{S(u)\setminus K}\kappa(x,\nu_u(x) ) \,d\hn(x)\le\\
&\displaystyle \le
\liminf_{k\to\infty} \int_{S(u_k)\setminus K} \kappa(x,\nu_{u_k}(x) ) 
\,d\hn(x)\le \\
&\displaystyle\le \liminf_{k\to\infty} \int_{S(u_k)\setminus E} 
\kappa(x,\nu_{u_k}(x) ) \,d\hn(x)+\kappa_2\hn (E\setmeno K)\le
\\
& \displaystyle\le \liminf_{k\to\infty} \int_{S(u_k)\setminus E} 
\kappa(x,\nu_{u_k}(x) ) \,d\hn(x)+\kappa_2\e\,. 
\end{eqnarray*}
As $\e\to 0$ we obtain (\ref{setE}).
\end{proof}

\begin{remark}\label{Su}
Let $E$ be an $\hn$-measurable set with $\hn(E)<+\infty$. 
Theorem~\ref{semigsbv} implies that,
if $u_k$ converges to $u$ weakly in $GSBV^p(U;\Rm)$ and 
$S(u_k)\subsethn E$ for every~$k$, then $S(u)\subsethn E$.
\end{remark}

\begin{remark}\label{union}
Let $\Ga$ be a rectifiable subset of $U$ with $\hn(\Ga)<+\infty$ and let $E$ be an $\hn$-measurable set with $\hn(E)<+\infty$. Since 
$(S(u)\cup\Ga)\setmeno E=
(S(u)\setmeno (\Ga\cup E))\cup(\Ga\setmeno E)$ and
$(S(u_k)\cup\Ga)\setmeno E=
(S(u_k)\setmeno (\Ga\cup E))\cup(\Ga\setmeno E)$,
Theorem~\ref{semigsbv} implies that,
if $u_k$ converges to $u$ weakly in $GSBV^p(U;\Rm)$, then
$$
\int_{(S(u)\cup\Ga)\setminus E}\kappa(x,\nu(x) ) \,d\hn(x)\le
\liminf_{k\to\infty}
\int_{(S(u_k)\cup\Ga)\setminus E} \kappa(x,\nu_k(x) ) \,d\hn(x)\,,
$$
where $\nu$ and $\nu_k$ are unit normal vector fields on $S(u)\cup\Ga$ and $S(u_k)\cup\Ga$, respectively.
\end{remark}

\end{section}

\begin{section}{FORMULATION OF THE PROBLEM}\label{formulation}

\subsection{The reference configuration}
Let $\Om$ be a bounded open set in $\Rn$ with Lipschitz boundary
$\partial\Om$, and let $\Om_B$ be an open subset of $\Om$ with
Lipschitz boundary. The set $\overline\Om$ represents the {\em
reference configuration\/} of an elastic body with cracks, while
$\overline\Om_B$ represents its {\em brittle part\/}, in
the sense that every crack in the reference configuration will be
contained in $\overline\Om_B$.

We fix a closed subset $\partial_N\Om$ of $\partial\Om$, called the
{\em Neumann part\/} of the boundary, on which we
will prescribe the boundary forces. On the {\em Dirichlet part\/} of the
boundary
$\partial_D\Om:=\partial\Om\setmeno\partial_N\Om$ we will prescribe the
boundary deformation, that will be attained only in the part of
$\partial_D\Om$ which is not contained in the crack.
We fix also a closed
subset $\partial_S\Om$ of $\partial_N\Om$, which will contain the
support of all boundary forces applied to the body. We assume that
\begin{equation}\label{separation}
\overline\Om_B\cap\partial_S\Om=\emptyset\,.
\end{equation}
The reason for such a condition will be explained later (see Remark~\ref{boundforces}).

\subsection{The cracks}\label{cracks}
A {\em crack \/} is represented in the reference configuration by a
rectifiable set $\Ga \subsethn \overline\Om_B$ with
$\hn(\Ga)<+\infty$. The collection of 
{\em admissible cracks\/} is given by
\begin{equation}\label{calR}
\Rp(\overline\Om_B):=
\{\Ga: \Ga\,\hbox{ is rectifiable, }\,
\Ga\subsethn\overline\Om_B\,,\
\hn(\Ga)<+\infty\}\,.
\end{equation}
The set
$\Ga\cap\partial_D\Om$ is interpreted as the part of
$\partial_D\Om$ where the prescribed boundary deformation is not
attained. On the contrary $\Ga\cap\partial_N\Om$ will not produce any effect, since $\overline\Om_B\cap\partial_N\Om$ is
traction free by (\ref{separation}).
 
According to Griffith's theory, we assume that the {\em energy spent to
produce the crack\/} $\Ga\in
\Rp(\overline\Om_B)$ is given by
\begin{equation}\label{calK}
\K(\Ga):=\int_{\Ga\setminus\partial_N\Om}\kappa(x,\nu_\Ga(x))
\,d\hn(x)\,,
\end{equation}
where $\nu_\Ga$ is a unit normal vector field on
$\Ga$ and $\kappa\colon\overline\Om_B{\times}\Rn \to\R$
is a lower semicontinuous function, which takes into account the {\em
toughness\/} of the material in different locations and in different
directions. Note that, since with our definition $\K(\Ga)=\K({\Ga\setmeno\partial_N\Om})$, there is no energy associated with the part of the crack that lies on $\partial_N\Om$. Nevertheless for mathematical convenience (see Section~\ref{preliminary}) it is sometimes useful to consider also cracks $\Ga$ with $\Ga\cap\partial_N\Om\neq\emptyset$. As in Theorem \ref{semigsbv}, we assume that
\begin{eqnarray}
&\kappa(x,\cdot) \hbox{ is a norm on }\Rn \hbox{ for every }x\in
\overline\Om_B\,,\label{K1}\\
& \kappa_1|\nu| \le\kappa(x,\nu)\le\kappa_2|\nu| \hbox{ for every }
(x,\nu)\in\overline\Om_B{\times}\Rn\,,\label{K2}
\end{eqnarray}
for some constants $\kappa_1>0$ and $\kappa_2>0$.

To simplify the exposition of  auxiliary results, we extend $\kappa$ to $\Rn\!{\times}\Rn$ by setting $\kappa(x,\nu)\!:= \kappa_2|\nu| $ if $x\notin\overline\Om_B$, and we define $\K(\Ga)$ by (\ref{calK}) for every rectifiable subset $\Ga$ of $\Rn$.
By (\ref{K2}) we have
\begin{equation}\label{K9}
\kappa_1\hn(\Ga\setmeno\partial_N\Om)\le \K(\Ga) \le
\kappa_2\hn(\Ga\setmeno\partial_N\Om)
\end{equation}
for every rectifiable subset $\Ga$ of $\Rn$.

\subsection{The body deformations and their bulk energy}\label{bulk}
Given an admissible crack $\Ga$, an {\em admissible
deformation\/} with crack $\Ga$ will be any function
$u\in GSBV(\Om;\Rm)$ with $S(u)\subsethn\Ga$. This
implies that $u$ has a representative $\tilde u$ which coincides with
$u$ a.e.\ on $\Om$, is defined at $\hn$-a.e.\ point of ${\Om\setmeno
\Ga}$, and is approximately continuous $\hn$-a.e.\ on ${\Om\setmeno
\Ga}$, in the sense that
$$
\tilde u(x)=\aplim_{y\to x,\, y\notin\Ga}\tilde u(y)
$$
for $\hn$-a.e.\ $x\in{\Om\setmeno \Ga}$. Note that when $m=n$ we are in the classical case of finite elasticity on $\Om\setmeno\Ga$, and when $m=1$ we are in the antiplane setting.

We assume that the uncracked part of the body is hyperelastic and that its
{\em bulk energy\/} relative to the deformation
$u\in GSBV(\Om;\Rm)$ can be written as
$$
\int_\Om W(x,\nabla u(x))\,dx\,,
$$
where $W\colon\Om\times\Mmn\to\R$ is a given Carath\'eodory
function such
that $\xi\mapsto W(x,\xi)$ is quasiconvex and $C^1$ on $\Mmn$
for every $x\in\Om$. As in Theorem \ref{semigsbv}, we assume that there exist three constants $p>1$,
$a_0^W>0$, $a_1^W>0$, and
two nonnegative functions $b_0^W$, $b_1^W\in L^1(\Om)$, such that
\begin{equation}\label{alphaW}
a_0^W|\xi|^p-b_0^W(x)\le W(x,\xi)\le a_1^W|\xi|^p+b_1^W(x)
\end{equation}
for every $(x,\xi)\in\Om{\times}\Mmn$. Since $\xi\mapsto
W(x,\xi)$ is rank-one convex on $\Mmn$ for every $x\in\Om$ (see, e.g.,
\cite{Dac}), we  deduce  from the
previous inequalities an estimate for the partial gradient
$\partial_\xi W\colon\Om\times\Mmn\to\Mmn$ of $W$
with respect to $\xi$ (see, e.g.,
\cite{Dac}). Specifically,
there exist a constant $a_2^W>0$
and a nonnegative function $b_2^W\in L^{p'}(\Om)$, with $p'=p/(p-1)$,
such that
\begin{equation}\label{alpha2W}
|\partial_\xi W(x,\xi)|\le a_2^W|\xi|^{p-1}+b_2^W(x)
\end{equation}
for every $(x,\xi)\in\Om{\times}\Mmn$.

Note that in the case $m=n$ the boundedness assumption (\ref{alphaW}) prohibits the introduction of the ``classical'' constraint that $W(\xi)\to\infty$ as $\det\xi\to 0$.

To shorten the notation we introduce the function $\W\colon
L^p(\Om;\Mmn)\to\R$ defined by
\begin{equation}\label{calW}
\W(\Phi):=\int_\Om W(x,\Phi(x))\,dx
\end{equation}
for every $\Phi\in L^p(\Om;\Mmn)$. By (\ref{alphaW}) and
(\ref{alpha2W}) the functional $\W$ is of class $C^1$ on
$L^p(\Om;\Mmn)$ and its differential $\partial \W\colon L^p(\Om;\Mmn)\to
L^{p'}(\Om;\Mmn)$ is given by
\begin{equation}\label{W7.5}
\langle \partial \W(\Phi), \Psi\rangle =\int_\Om \partial_\xi
W(x,\Phi(x))\Psi(x)\,dx\,,
\end{equation}
for every $\Phi,\, \Psi\in L^p(\Om;\Mmn)$, where
$\langle\cdot,\cdot\rangle$ denotes the duality pairing between the
spaces $L^{p'}(\Om;\Mmn)$ and $L^{p}(\Om;\Mmn)$.

By (\ref{alphaW}) and (\ref{alpha2W})
there exist six constants $\alpha_0^\W>0$, $\alpha_1^\W>0$,
$\alpha_2^\W>0$, $\beta_0^\W\ge 0$, $\beta_1^\W\ge0$, $\beta_2^\W\ge0$
such that
\begin{eqnarray}
& \alpha_0^\W\|\Phi\|_{p}^{p}
-\beta_0^\W\le \W(\Phi)\le
\alpha_1^\W\|\Phi\|_{p}^{p}
+\beta_1^\W\,,\label{W9} \\
& |\langle \partial \W(\Phi), \Psi \rangle|\le
(\alpha_2^\W\|\Phi\|_{p}^{p-1}
+\beta_2^\W)\|\Psi\|_{p}^{} \,, \label{W10}
\end{eqnarray}
for every $\Phi,\, \Psi\in L^p(\Om;\Mmn)$.

If $u\in GSBV(\Om;\Rm)$ is an admissible
deformation for some crack $\Ga$ and $u$ has finite bulk
energy, then
$u$ belongs to $GSBV^p(\Om;\Rm)$ by (\ref{alphaW}) and its bulk energy
is given by
$\W(\nabla u)$.

\subsection{The body forces}\label{body}
We assume that at each time $t\in[0,T]$ the applied {\em body
forces\/} depend on the deformation $u$ and
are conservative. This means that there exists a function
$F\colon[0,T]{\times}\Om{\times}\Rm\to\R$ such that the density of
the applied body forces per unit volume in the reference configuration
is given by
$\partial_zF(t,x,u(x))$, where $\partial_zF(t,x,z)$ denotes the
partial gradient of $F$ with respect to $z$. We always assume that for every $t\in[0,T]$ the function $(x,z)\mapsto F(t,x,z)$ is ${\mathcal L}^n$-measurable in $x$ and $C^1$ in $z$.

Rather than imposing further regularity conditions on $F$, we prefer to impose appropriate conditions on the associated work, corresponding to the deformation $u$, given by
\begin{equation}\label{F7}
\F(t)(u):=\into F(t,x,u(x))\,dx\,. 
\end{equation}
This is because only $\F(t)$ enters in the expression for the energy. 
 
First of all we assume that there exists $q>1$ such that for every $t\in [0,T]$ the function $\F(t)$ is of class $C^1$ on $L^q(\Om;\Rm)$, with differential
$\partial \F(t)\colon L^{q{}}(\Om;\Rm)\to L^{q{}'}(\Om;\Rm)$ given by 
\begin{equation}\label{F7.5}
\langle\partial \F(t)(u),v\rangle=\int_\Om\partial_zF(t,x,u(x))\,v(x)\,dx
\end{equation}
for every $u,\, v\in L^{q{}}(\Om;\Rm)$, where $\langle\cdot,\cdot\rangle$ is the duality pairing between $L^{q'}(\Om;\Rm)$ and $L^{q}(\Om;\Rm)$, and $q':=q/(q-1)$. We assume also that
\begin{equation}\label{(2)}
\F(t)(u)\ge\limsup_{k\to\infty}\F(t)(u_k)
\end{equation}
whenever $u_k,\, u\in L^q(\Om;\Rm)$ and $u_k\to u$ a.e.\ on $\Om$. Notice that this inequality follows from  Fatou's Lemma and from the continuity hypothesis on $F$, 
provided suitable upper bounds are satisfied.
 
About the regularity in $t$ we assume that there exist a constant $\dot q\in{[1,q)}$ and, for a.e.\ $t\in[0,T]$, a function 
$\dot\F(t)\colon L^{\dot q}(\Om;\Rm)\to \R$ of class $C^1$, with differential 
$\partial \dot\F(t)\colon L^{\dot q}(\Om;\Rm)\to L^{\dot q'}(\Om;\Rm)$, $\dot q':=\dot q/(\dot q-1)$, such that
\begin{eqnarray}
&\displaystyle \F(t)(u)=\F(0)(u)+\int_0^t\dot\F(s)(u)\,ds\,, \label{F8}\\
&\displaystyle \langle \partial \F(t)(u),v\rangle = \langle \partial
\F(0)(u),v\rangle + \int_0^t \langle \partial
\dot\F(s)(u),v\rangle\, ds \label{F8.5}
\end{eqnarray}
for every $u,\, v\in L^{q}(\Om;\Rm)$ and for every $t\in[0,T]$. In order for (\ref{F8}) and (\ref{F8.5}) to make sense, we assume in addition that $t\mapsto\dot\F(t)(u)$ and 
$t\mapsto\langle \partial\dot\F(t)(u),v\rangle$ are integrable on $[0,T]$ for every $u,\,v\in L^q(\Om;\Rm)$. Of course, under these assumptions the functions
$t\mapsto \F(t)(u)$ and $t\mapsto\langle \partial \F(t)(u),v\rangle$ are absolutely continuous on $[0,T]$ for every $u,\,v\in L^q(\Om;\Rm)$, and their
time derivatives coincide a.e.\ with $t\mapsto\dot\F(t)(u)$ and $t\mapsto\langle \partial\dot\F(t)(u),v\rangle$, which justifies  the notation.
 
We also require some growth conditions on $\F(t)$, $\partial \F(t)$, $\dot\F(t)$, and $\partial \dot\F(t)$  to pass to the limit in our approximating sequences. To be explicit we assume that there exist six constants
$\alpha_0^\F>0$, $\alpha_1^\F>0$,
$\alpha_2^\F>0$, $\beta_0^\F\ge 0$, $\beta_1^\F\ge0$, $\beta_2^\F\ge0$
and four nonnegative functions
$\alpha_3^\F$, $\alpha_4^\F$,
$\beta_3^\F$, $\beta_4^\F\in L^1([0,T])$
such that
\begin{eqnarray}
& \alpha_0^\F\|u\|_{q{}}^{q{}}
-\beta_0^\F\le - \F(t)(u)\le
\alpha_1^\F\|u\|_{q{}}^{q{}}
+\beta_1^\F\,,\label{F9} \\
& |\langle \partial \F(t)(u), v\rangle|\le
(\alpha_2^\F\|u\|_{q{}}^{q{}-1}
+\beta_2^\F) \|v\|_{q{}}^{}\,, \label{F10} \\
& |\dot\F(t)(u)|\le
\alpha_3^\F(t)\|u\|_{\dot q}^{\dot q}
+\beta_3^\F(t)\,, \label{F11}\\
&|\langle \partial \dot\F(t)(u), v\rangle|\le
(\alpha_4^\F(t)\|u\|_{\dot q}^{\dot q-1}
+\beta_4^\F(t)) \|v\|_{\dot q}^{} \label{F11.5}
\end{eqnarray}
for a.e.\ $t\in[0,T]$ and for every $u,\, v\in L^{q{}}(\Om;\Rm)$. Note that by continuity (\ref{F9}) and (\ref{F10}) hold for every $t\in[0,T]$.

\begin{remark}\label{rem0}
The conditions $q{}>1$ and $\alpha_0^F>0$ in (\ref{F9}) play a
crucial role in our results. {}From a mathematical viewpoint, they
ensure that every sequence of deformations with bounded total energy
is bounded in $L^{q{}}(\Om;\Rm)$.
{}From a mechanical viewpoint, these conditions ensure that, even if
the cracks divide the body into several components, no part of the
body is sent to infinity by the applied forces. Unfortunately, they exclude the case of a constant body force, which corresponds to a potential $F$ which is linear with respect to~$z$.
\end{remark}

\begin{remark}\label{rem100}
Let us explain the roles of the different exponents $q$ and $\dot q$
in (\ref{F9})--(\ref{F11.5}). By (\ref{F9}) we obtain the exponent
$q{}$ in (\ref{coerc1}), so that, when we apply the $GSBV$
compactness theorem to a sequence $u_k$ of functions with bounded energy, we obtain a subsequence which converges pointwise a.e.\ on $\Om$ and is bounded in
$L^{q{}}(\Om;\Rm)$, but, in general, does not converge strongly in
$L^{q{}}(\Om;\Rm)$. In
some estimates we need to pass to the limit in
sequences like $\dot\F(t)(u_k)$, and this is made possible by
(\ref{F11}), since
$\dot q<q{}$ and, therefore, $u_k$
converges strongly in $L^{\dot q}(\Om;\Rm)$.
\end{remark}

 \begin{remark}\label{0109}
 All the conditions for $\F(t)$ and $\dot\F(t)$ listed above are satisfied whenever it is assumed that
 \begin{eqnarray}
&&\hbox{for every $(t,z)\in[0,T]{\times}\Rm$ the function 
 $x\mapsto F(t,x,z)$ is
integrable on
$\Om$,\ \ \ \ \ \ \ \ \ }\label{F1}\\
&&\hbox{for every $x\in\Om$ the function $(t,z)\mapsto
F(t,x,z)$
belongs to $C^2([0,T]{\times}\Rm)$,}\label{F2}
\end{eqnarray}
and that there exist seven
constants $q>\dot q\ge 1$, $a_0^F>0$, $a_1^F>0$, $a_2^F>0$, $a_3^F\ge0$, $a_4^F\ge0$ and
five nonnegative functions
$b_0^F,\, b_1^F\in C^0([0,T];L^1(\Om))$,
$b_2^F\in C^0([0,T];L^{q'}(\Om))$,
$b_3^F\in L^1(\Om)$, and $b_4^F\in L^{\dot q'}(\Om)$
such that
\begin{eqnarray}
&a_0^F |z|^{q{}} - b_0^F(t,x)\le - F(t,x,z) \le
a_1^F |z|^{q{}} + b_1^F(t,x)\,,\label{alphaF}
\\
& |\partial_z F(t,x,z)| \le a_2^F |z|^{q{}-1} + b_2^F(t,x)\,,
\label{alphaF2}\\
&|\partial_t F(t,x,z)|\le a_3^F|z|^{\dot q}+b_3^F(x)\,,\label{F6}\\
&|\partial_z \partial_t F(t,x,z)|\le a_4^F|z|^{\dot q-1}+b_4^F(x)
\label{F6.5}
\end{eqnarray}
for every $(t,x,z)\in [0,T]{\times}\Om{\times}\Rm$, where $\partial_t$
denotes the partial derivative with respect to~$t$. In this case (\ref{F8}) and (\ref{F8.5}) are satisfied with
\begin{equation}\label{F20}
\dot\F(t)(u):=\int_\Om \partial_t F(t,x,u(x))\,dx\,,
\end{equation}
since for every $u,\, v\in L^q(\Om)$ the functions $t\mapsto\dot\F(t)(u)$ and $t\mapsto\langle\partial\dot\F(t)(u),v\rangle$ are the continuous time derivatives of the functions $t\mapsto\F(t)(u)$ and $t\mapsto\langle\partial\F(t)(u),v\rangle$. Weaker hypotheses on $F$ will be considered in Section~\ref{ipgen}.
\end{remark}
\subsection{The surface forces}\label{surfacef}
We assume that at
each time $t\in[0,T]$ the {\em surface forces\/} applied on
$\partial_S\Om$ depend on the deformation $u$,
are conservative, and can be expressed by means of a potential
function $G\colon[0,T]{\times}\partial_S\Om{\times}\Rn\to\R$. More
precisely, we assume that for a given deformation $u$
the density of the applied surface forces
per unit area in the reference configuration is given by
$\partial_zG(t,x,u(x))$, where $\partial_zG(t,x,z)$ denotes the
partial gradient of $G$ with respect to $z$.
We assume also that for every $t\in[0,T]$ the function $(x,z)\mapsto G(t,x,z)$ is $\hn$-measurable in $x$ and $C^1$ in $z$.

\begin{remark}\label{surface} These assumptions are natural when the
surface forces applied to the deformed body depend on a conservative
field acting on a charge distribution which is deformed with the
body, i.e., the charge density per unit area in the reference
configuration does not depend on the deformation. Indeed, in this
case the change in the area elements between the reference and the
deformed configuration is compensated by the corresponding change in
the charge densities, so that the surface force applied to the
deformed body has a density per unit area in the reference
configuration which depends on the position $x$ and on the
deformation $u(x)$, but not on the deformation gradient $\nabla u(x)$.

Unfortunately,  pressure forces do not satisfy this assumption, so
they cannot be treated directly in the framework of this paper.
\end{remark}

As we did for body forces, we will impose appropriate conditions on the associated work, corresponding to the deformation $u$, which in this case is given by
\begin{equation}\label{G7}
\G(t)(u):=\int_{\partial_S\Om} G(t,x,u(x))\,d\hn(x)\,. 
\end{equation}

First of all we assume that there exists $r>1$ such that for every $t\in [0,T]$ the function $\G(t)$ is of class $C^1$ on $L^r(\partial_S\Om;\Rm)$, with differential
$\partial\G(t)\colon L^r(\partial_S\Om;\Rm)\to L^{r'}(\partial_S\Om;\Rm)$ given by 
\begin{equation}\label{G7.5}
\langle\partial \G(t)(u),v\rangle=\int_{\partial_S\Om}\partial_zG(t,x,u(x))\,v(x)\,d\hn(x)
\end{equation}
for every $u,\, v\in L^{r}(\partial_S\Om;\Rm)$, where
$\langle\cdot,\cdot\rangle$ denotes the duality pairing between
$L^{r'}(\partial_S\Om;\Rm)$ and $L^{r}(\partial_S\Om;\Rm)$, and $r':=r/(r-1)$. If the exponent $p$ which appears in (\ref{alphaW}) is less than the dimension
$n$ of $\Om$, we suppose that $p \le r< p(n-1)/(n-p)$. If $p\ge n$, we just suppose that $p \le r$.

Let us fix an open set
$\Om_S\subset \Om\setmeno \overline\Om_B$ with Lipschitz
boundary and such that $\partial_S\Om\subset \partial\Om_S$.
Under our hypothesis on $p$ and $r$ the trace operator from $W^{1,p}(\Om_S;\Rm)$ into
$L^{r}(\partial_S\Om;\Rm)$ is compact
(see, e.g., \cite{Nec}). Therefore there exists a constant $\gamma_S>0$ such
that
\begin{equation}\label{trace}
\|u\|_{r,\partial_S\Om}\le \gamma_S
(\|\nabla u\|_{p,\Om_S} +
\|u\|_{p,\Om_S})
\end{equation}
for every $u\in W^{1,p}(\Om_S;\Rm)$.
Here and in the rest of the paper we use
the same symbol to denote a function defined on (a set containing)
$\Om_S$ and its trace on $\partial_S\Om$. In particular,
if $u\in GSBV(\Om;\Rm)$ is a deformation with
$S(u)\subsethn \overline\Om_B$, $\W(\nabla u)<+\infty$, and
$\F(t)(u)<+\infty$ for some $t\in[0,T]$, then by (\ref{W9}) and (\ref{F9}) 
$u\in GSBV^p_{q{}}(\Om;\Rm):=GSBV^p(\Om;\Rm)\cap L^{q{}}(\Om;\Rm)$ and $S(u)\cap\Om_S\eqhn\emptyset$, so
that by Lemma~\ref{rem3} $u\in W^{1,p}(\Om_S;\Rm)\cap
L^{q{}}(\Om_S;\Rm)$. Therefore (the trace of) $u$ belongs to
$L^{r}(\partial_S\Om;\Rm)$
and $\G(t)(u)$ is well defined.

As for the regularity in $t$ we assume that for a.e.\ $t\in[0,T]$ there exists a function $\dot\G(t)\colon L^r(\partial_S\Om;\Rm)\to \R$ of class $C^1$, with differential 
$\partial \dot\G(t)\colon L^r(\partial_S\Om;\Rm)\to L^{r'}(\partial_S\Om;\Rm)$,
such that
\begin{eqnarray}
&\displaystyle \G(t)(u)=\G(0)(u)+\int_0^t\dot\G(s)(u)\,ds\,, \label{G8}\\
&\displaystyle \langle \partial \G(t)(u),v\rangle = \langle \partial
\G(0)(u),v\rangle + \int_0^t \langle \partial
\dot\G(s)(u),v\rangle\, ds \label{G8.5}
\end{eqnarray}
for every $u,\, v\in L^r(\partial_S\Om;\Rm)$ and for every $t\in[0,T]$. In order for (\ref{G8}) and (\ref{G8.5}) to make sense, we also assume that $t\mapsto\dot\G(t)(u)$ and $t\mapsto\langle \partial\dot\G(t)(u),v\rangle$ are integrable on $[0,T]$ for every  $u,\,v\in L^r(\partial_S\Om;\Rm)$. This implies that
the functions
$t\mapsto \G(t)(u)$ and $t\mapsto\langle \partial \G(t)(u),v\rangle$ are absolutely continuous on $[0,T]$ for every $u,\,v\in L^r(\partial_S\Om;\Rm)$.

As we did for the body forces, we require some growth conditions on $\G(t)$, $\partial \G(t)$, $\dot\G(t)$, and $\partial \dot\G(t)$. To be explicit we assume that
 there exist six nonnegative constants $\alpha_0^\G$, $\alpha_1^\G$,
$\alpha_2^\G$, $\beta_0^\G$, $\beta_1^\G$, $\beta_2^\G$
and four nonnegative functions
$\alpha_3^\G$, $\alpha_4^\G$, $\beta_3^\G$, $\beta_4^\G\in L^1([0,T])$
such that
\begin{eqnarray}
& -\alpha_0^\G\|u\|_{r,\partial_S\Om}
-\beta_0^\G\le - \G(t)(u)\le
\alpha_1^\G\|u\|_{r,\partial_S\Om}^{r}
+\beta_1^\G\,,\label{G9} \\
& |\langle \partial \G(t)(u), v\rangle|\le
(\alpha_2^\G\|u\|_{r,\partial_S\Om}^{r-1}
+\beta_2^\G)\|v\|_{r,\partial_S\Om}^{}\,, \label{G10} \\
& |\dot\G(t)(u)|\le
\alpha_3^\G(t)\|u\|_{r,\partial_S\Om}^{r}
+\beta_3^\G(t) \label{G11}\\
&|\langle \partial \dot\G(t)(u), v\rangle|\le
(\alpha_4^\G(t)\|u\|_{r,\partial_S\Om}^{r-1}
+\beta_4^\G(t)) \|v\|_{r,\partial_S\Om}^{} \label{G11.5}
\end{eqnarray}
for a.e.\ $t\in[0,T]$ and for every $u,\, v\in
L^{r}(\partial_S\Om;\Rm)$. By continuity (\ref{G9}) and (\ref{G10}) hold for every $t\in[0,T]$. Note that
in the first inequality of (\ref{G9}) the term $\|u\|_{r,\partial_S\Om}$ appears with
exponent~$1$. This is because in (\ref{coerc1}) we want a constant
$\alpha_0^\E>0$.

\begin{remark}\label{01091}
 All the conditions for $\G(t)$ and $\dot\G(t)$ listed above are satisfied whenever it is assumed that
\begin{eqnarray}
&&\hbox{\hspace{-40pt}for every $(t,z)\in[0,T]{\times}\Rm$ the function
$x\!\mapsto \!G(t,x,z)$ is $\hn\!$-integrable on $\partial_S\Om$,}\label{G1}\\
&&\hbox{\hspace{-40pt}for every
$x\in\partial_S\Om$ the function $(t,z)\mapsto
G(t,x,z)$
belongs to $C^2([0,T]{\times}\Rm)$,}\label{G2}
\end{eqnarray}
and that there exist four constants $a_1^G\ge 0$, $a_2^G\ge0$, $a_3^G\ge0$, $a_4^G\ge0$ and
six nonnegative functions $a_0^G$, $b_2^G\in
C^0([0,T];L^{r'}(\partial_S\Om))$, 
$b_0^G,\, b_1^G\in C^0([0,T];L^1(\partial_S\Om))$,
 $b_3^G\in L^1(\partial_S\Om)$, and
$b_4^G\in  L^{r'}(\partial_S\Om)$
such that
\begin{eqnarray}
& -a_0^G(t,x) |z| - b_0^G(t,x)\le - G(t,x,z) \le
a_1^G |z|^{r} + b_1^G(t,x)\,,\label{alphaG}
\\
& |\partial_z G(t,x,z)| \le a_2^G |z|^{r-1} +
b_2^G(t,x)\,, \label{alphaG2} \\
&|\partial_t G(t,x,z)|\le
a_3^G|z|^{r}+b_3^G(x)\,,\label{G6}\\
&|\partial_z \partial_t G(t,x,z)|\le a_4^G|z|^{r-1}+b_4^G(x) \label{G6.5}
\end{eqnarray}
for every $(t,x,z)\in [0,T]{\times}\partial_S\Om{\times}\Rm$, where $\partial_t$
denotes the partial derivative with respect to~$t$. 
As in Remark~\ref{0109} we can prove that (\ref{G8}) and (\ref{G8.5}) are satisfied with
\begin{equation}\label{G20}
\dot\G(t)(u):=\int_{\partial_S\Om} \partial_t G(t,x,u(x))\,d\hn(x)\,.
\end{equation}
Weaker hypotheses on $G$ will be considered in Section~\ref{ipgen}.
\end{remark}

\begin{remark}\label{boundforces}
In this model $\overline\Om_B$ represents the reference configuration
of the brittle part of the material, while
$\Om\setmeno\overline\Om_B$ can be considered as the reference
configuration of an elastic unbreakable part attached to it through
the interface $\Om\cap\partial\Om_B$. Since $W$ is not assumed to be
continuous with respect to $x$, it may happen that the bulk energy
density is discontinuous across $\Om\cap\partial\Om_B$, so that we
can interpret $\overline\Om_B$ and $\Om\setmeno\overline\Om_B$ as
representing two bodies with different material properties. In other
words, in this model the surface forces can act on the brittle body $\Om_B$ only
through the layer of unbreakable material
$\Om\setmeno\overline\Om_B$. At the moment it is not known what
happens if the thickness of this layer tends to $0$.
\end{remark}

\subsection{The prescribed boundary deformations}\label{bounddef}
In this paper we do not
consider the case of imposed ``discontinuous" 
boundary deformations, but only boundary deformations
that are traces on $\partial_D\Om$ of functions $\psi\in
W^{1,p}(\Om;\Rm)\cap L^{q{}}(\Om;\Rm)$, so that there is always a
configuration with finite energy without cracks which satisfies the
boundary conditions. The choice of the exponents is determined by
(\ref{W9}) and (\ref{F9}).

The set $AD(\psi,\Ga)$ of {\em admissible deformations
in $\Om$ with finite energy\/}, corresponding to a crack
$\Ga\in \Rp(\overline\Om_B)$ and to a
boundary deformation $\psi\in
W^{1,p}(\Om;\Rm)\cap L^q(\Om;\Rm)$, is defined by
\begin{equation}\label{Vpk}
AD(\psi,\Ga):=\{u\in GSBV^p_{q{}}(\Om;\Rm):
S(u)\subsethn\Ga,\ u=\psi\ \,\, \hn\hbox{-a.e.\ on }\,
\partial_D\Om\setmeno \Ga\}\,,
\end{equation}
where the last equality in the previous formula refers to the traces
of $u$ and $\psi$ on $\partial\Om$ introduced in
Proposition~\ref{traces}.

Note that if $\Ga$ is closed, then $AD(\psi,\Ga)$ coincides with the space of all functions $u\in L^q(\Om;\Rm)$ whose distributional gradient on $\Om\setmeno\Ga$ belongs to $L^p(\Om\setmeno\Ga;\Mmn)$ and which agree with $\psi$ on $\partial_D\Om\setmeno\Ga$ in the standard sense of Sobolev spaces (to prove this fact we can use Lemma~\ref{rem3} and \cite[Lemma 2.3]{DG-C-L}). This space is frequently used in the variational approach to nonlinear elasticity. Our ``non-conventional" definition of $AD(\psi,\Ga)$ stems from the potential
failure of the crack $\Ga$ to remain closed  in our existence theorem.

We assume that the boundary deformation $\psi(t)$ depends on time and
that the function $t\mapsto\psi(t)$ is absolutely continuous from
$[0,T]$ into $W^{1,p}(\Om;\Rm)\cap L^{q{}}(\Om;\Rm)$ (endowed with
the sum of the norms), so that the time derivative
$t\mapsto\dot\psi(t)$ belongs to the space $L^1([0,T];
W^{1,p}(\Om;\Rm)\cap L^{q{}}(\Om;\Rm))$ and its spatial gradient
$t\mapsto\nabla \dot \psi(t)$ belongs to the space
$L^1([0,T];L^p(\Om;\Mmn))$.

\subsection{The admissible configurations and their total energy}
An {\em admissible configuration\/} is a pair $(u,\Ga)$,
where $\Ga\in \Rp(\overline\Om_B)$ is an
admissible crack and $u\in GSBV^p_{q{}}(\Om;\Rm)$ is a
deformation with finite energy and with
$S(u)\subsethn\Ga$.

For every
$t\in[0,T]$, $\Ga\in \Rp(\overline\Om_B)$, 
and $u\in AD(\psi(t),\Ga)$,
the {\em total energy\/} of the admissible configuration
$(u,\Ga)$
at time $t$ is
given by
\begin{equation}\label{E}
\E(t)(u,\Ga):=\Ec(t)(u) +\K(\Ga)\,,
\end{equation}
where for every $u\in GSBV^p_{q{}}(\Om;\Rm)$ the {\em elastic energy\/}
is defined by
\begin{equation}\label{E0}
\Ec(t)(u):=\W(\nabla u)-\F(t)(u)-\G(t)(u)\,.
\end{equation}
Note that $u\in
W^{1,p}(\Om_S;\Rm)$ by Lemma~\ref{rem3}, so that $u\in
L^{r}(\partial_S\Om,\Rm)$ and $\G(t)(u)$ is
well defined. We will sometimes write
\begin{equation}\label{Etot}
\E(t)(u,\Ga):=\E^{in}(u,\Ga)-{\mathcal L}(t)(u)
\,,
\end{equation}
where
\begin{equation}\label{Eint}
\E^{in}(u,\Ga):=\W(\nabla u) +\K(\Ga)
\end{equation}
is the {\em internal energy\/}, while
\begin{equation}\label{Loads}
{\mathcal L}(t)(u):=\F(t)(u)+\G(t)(u)
\end{equation}
is the {\em work done by the applied loads\/}.

There exist four constants
$\alpha_0^\E>0$, $\alpha_1^\E>0$, $\beta_0^\E\ge0$,
and $\beta_1^\E\ge0$ such that
\begin{eqnarray}
&\Ec(t)(u)\ge \alpha_0^\E(\|\nabla u\|_p^p +
\|u\|_{q{}}^{q{}})
- \beta_0^\E\,, \label{coerc1}\\
&\Ec(t)(u)\le
\alpha_1^\E(\|\nabla u\|_p^{p}+ \|u\|_{q{}}^{q}+\|u\|^r_{r,\partial_S\Om})
+\beta_1^\E \label{bound1}
\end{eqnarray}
for every $t\in[0,T]$ and for every $u\in GSBV^p_{q{}}(\Om;\Rm)$.

To prove this fact let us fix $t$ and $u$. By Lemma~\ref{rem3} we have
$u\in W^{1,p}(\Om_S,\Rm)$, and by (\ref{W9}), (\ref{F9}), and
(\ref{G9})
we have
\begin{equation}\label{sett1}
\Ec(t)(u)\ge \alpha_0^\W \|\nabla
u\|_p^p - \beta_0^\W +
\alpha_0^\F \|u\|_{q{}}^{q{}} -
\beta_0^\F
-\alpha_0^\G \|u\|_{r,\partial_S\Om}-
\beta_0^\G\,.
\end{equation}
Since $\Om_S$ has a Lipschitz boundary, there exists a constant
$k_S>0$, depending only on $p$, $q{}$, and $\Om_S$, such that
\begin{equation}\label{pq}
\|u\|_{p,\Om_S}\le k_S
(\|\nabla u\|_{p,\Om_S} +
\|u\|_{q{},\Om_S})
\end{equation}
for every
$u\in W^{1,p}(\Om_S;\Rm)\cap L^{q{}}(\Om_S;\Rm)$.

{}Using  Young's inequality,  it follows from (\ref{trace}) and (\ref{pq}) that 
there exists a 
constant $\lambda\ge 0$, depending only on $p$,
$r$, $q{}$, $\alpha_0^\W$, $\alpha_0^\F$, $\alpha_0^\G$, and $\Om_S$, such that
\begin{equation}\label{tracepq}
\alpha_0^\G
\|u\|_{r,\partial_S\Om}\le
\frac{\alpha_0^\W}{2} \|\nabla u\|_{p,\Om_S}^p +
\frac{\alpha_0^\F}{2}\|u\|_{q,\Om_S}^q+ \lambda
\end{equation}
for every
$u\in W^{1,p}(\Om_S;\Rm)\cap L^{q{}}(\Om_S;\Rm)$.
Therefore (\ref{coerc1}) follows from (\ref{sett1}), with $\alpha_0^{\E}:=\min\{\frac{1}{2}\alpha_0^\W,\frac{1}{2}\alpha_0^\F\}$, and $\beta_0^{\E}:=\beta_0^\W +\beta_0^\F+\beta_0^\G+\lambda$. 
 
By (\ref{W9}), (\ref{F9}), and
(\ref{G9}) we also have 
$$
\Ec(t)(u)\le \alpha_1^\W \|\nabla
u\|_p^p + \beta_1^\W +
\alpha_1^\F \|u\|_{q{}}^{q{}} +
\beta_1^\F
+\alpha_1^\G \|u\|_{r,\partial_S\Om}^{r}+
\beta_1^\G\,,
$$
which gives (\ref{bound1}) with $\alpha_1^\E:=\max\{\alpha_1^\W,\alpha_1^\F,\alpha_1^\G\}$ and $\beta_1^\E:=\beta_1^\W+\beta_1^\F+\beta_1^\G$.

\subsection{Minimum energy configurations}\label{minen}
For a given time $t\in[0,T]$ and a given crack
$\Ga(t)\in\Rp(\overline\Om_B)$, a deformation
$u$ which corresponds to static equilibrium is a
critical
point of the functional $\Ec(t)$ on the set
$AD(\psi(t),\Ga(t))$ defined in (\ref{Vpk}).
Among such critical points, the minimum points of the problem
\begin{equation}\label{min00}
\min_{u\in AD(\psi(t),\Ga(t))} \Ec(t)(u)
\end{equation}
play an important role and can be regarded as the most stable equilibria. We will call them
{\em minimum energy deformations\/} at time $t$ with crack
$\Ga(t)$.
The existence of these minimizers is guaranteed by the following
theorem, that will be proved in Section~\ref{preliminary}.
\begin {theorem}\label{esistmin00}
For every $t\in[0,T]$ and every
$\Ga(t)\in\Rp(\overline\Om_B)$
the minimum problem (\ref{min00})
has a solution.
\end{theorem}
Let $u(t)$ be a minimum energy deformation at time $t$ with crack
$\Ga(t)\in\Rp(\overline\Om_B)$. For every $v\in AD(0,\Ga(t))$ and every $\e\in\R$ the 
function $u(t)+\e v$ belongs to $AD(\psi(t),\Ga(t))$. Therefore
$\Ec(t)(u(t))\le \Ec(t)(u(t)+\e v)$ for every $\e\in\R$. By taking 
the derivative with respect to $\e$ at $\e=0$, 
 we obtain the weak formulation of the 
{\em Euler equation\/}
\begin{equation}\label{Euler}
\langle \partial\W(\nabla u(t)),\nabla v\rangle= \langle 
\partial\F(t)(u(t)),v\rangle+ \langle \partial\G(t)(u(t)),v\rangle
\end{equation}
for every $v\in AD(0,\Ga(t))$.
The critical
points of the functional $\Ec(t)$ on
$AD(\psi(t),\Ga(t))$ are, by definition, the solutions $u(t)\in 
AD(\psi(t),\Ga(t))$ of (\ref{Euler}), which turns out to be the 
equation of equilibrium with prescribed crack $\Ga(t)$, and coincides with the classical equilibrium equation considered in nonlinear elasticity when $\Ga(t)$ is closed.

If $\Ga(t)\in\Rp(\overline\Om_B)$ and $u(t)\in AD(\psi(t),\Ga(t))$ is a solution of (\ref{Euler}), we introduce
the linear functional  $g(t)$ on $GSBV^p_q(\Om;\Rm)$ defined by
\begin{equation}\label{g(t)}
\langle g(t),v\rangle:=\langle \partial\W(\nabla u(t)),\nabla 
v\rangle-\langle \partial\F(t)(u(t)),v\rangle- \langle 
\partial\G(t)(u(t)),v\rangle
\end{equation}
for every $v\in GSBV^p_q(\Om;\Rm)$.
By the Euler equation we have $\langle g(t),v_1\rangle=\langle 
g(t),v_2\rangle$ for every $v_1,\, v_2\in GSBV^p_q(\Om;\Rm)$ with $S(v_1)\subsethn \Ga(t)$, $S(v_2)\subsethn \Ga(t)$, and 
$v_1=v_2$ $\hn$-a.e.\ on $\partial_D\Om\setmeno\Ga(t)$,
since in this case $v_1-v_2\in AD(0,\Ga(t))$. In other words
$\langle g(t),v\rangle$ depends only on the trace of $v$ on 
${\partial_D\Om\setmeno\Ga(t)}$, provided $S(v)\subsethn\Ga(t)$.

Under suitable regularity assumptions we have
\begin{equation}\label{sforce}
\langle g(t),v\rangle=\int_{\partial_D\Om\setminus\Ga(t)} \partial_\xi 
W(\nabla u(t))\nu\, v\, d\hn\,,
\end{equation}
where $\nu$ is the outer unit normal to $\partial\Om$, so that $g(t)$ 
can be identified with the function $ \partial_\xi W(\nabla u(t))\nu$ 
defined on $\partial_D\Om\setmeno\Ga(t)$, which represents the 
density per unit area of the surface force acting on 
$\partial_D\Om\setmeno\Ga(t)$ at time $t$.

Returning to the general case considered at the beginning, the expression
$\langle g(t),v\rangle$ can always be interpreted as the work done by 
the surface forces acting on $\partial_D\Om\setmeno\Ga(t)$ at time 
$t$ under the deformation $v$.

In the spirit of Griffith's theory, an {\em
equilibrium configuration\/} at time $t\in[0,T]$ is an admissible
configuration $(u(t),\Ga(t))$ which is a
``critical point", in a sense that has not yet been made mathematically
precise, of the functional
$\E(t)(u,\Ga)$ on the set of configurations
$(u,\Ga)$ with $\Ga\in\Rp(\overline\Om_B)$, 
$\Ga(t)\subsethn\Ga$, and
$u\in AD(\psi(t),\Ga)$.
Following \cite{F-M}, we will consider only {\em minimum energy
configurations at time\/} $t$, which are defined as the
admissible configurations $(u(t),\Ga(t))$, with 
$\Ga(t)\in\Rp(\overline\Om_B)$ and
$u(t)\in AD(\psi(t),\Ga(t))$, such that
$$
\E(t)(u(t),\Ga(t))\le \E(t)(u,\Ga)
$$
for every
$\Ga\in\Rp(\overline\Om_B)$, with 
$\Ga(t)\subsethn\Ga$, and every
$ u\in AD(\psi(t),\Ga)$.
These are regarded as the most stable equilibrium configurations.

The following theorem, which will be proved in
Section~\ref{preliminary}, ensures that for every $t\in[0,T]$ and for
every $\Ga_0\in\Rp(\overline\Om_B)$ there
exists at least a minimum energy configuration
$(u(t),\Ga(t))$ such that
$\Ga_0\subsethn\Ga(t)$.

\begin {theorem}\label{esistmin}
Let $t\in[0,T]$ and let
$\Ga_0\in\Rp(\overline\Om_B)$.
Then the minimum problem
\begin{equation}\label{min}
\min\,\{ \E(t)(u,\Ga): \Ga\in\Rp(\overline\Om_B),\ \Ga_0\subsethn\Ga,\ u\in AD(\psi(t),\Ga)\,\}
\end{equation}
has a solution.
\end{theorem}

\subsection{Quasistatic evolution}\label{quasistatic}
An {\em irreversible quasistatic evolution of minimum energy 
configurations\/} is a function $t\mapsto(u(t),\Ga(t))$
which satisfies the following conditions:
\begin{itemize}
\smallskip
\item[{\rm (a)}] {\em static equilibrium}: for every $t\in [0,T]$ 
the pair $(u(t),\Ga(t))$
is a minimum energy configuration at time $t$, i.e.,
$ \Ga(t)\in \Rp(\overline\Om_B)$, $u(t)\in 
AD(\psi(t),\Ga(t))$, and
$$
\E(t)(u(t) ,\Ga(t))\le \E(t )(v,\Ga)
$$
for every
$ \Ga\in \Rp(\overline\Om_B)$, with $\Ga(t) 
\subsethn\Ga$, and every $v\in AD(\psi(t),\Ga)$;
\smallskip 
\item[{\rm (b)}] {\em irreversibility\/}:
$\Ga(s)\subsethn\Ga(t)$ for $0\le s<t\le T$;
\smallskip
\item[{\rm (c)}]
{\em nondissipativity\/}:
the function $t\mapsto E(t):=\E(t)(u(t),\Ga(t))$
is absolutely continuous on $[0,T]$ and its time derivative
$\dot E(t)$ satisfies
\begin{equation}\label{dote}
\dot E(t)= \langle g(t), \dot\psi(t) \rangle
- \dot\F(t)(u(t))- \dot\G(t)(u(t))
\end{equation}
for a.e.\ $t\in[0,T]$.
\smallskip
\end{itemize}

\begin{remark}\label{conservation}
To explain why condition (c) can be interpreted
as the conservation of energy in this model, let us consider the very special case
where $\Ga(t)\eqhn\Ga_0$ for every $t\in[0,T]$,  with $\Ga_0$ a closed set, and the function $t\mapsto u(t)$ is absolutely continuous from $[0,T]$ into $W^{1,p}(\Om\setminus\Ga_0;\Rm)\cap L^q(\Om;\Rm)$ (endowed with the sum of the norms). Then  the time derivative $t\mapsto \dot u(t)$ belongs to the space $L^1([0,T];W^{1,p}(\Om\setmeno\Ga_0;\Rm)\cap L^q(\Om;\Rm))$. Since $u(t)-\psi(t)=0$ $\hn$-a.e.\ on $\partial_D\Om\setmeno\Ga_0$ for every $t\in[0,T]$, we obtain that $\dot u(t)-\dot\psi(t)=0$ $\hn$-a.e.\ on $\partial_D\Om\setmeno\Ga_0$
for a.e.\ $t\in[0,T]$, so that $\dot u(t)-\dot\psi(t)\in AD(0,\Ga_0)$ for a.e.\ $t\in[0,T]$.
{}From the Euler equation (\ref{Euler}) and from (\ref{g(t)}) we obtain 
$\langle g(t),\dot u(t)\rangle=\langle g(t),\dot\psi(t)\rangle$, 
so that (\ref{dote}) yields
\begin{equation}\label{dote2}
\dot E(t)= \langle g(t), \dot u(t) \rangle- \dot\F(t)(u(t))- 
\dot\G(t)(u(t))
\end{equation}
for a.e.\ $t\in[0,T]$. On the other hand, we have 
\begin{eqnarray}
& \frac{d}{dt}\F(t)(u(t))=\langle \partial\F(t)(u(t)), \dot u(t) \rangle+\dot\F(t)(u(t))\,,
\label{dtf}\\
& \frac{d}{dt}\G(t)(u(t))=\langle \partial\G(t)(u(t)), \dot u(t) \rangle+\dot\G(t)(u(t))\label{dtg}
\end{eqnarray}
for a.e.\ $t\in[0,T]$. This follows from our qualitative hypotheses on $\F$, $\dot\F$, $\G$, $\dot\G$, together with the estimates given by (\ref{F8}), (\ref{F8.5}), (\ref{F11}), (\ref{F11.5}), (\ref{G8}), (\ref{G8.5}), (\ref{G11}), and~(\ref{G11.5}).

Let
$$
E^{in}(t):=\E^{in}(u(t),\Ga_0)=\E(t)(u(t),\Ga_0)+\F(t)(u(t))+\G(t)(u(t))
$$
be the interior energy of the solution at time~$t$.
By (\ref{dote2}), (\ref{dtf}), and (\ref{dtg}) we have
$$
\dot E^{in}(t)= \langle g(t), \dot u(t)\rangle +
\langle \partial\F(t)(u(t)), \dot u(t) \rangle+
\langle \partial\G(t)(u(t)), \dot u(t) \rangle\,,
$$
where the right hand side is the power of the exterior forces applied
to the body at time~$t$, including the surface forces acting on 
${\partial_D\Om\setmeno\Ga_0}$ (see (\ref{sforce})). Similar results can be obtained under weaker regularity conditions on $u(t)$ and~$\Ga(t)$.
\end{remark}

\begin{remark}\label{partialN}
 If $t\mapsto(u(t),\Ga(t))$ is an irreversible quasistatic evolution of minimum energy configurations, so is $t\mapsto(u(t),\Ga(t)\setmeno\partial_N\Om)$. 
 \end{remark}
 
\begin{remark}\label{meas}
 In the definition of quasistatic evolution we make no measurability assumption on the function $t\mapsto u(t)$. However, the nondissipativity condition (c) implies that the function 
 $t\mapsto\langle g(t), \dot\psi(t) \rangle- \dot\F(t)(u(t))- \dot\G(t)(u(t))$ 
is measurable and belongs to $L^1([0,T])$ being the a.e.\ derivative of an absolutely continuous function.
\end{remark}

Given an initial admissible configuration $(u_0,\Ga_0)$, we
look for an irreversible quasistatic evolution such that
$(u(0),\Ga(0))=(u_0,\Ga_0)$.
{}From the definition it
follows that
$(u(0),\Ga(0))$ is a
minimum energy configuration at time~$0$. Therefore
a necessary condition for the solvability of the initial value problem
is that
$(u_0,\Ga_0)$ is a
minimum energy configuration at time~$0$, i.e., $\Ga_0\in 
\Rp(\overline\Om_B)$,
$u_0\in AD(\psi(0),\Ga_0)$, and
\begin{equation}\label{mininit}
\E(0)(u_0,\Ga_0)\le \E(0)(u,\Ga)
\end{equation}
for every $ \Ga\in \Rp(\overline\Om_B)$, with
$\Ga_0\subsethn\Ga$, and every $u\in AD(\psi(0),\Ga)$.

The following remark shows that, if no
forces are applied at time~$0$, then there are minimum energy 
configurations at time~$0$ with an arbitrary crack
$\Ga_0$, provided that some very weak
additional conditions are satisfied.

\begin{remark}\label{identical} Let $m=n$ and 
let $u_{id}(x):=x$ be the identical deformation, so that $\nabla
u_{id}(x)$ is the identity matrix in $\Mnn$. If $\psi(0)=u_{id}$ and
\begin{eqnarray*}
&W(x,\xi)\ge W(x,\nabla u_{id}(x))\,,\\
& F(0,x,z)\le 0=F(0,x,u_{id}(x))\,,\\
&G(0,x,z)\le 0=G(0,x,u_{id}(x))\hphantom{\,,}
\end{eqnarray*}
for every $x\in\Om$, $z\in\Rn$, $\xi\in\Mnn$, then no force is
applied to the body at time $0$ and the reference configuration
$u_{id}$ belongs to $AD(\psi(0),\emptyset)$ and is stress free.
Moreover for every admissible crack
$\Ga\in \Rp(\overline\Om_B)$ the 
configuration
$(u_{id},\Ga)$ is a minimum energy configuration
at time~$0$.
\end{remark}

The main result of this paper is the following theorem, which will be
proved in Section~\ref{proof}.

\begin{theorem}\label{main}
Let $(u_0,\Ga_0)$ be a minimum energy
configuration at time~$0$, i.e., assume 
$\Ga_0\in \Rp(\overline\Om_B)$,
$u_0\in AD(\psi(0),\Ga_0)$, and
(\ref{mininit}). Then there exists an irreversible quasistatic
evolution
$t\mapsto(u(t),\Ga(t))$ with
$(u(0),\Ga(0))=(u_0,\Ga_0)$.
\end{theorem}
\end{section}

\begin{section}{A FEW TOOLS}\label{tools}

In this section we develop a few tools which will be very useful 
in the proof of
Theorem~\ref{main}. Subsections~\ref{sets} and~\ref{compa} introduce a weak notion of set convergence, which plays the role of Hausdorff convergence when no restriction is placed on either dimensionality or connectedness. Hausdorff convergence, which was instrumental in \cite{DM-T} and \cite{Ch}, can not be used in the present setting for two reasons: first, $\hn$ is not lower semicontinuous with respect to Hausdorff convergence; then, Hausdorff convergence does not imply convergence of the associated minimum energy deformations.
 
 Subsection~\ref{mth} establishes the weak convergence of the stresses associated to the minimum energy deformations, while Subsection~\ref{aRiem} deals with a technical result concerning the approximation of Bochner integrals with Riemann sums.

As in Section~\ref{spaces},
let $U$ be a bounded open set in $\Rn$ and let $1<p<+\infty$.

\subsection{A convergence of sets}\label{sets}
We introduce a notion of convergence of sets based on the weak 
convergence in $SBV^p(U)$.

\begin{definition}\label{setsbv}
We say that $\Ga_k$ $\sigmap$-converges to $\Ga$ in $U$ if $\Ga_k,\, 
\Ga\subset U$, $\hn(\Ga_k)$ is bounded uniformly with respect to $k$,  and the following conditions are satisfied:
\begin{itemize}
\item[{\rm(a)}] if $u_j$ converges weakly to $u$ in $SBV^p(U)$ and 
$S(u_j)\subsethn \Ga_{k_j}$ for some sequence $k_j\to\infty$, then 
$S(u)\subsethn\Ga$;
\item[{\rm(b)}] there exist a function $u\in SBV^p(U)$ and a 
sequence $u_k$ converging to $u$ weakly in $SBV^p(U)$ such that 
$S(u)\eqhn\Ga$ and $S(u_k)\subsethn \Ga_k$ for every $k$.
\end{itemize}
\end{definition}

\begin{remark}\label{sigmap}
The rectifiability of the $\sigmap$-limit of any sequence of sets follows 
from the rectifiability of $S(u)$ for any $u\in SBV^p(U)$.

It is clear from the definition that, if a sequence $\sigmap$-converges, then every subsequence $\sigmap$-converges to the same limit. If $\Ga_k$ and $\Ga'_k$ $\sigmap$-converge to $\Ga$ and $\Ga'$, respectively, then
$$
\Ga_k\subsethn \Ga'_k\ \,\hbox{for every}\ \, k
\quad\Longrightarrow\quad
\Ga\subsethn \Ga'\,.
$$

Let us consider now the special case where $\Ga_k=\Ga_0$ for every $k$, with $\hn(\Ga_0)<+\infty$.
Then it is not always true that
$\Ga_k$ $\sigmap$-converges to $\Ga_0$ in $U$. Indeed, using Theorem~\ref{compattezza} and Remark~\ref{Su} we can prove that in this case $\Ga_k$ $\sigmap$-converges to the set $\Ga$ characterized by the following properties:
\begin{itemize}
\item[{\rm(a)}] if $v\in SBV^p(U)$ and 
$S(v)\subsethn \Ga_0$, then 
$S(v)\subsethn\Ga$;
\item[{\rm(b)}] there exists a function $u\in SBV^p(U)$ such that 
$S(u)\eqhn\Ga \subsethn \Ga_0$.
\end{itemize}
Therefore $\Ga_k$ $\sigmap$-converges to $\Ga_0$ in $U$ if and only if 
\begin{equation}\label{exists}
\hbox{there exists }\, u\in SBV^p(U) \, \hbox{ with }\, S(u)\eqhn \Ga_0 \,.
\end{equation}
Note that (\ref{exists}) is not always true. For instance, if $\Ga_0$ is contained in a smooth manifold $M$ of dimension $n-1$, then (\ref{exists}) implies that there exists $u\in W^{1,p}({\Om\setmeno M})$ such that $\Ga_0\eqhn \{{x\in M}: u^-(x)\neq u^+(x)\}$, where $u^-$ and $u^+$ are the traces of $u$ on both sides of $M$. If $p>n$, by the Sobolev embedding theorem $u^-$ and $u^+$ are continuous on $M$, hence a set $\Ga_0\subset M$ which satisfies (\ref{exists}) must be open.
Using the notion of capacity associated to $W^{1,p}$ it is possible to prove that also in the case $1<p\le n$ there are sets $\Ga_0\subset M$ which do not satisfy~(\ref{exists}).
\end{remark}

The following lower semicontinuity theorem is an easy consequence of 
Theorem~\ref{semigsbv}.
\begin{theorem}\label{semigamma}
Let $\kappa$ be as in Theorem~\ref{semigsbv}, let
$\Ga_k$, $\Ga$, and $\Ga^\prime$ be rectifiable subsets of $U$ with $\hn(\Ga^\prime)<+\infty$, and let $E$ be an $\hn$-measurable set with $\hn(E)<+\infty$. If $\Ga_k$ $\sigmap$-converges to $\Ga$ in $U$, then
$$
\int_{(\Ga\cup\Ga^\prime)\setminus E}\kappa(x,\nu)\,d\hn\le
\liminf_{k\to\infty}\int_{(\Ga_k\cup\Ga^\prime)\setminus E} \kappa(x,\nu_k)\,d\hn\,,
$$
where $\nu$ and $\nu_k$ are unit normal vector fields on
$\Ga\cup\Ga^\prime$ and $\Ga_k\cup\Ga^\prime$, respectively.
\end{theorem}
\begin{proof}
By condition (b) there exist a function $u\in SBV^p(U)$ and a 
sequence $u_k$ converging to $u$ weakly in $SBV^p(U)$ such that 
$S(u)\eqhn\Ga$ and $S(u_k)\subsethn \Ga_k$ for every $k$. The 
conclusion follows by applying Theorem~\ref{semigsbv} and Remark~\ref{union}.
\end{proof}

\begin{remark}\label{bgk}
Assume that $\Ga_k$ $\sigmap$-converges to $\Ga$, then 
\begin{itemize}
\item[{\rm(a)}]
$\hn(\Ga)<+\infty$ by Theorem~\ref{semigamma};
\item[{\rm(b)}]
if further, $K$ is compact and $\Ga_k \subsethn K$ for every $k$, then $\Ga  \subsethn K$
by (b) of Definition~\ref{setsbv}, together with the fact that
$\hn(S(u)\setminus K)\leq \liminf_k\hn(S(u_k)\setminus K)$.
\end{itemize}
\end{remark}

\begin{lemma}\label{cond2}
Let $\Ga\subset U$ be a set with $\hn(\Ga)<+\infty$. If there exists 
a sequence $u^i\in SBV^p(U)\cap L^\infty(U)$ such that 
$\Ga\eqhn\bigcup_i S(u^i)$, then there exists $u\in SBV^p(U)\cap 
L^\infty(U)$ such that $S(u)\eqhn\Ga$. 
Let $\Ga_k$ be a 
sequence of subsets of $U$ 
with $\hn(\Ga_k)$ bounded uniformly with respect to $k$. If, in addition to the previous hypotheses, each 
 function $u^i$ is the weak 
limit in $SBV^p(U)$ of a sequence $u_k^i$ with 
$S(u_k^i)\subsethn\Ga_k$ for every~$k$, then
condition (b) of Definition~\ref{setsbv} is satisfied.
\end{lemma}

\begin{proof}
Suppose that there exists a sequence $u^i\in SBV^p(U)\cap 
L^\infty(U)$ such that $\Ga\eqhn\bigcup_i S(u^i)$. It is not 
restrictive to
assume that
\begin{equation}\label{limuik}
\|u^i\|_\infty\le 1\,,\qquad
\|\nabla u^i\|_p\le 1\,.
\end{equation}
Given a sequence of real numbers $c_i>0$, with 
$\sum_ic_i<+\infty$, we define
$$
u:=\sum_{i=1}^\infty c_iu^i\,,\qquad v^\ell:=\sum_{i=1}^\ell c_iu^i\,.
$$
Since $v^\ell\in SBV^p(U)\cap L^\infty(U)$ and $S(v^\ell)\subsethn \Ga$,
from Theorem~\ref{compsbv} we obtain that $u\in SBV^p(U)\cap 
L^\infty(U)$ and $v^\ell$ converges to $u$ weakly in $SBV^p(U)$. By 
Remark~\ref{Su} we have
$S(u)\subsethn\Ga$.

It remains to choose the sequence $c_i$ so that $\Ga\subsethn S(u)$.
First of all we fix a (Borel) orientation of $\Ga$ and for every 
$v\in SBV^p(U)$ we define the jump of $v$ on $\Ga$ as 
$[v]:=v^+-v^-\in L^1(\Ga,\hn)$.
We define two sequences $c_i$ and $\e_i$ inductively. We set $c_1=1$.
Suppose that $c_\ell$ and $\e_{\ell-1}$ have already been defined. We 
choose $\e_\ell$ such that $0<\e_\ell<\e_{\ell-1}$ (for $\ell=1$ we 
require only $0<\e_1$) and
\begin{equation}\label{em}
\hn(\{x\in S(v^\ell):|[v^\ell](x)|<\e_\ell\})<2^{-\ell}\,.
\end{equation}
To choose $c_{\ell+1}$, for every $c>0$ we consider the set
$$
A^{\ell}_c:=\{x\in S(v^\ell): [v^\ell](x)+c[u^{\ell+1}](x)=0\}\,,
$$
as in the proof of \cite[Lemma~3.1]{Fra-Lar}.
The family $(A^{\ell}_c)^{}_{c>0}$ is composed of  pairwise disjoint subsets of $S(v^\ell)$. As $\hn(S(v^\ell))<+\infty$, we 
have that $\hn(A^{\ell}_c)=0$ except for a countable number of $c$. We 
choose $c_{\ell+1}$ such that $0<c_{\ell+1}<\e_\ell2^{-\ell-1}$ and 
$\hn(A^{\ell}_{c_{\ell+1}})=0$. Since 
$S(v^\ell+c_{\ell+1}u^{\ell+1})\cup A^{\ell}_{c_{\ell+1}}\eqhn S(v^\ell)\cup 
S(u^{\ell+1})$,
we have
$$
S(v^{\ell+1})=S(v^\ell+c_{\ell+1}u^{\ell+1})\eqhn S(v^\ell)\cup 
S(u^{\ell+1})\,.
$$
We deduce by induction that
\begin{equation}\label{induction}
S(v^\ell)\eqhn \bigcup_{i=1}^\ell S(u^i)
\end{equation}
for every $\ell$.
Let
$$
E_\ell:=\{x\in S(v^\ell):|[v^\ell]|<\e_\ell\}\quad\hbox{and}\quad
F_\ell:=\bigcup_{i=\ell}^\infty E_i\,.
$$
By (\ref{em}) we have $\hn(E_\ell)<2^{-\ell}$, hence $\hn(F_\ell)<2^{1-\ell}$.

Let us prove that
\begin{equation}\label{svm}
S(v^\ell)\subsethn S(u)\cup E_\ell\,.
\end{equation}
Indeed, if $x\in S(v^\ell)\setmeno E_\ell$, then 
$|[v^\ell](x)|\ge\e_\ell$. On the other hand, for $\hn$-a.e.\ 
$x\in\Ga$ we have
$$
[u](x)=[v^\ell](x)+\sum_{i=\ell+1}^\infty c_i[u^i](x)\,,
$$
hence
\begin{equation}\label{vx}
|[u](x)|\ge \e_\ell-\Big|\sum_{i=\ell+1}^\infty c_i[u^i](x)\Big|\,.
\end{equation}
Since $|[u^i](x)|<1$ for $\hn$-a.e.\ $x\in\Ga$ we have
$$
\Big|\sum_{i=\ell+1}^\infty c_i[u^i](x)\Big|\le 
\sum_{i=\ell+1}^\infty c_i<\e_\ell\sum_{i=\ell+1}^\infty 2^{-i}\le 
\e_\ell\,,
$$
and we deduce from (\ref{vx}) that $[u](x)\neq0$, hence $x\in S(u)$ 
for $\hn$-a.e.\ $x\in S(v^\ell)\setmeno E_\ell$. We conclude that 
(\ref{svm}) is satisfied.
By (\ref{induction}) and (\ref{svm}) we obtain that
$$
\bigcup_{i=1}^k S(u^i)\subsethn S(u)\cup F_\ell
$$
for every $k\ge \ell$. Taking the union with respect to $k$ we get
$\Ga\subsethn S(u)\cup F_\ell$, which implies $\hn(\Ga\setmeno 
S(u))\le 2^{1-\ell}$. By the arbitrariness of $\ell$ we obtain 
$\Ga\subsethn S(u)$.

Under the additional hypotheses of the second part of the lemma
we may assume that
\begin{equation}\label{limuik2}
\|u_k^i\|_\infty\le 1\,,\qquad
\|\nabla u_k^i\|_p\le 1\,.
\end{equation}
and we define
$$
v_k^\ell:=\sum_{i=1}^\ell c_iu_k^i\,,\qquad u_k:=\sum_{i=1}^\infty c_iu_k^i\,.
$$
Since $v_k^\ell\in SBV^p(U)\cap L^\infty(U)$ and $S(v_k^\ell)\subsethn \Ga_k$,
from Theorem~\ref{compsbv} we obtain that $u_k\in SBV^p(U)\cap 
L^\infty(U)$ and $v_k^\ell$ converges to $u_k$ weakly in $SBV^p(U)$ 
as $\ell\to\infty$ and
$u_k$ converges to $u$ weakly in $SBV^p(U)$ as $k\to \infty$. By 
Remark~\ref{Su} we have $S(u_k)\subsethn\Ga_k$ for every~$k$.
\end{proof}

\begin{proposition}\label{gsbvrp}
Assume that $\Ga_k$ $\sigmap$-converges to $\Ga$ in $U$ and that 
$u_k$ converges to $u$ weakly in $GSBV^p(U;\Rm)$. If $S(u_k)\subsethn 
\Ga_k$ for every~$k$, then $S(u)\subsethn\Ga$.
\end{proposition}

\begin{proof}
Assume that $S(u_k)\subsethn \Ga_k$ for every~$k$. Let 
$\psi_j\in C^\infty_c(\Rm;\Rm)$ with $\psi_j(z)=z$ for 
$|z|\le j$. Then $\psi_j(u_k)$ converges to 
$\psi_j(u)$ weakly in $SBV^p(U;\Rm)$. Since 
$S(\psi_j(u_k))\subsethn S(u_k)\subsethn\Ga_k$, if we apply 
condition (a) of Definition~\ref{setsbv} to each component we obtain
$S(\psi_j(u))\subsethn \Ga$. The conclusion follows from the 
fact that $S(u)\eqhn\bigcup_j S(\psi_j(u))$.
\end{proof}

\subsection{Compactness properties}\label{compa}
We now prove some compactness properties of the $\sigmap$-con\-vergence.

\begin{theorem}\label{compattezza}
Every sequence $\Ga_k\subset U$, with $\hn(\Ga_k)$ uniformly bounded, 
has a $\sigmap$-con\-vergent subsequence.
\end{theorem}

\begin{proof}Let $\Ga_k$ be a sequence of subsets of $U$ with 
$\hn(\Ga_k)\le C<+\infty$ for every~$k$. Let $w_h$ be a sequence in 
$L^\infty(U)$ with the following density property: for every $w\in 
L^\infty(U)$ there exists a subsequence $w_{h_i}$ which converges to 
$w$ strongly in $L^p(U)$ and satisfies the inequality 
$\|w_{h_i}\|_\infty\le\|w\|_\infty$ for every $i$.
For every positive integers $\ell$, $h$, and $k$ let $u_k^{\ell,h}$ 
be the solution of the following minimization problem
$$
\min\{\|\nabla u\|_p^p+\ell\|u-w_h\|_p^p: u\in SBV^p(U),\ 
S(u)\subsethn\Ga_k\}\,.
$$
To prove the existence of a solution it is enough to apply 
Theorem~\ref{compsbv} and Remark~\ref{Su} to a minimizing sequence; the uniqueness 
follows from the strict convexity of the functional. By a truncation argument we obtain $\|u_k^{\ell,h}\|_{\infty}\le\|w_h\|_{\infty}$. 
By Theorem~\ref{compsbv} and by a diagonal argument, passing to a 
subsequence, we may assume that $u_k^{\ell,h}$ converges weakly in 
$SBV^p(U)$ to a function $u^{\ell,h}\in SBV^p(U)$ as $k\to\infty$.
Let
$$
\Ga:=\bigcup_{\ell,h=1}^\infty S(u^{\ell,h})\,.
$$
To prove that $\hn(\Ga)<+\infty$, for every integer $r>0$ we consider 
the sequence $v^r_k:=(u^{\ell,h}_k)_{1\le\ell, h\le r}$ and the 
function $v^r:=(u^{\ell,h})_{1\le\ell, h\le r}$. As $v^r_k$ converges 
to $v^r$ weakly in $SBV^p(U;{\mathbb M}^{r\times r})$ we have
\begin{eqnarray*}
&\displaystyle \hn(\bigcup_{\ell,h=1}^r S(u^{\ell,h}))=
\hn(S(v^r)) \le \liminf_{k\to\infty} \hn(S(v^r_k)) = \\
&\displaystyle = \liminf_{k\to\infty} \hn(\bigcup_{\ell,h=1}^r 
S(u_k^{\ell,h}))
\le \liminf_{k\to\infty} \hn(\Ga_k)\le C\,.
\end{eqnarray*}
Passing to the limit as $r\to\infty$ we obtain $\hn(\Ga)\le C$.

We want to prove that $\Ga_k$ $\sigmap$-converges to $\Ga$. By 
construction $\Ga$ satisfies all hypotheses of Lemma~\ref{cond2}, 
which implies that condition (b) of Definition~\ref{setsbv} is 
satisfied.

To prove condition (a) let $v_j$ be a sequence which converges weakly 
in $SBV^p(U)$ to a function $v$, and with $S(v_j)\subsethn\Ga_{k_j}$ 
for some sequence $k_j\to\infty$. We have to show that 
$S(v)\subsethn\Ga$.
By the density property of $w_h$ there exists a subsequence $w_{h_i}$ 
which converges to $v$ in $L^p(U)$ and such that 
$\|w_{h_i}\|_\infty\le\|v\|_\infty<+\infty$. Let 
$\ell_i\to\infty$ such that $\ell_i\|w_{h_i}-v\|_p^p\to0$.
{}From the minimality of $u_{k_j}^{\ell_i,h_i}$ we obtain 
$\|u_{k_j}^{\ell_i,h_i}\|_\infty\le
\|w_{h_i}\|_\infty\le\|v\|_\infty$ and
$$
\|\nabla 
u_{k_j}^{\ell_i,h_i}\|_p^p+\ell_i\|u_{k_j}^{\ell_i,h_i}-w_{h_i}\|^p_p\le
\|\nabla v_j\|^p_p+\ell_i\|v_j-w_{h_i}\|_p^p\,.
$$
Passing to the limit as $j\to\infty$ we get 
$\|u^{\ell_i,h_i}\|_\infty\le\|v\|_\infty$ and
\begin{equation}\label{Yosida}
\|\nabla u^{\ell_i,h_i}\|_p^p+\ell_i\|u^{\ell_i,h_i}-w_{h_i}\|^p_p\le
M+\ell_i\|v-w_{h_i}\|_p^p\,,
\end{equation}
where $M:=\sup_j\|\nabla v_j\|^p_p<+\infty$.
This inequality implies that $\nabla u^{\ell_i,h_i}$ is bounded in 
$L^p(U;\R^n)$ uniformly with respect to $i$ and that 
$u^{\ell_i,h_i}-w_{h_i}$ tends to $0$ in $L^p(U)$. Since $w_{h_i}$ 
converge to $v$ in $L^p(U)$, we conclude that $u^{\ell_i,h_i}$ 
converge to $v$ in $L^p(U)$. We now apply Remark~\ref{Su} to the 
sequence $u^{\ell_i,h_i}$ and conclude that $S(v)\subsethn\Ga$.
\end{proof}

We shall use the following extension of the compactness theorem, that 
can be proved by adapting the arguments of  Helly's theorem.

\begin{theorem}\label{helly}
Let $t\mapsto\Ga_k(t)$ be a sequence of increasing set functions 
defined on an interval $I\subset\R$ with values contained in $U$, i.e.,
$$
\Ga_k(s)\subsethn\Ga_k(t)\subset U\quad\hbox{for every $s,\, t\in I$ 
with $s<t$.}
$$
Assume that the measures
$\hn(\Ga_k(t))$ are bounded uniformly with respect to $k$ and $t$. 
Then there exist a subsequence $\Ga_{k_j}$ and an increasing set 
function $t\mapsto\Ga(t)$ on $I$ such that
\begin{equation}\label{gkj}
\Ga_{k_j}(t)\quad\sigmap\hbox{-converges to }\Ga(t) \hbox{ in }U
\end{equation}
for every $t\in I$.
\end{theorem}
\begin{proof}
Let $D$ be a countable dense set in $I$. Using 
Theorem~\ref{compattezza} and a diagonal argument we can extract a 
subsequence $\Ga_{k_j}$ such that
$\Ga_{k_j}(t)$ $\sigmap$-converges to some set $\Ga(t)$ in $U$
for every $t\in D$. {}From Remark~\ref{sigmap} we have 
$\Ga(s)\subsethn\Ga(t)$
for every $s,\, t\in D$ with $s<t$. By Theorem~\ref{semigamma} the measures
$\hn(\Ga(t))$ are bounded uniformly with respect to $t\in D$.
For every $t\in I\setmeno D$ we define
$$
\Ga_{-}(t):=\bigcup_{s<t,\,s\in D}\Ga(s)\,,\qquad 
\Ga_{+}(t):=\bigcap_{s>t,\,s\in D}\Ga(s)\,.
$$
Then 
$\Ga_{+}(s)\subsethn\Ga_{-}(t)\subsethn\Ga_{+}(t)$ for every $s$, 
$t\in I\setmeno D$ with $s<t$, and the measures 
$\hn(\Ga_{+}(t))$ are bounded uniformly with respect to $t\in 
I\setmeno D$. This implies that the set $D_1:=\{t\in I\setmeno 
D:\hn(\Ga_{+}(t)\setmeno\Ga_{-}(t))>0\}$ is at most countable.

For every $t\in I\setmeno (D\cup D_1)$ we define
$\Ga(t):= \Ga_{+}(t)\eqhn\Ga_{-}(t)$.
Given $t\in I\setmeno (D\cup D_1)$, let us prove that $\Ga_{k_j}(t)$ 
$\sigmap$-converges to $\Ga(t)$. If $u_i$ converges weakly to $u$ in 
$SBV^p(U)$ and $S(u_i)\subsethn\Ga_{k_{j_i}}(t)$ for some sequence 
$j_i\to\infty$, then
$S(u_i)\subsethn\Ga_{k_{j_i}}(s)$ for every $s>t$, $s\in D$; by 
$\sigmap$-convergence, this implies $S(u)\subsethn \Ga(s)$, and 
taking the intersection for $s>t$, $s\in D$,
we obtain
$S(u)\subsethn \Ga_{+}(t)\eqhn\Ga(t)$,
so that condition (a) of Definition~\ref{setsbv} is satisfied.

To prove condition (b), we observe that for every $s<t$, $s\in D$ 
there exists a function $u(s)\in SBV^p(U)$ and a sequence $u_k(s)$ 
converging to $u(s)$ weakly in $SBV^p(U)$ such that $S(u(s))\eqhn 
\Ga(s)$ and $S(u_k(s))\subsethn\Ga_{k_j}(s)\subsethn\Ga_{k_j}(t)$. 
Condition (b) of Definition~\ref{setsbv} follows now from 
Lemma~\ref{cond2}.

Since the set $D_1$ is at most countable, using 
Theorem~\ref{compattezza} and a diagonal argument it is possible to 
extract a further subsequence, still denoted $\Ga_{k_j}$, such that
$\Ga_{k_j}(t)$ $\sigmap$-converges to some set $\Ga(t)$ for every $t\in D_1$.
This concludes the proof of (\ref{gkj}) for every $t\in I$. The 
monotonicity of $t\mapsto\Ga(t)$ on $I$ follows from Remark~\ref{sigmap}.
\end{proof}

\subsection{Some results in measure theory}\label{mth}

We begin with a lemma concerning perturbations of bounded sequences 
in $L^p$ spaces.

\begin{lemma}\label{perturbation}
Let $(X,{\mathcal A}, \mu)$ be a finite measure space, let $p>1$, let $m,\,n\ge 1$,
and let $H\colon {X{\times}\Rn}\to \Rm$ be a Carath\'eodory function. 
Assume that there exist a constant $a\ge 0$ and a nonnegative 
function $b\in L^{p'}(X)$, with $p'=p/(p-1)$, such that
\begin{equation}\label{H1}
|H(x,\xi)|\le a|\xi|^{p-1}+b(x)
\end{equation}
for every $(x,\xi)\in {X{\times}\Rn}$. Let $\Phi_k$ and $\Psi_k$ be 
two sequences in $L^p(X;\Rn)$. Assume that $\Phi_k$ is bounded in 
$L^p(X;\Rn)$ and $\Psi_k$ converges to $0$ strongly in $L^p(X;\Rn)$. 
Then
\begin{equation}\label{H2}
\int_X [H(x,\Phi_k(x)+\Psi_k(x))-H(x,\Phi_k(x))]\,\Phi(x)\,d\mu(x)\to 0
\end{equation}
for every $\Phi\in L^p(X;\Rm)$.
\end{lemma}

\begin{proof}
There exists a constant $C\ge0$ such that
\begin{equation}\label{H3}
\|\Phi_k\|_p \le C \qquad\hbox{and}\qquad \|\Phi_k+\Psi_k\|_p\le C
\end{equation}
for every $k$. Let us fix $\Phi\in L^p(X;\Rm)$. By the absolute 
continuity of the integral, for every $\e>0$ there exists $\delta>0$ 
such that
\begin{equation}\label{H4}
A\in{\mathcal A},\,\ \mu(A)<\delta\qquad \Longrightarrow \qquad 
\int_A|\Phi(x)|^pd\mu(x)<\e^p \,.
\end{equation}

Let us fix $M>0$ such that
\begin{equation}\label{H3.5}
2\,C^p/M^p<\delta\,.
\end{equation}
For every $x\in X$ and $\eta>0$ let
\begin{equation}\label{H5}
\om(x,\eta):=\max\,\{|H(x,\xi_1)-H(x,\xi_2)|: |\xi_1|\le M, \ |\xi_2|\le M, \
|\xi_1-\xi_2| \le\eta\}\,.
\end{equation}
Since $H$ is a Carath\'eodory function and satisfies (\ref{H1}), it 
turns out that $\om$ is a Carath\'eodory function, $\om(x,0)=0$, and
$$
0\le \om(x,\eta)\le 2aM^{p-1}+2b(x)
$$
for every $x\in X$ and every $\eta>0$. As $\Psi_k$ converges to $0$ 
strongly in $L^p(X;\Rn)$ and $\om(x,\eta)\to0$ as $\eta\to0$, we have
\begin{equation}\label{H7}
\int_X \om(x,|\Psi_k(x)|)\, |\Phi(x)|\, d\mu(x) \to 0\,.
\end{equation}

Let
\begin{eqnarray}
&A_k:=\{x\in X: |\Phi_k(x)+\Psi_k(x)|> M\}\cup \{x\in X: |\Phi_k(x)|> M\}\,,
\label{H8}\\
&B_k:=\{x\in X: |\Phi_k(x)+\Psi_k(x)|\le M\}\cap \{x\in X: 
|\Phi_k(x)|\le M\}\,.
\label{H8.5}
\end{eqnarray}
By (\ref{H1}) we have
\begin{equation}\label{H9}
|H(x,\Phi_k(x)+\Psi_k(x))-H(x,\Phi_k(x))|\le
a|\Phi_k(x)+\Psi_k(x)|^{p-1}+ a|\Phi_k(x)|^{p-1} + 2b(x)
\end{equation}
for every $x\in A_k$.
By (\ref{H5}) and (\ref{H8.5}) we have
\begin{equation}\label{H10}
|H(x,\Phi_k(x)+\Psi_k(x))-H(x,\Phi_k(x))|\le \om(x,|\Psi_k(x)|)
\end{equation}
for every $x\in B_k$.
Using the H\"older inequality, from (\ref{H9}) and (\ref{H10})
we obtain
\begin{eqnarray}
&\displaystyle
\int_X |H(x,\Phi_k(x)+\Psi_k(x))-H(x,\Phi_k(x))|\, 
|\Phi(x)|\,d\mu(x)\le \label{H11}\\
&\displaystyle \le
K \Big(\int_{A_k}|\Phi(x)|^p d\mu(x)\Big)^\frac{1}{p}
+\int_{B_k} \om(x,|\Psi_k(x)|)\, |\Phi(x)|\, d\mu(x) \,, \nonumber
\end{eqnarray}
where $K:=2(aC^{p-1}+\|b\|_{p'}) $.

By  Chebyshev's inequality, (\ref{H3}) and (\ref{H8}) imply 
$\mu(A_k)\le 2\,C^p/M^p$, and thus, by 
 (\ref{H3.5}), $\mu(A_k)<\delta$. Therefore,
from (\ref{H4}) and (\ref{H11}) we deduce that
\begin{eqnarray*}
&\displaystyle
\int_X |H(x,\Phi_k(x)+\Psi_k(x))-H(x,\Phi_k(x))|\, |\Phi(x)|\,d\mu(x)\le \\
&\displaystyle \le K\e
+ \int_{X} \om(x,|\Psi_k(x)|)\, |\Phi(x)|\, d\mu(x)\,.
\end{eqnarray*}
Taking (\ref{H7}) into account, we 
obtain (\ref{H2}) upon passing to the limit  in the previous inequality first as ${k\to\infty}$ and 
then as $\e\to 0$.
\end{proof}

\begin{remark}\label{nonat}
Let $(X,{\mathcal A}, \mu)$ be a finite nonatomic measure space, let $p>1$, let $m,\, n\ge 1$,
and let $H\colon {X{\times}\Rn}\to \Rm$ be a Carath\'eodory function. Assume that the function $x\mapsto H(x,\Phi(x))\,\Psi(x)$ is $\mu$-integrable for every $\Phi\in L^p(X;\Rn)$ and $\Psi\in L^p(X;\Rm)$.  Then,  the function $x\mapsto H(x,\Phi(x))$ belongs to $L^{p'}(X;\Rm)$ for every $\Phi\in L^p(X;\Rn)$, and this implies (\ref{H1}) by a classical result in the theory of integral operators (see \cite[Theorem~2.3 in 
Chapter~I]{Kra}).

In particular,  the conclusion of Lemma~\ref{perturbation} still holds.

\end{remark}

Let $U$ and $W$ be as in Theorem~\ref{semigsbv}. Assume in addition that
$\xi\mapsto W(x,\xi)$ belongs to $C^1(\Mmn)$ for every $x\in U$.
Since $\xi\mapsto
W(x,\xi)$ is rank-one convex on $\Mmn$ for every $x\in U$ (see, e.g.,
\cite{Dac}), from (\ref{growth}) we can deduce that
there exist a constant $a_2>0$
and a nonnegative function $b_2\in L^{p'}(U)$ such that
\begin{equation}\label{W2}
|\partial_\xi W(x,\xi)|\le a_2|\xi|^{p-1}+b_2(x)
\end{equation}
for every $(x,\xi)\in U{\times}\Mmn$.

Let us consider the $C^1$ functional $\W\colon
L^p(U;\Mmn)\to\R$ defined by
$$
\W(\Phi):=\int_U W(x,\Phi(x))\,dx
$$
whose differential $\partial \W\colon L^p(U;\Mmn)\to
L^{p'}(U;\Mmn)$ is given by
$$
\langle \partial \W(\Phi), \Psi\rangle =\int_U \partial_\xi
W(x,\Phi(x))\Psi(x)\,dx\,,
$$
for every $\Phi,\, \Psi\in L^p(U;\Mmn)$, where
$\langle\cdot,\cdot\rangle$ denotes the duality pairing between the
spaces $L^{p'}(U;\Mmn)$ and $L^{p}(U;\Mmn)$.

\begin{lemma}\label{stressconv}
Assume that $u_k$ converges to $u$ weakly in $GSBV^p(U;\Rm)$ and that $\W(\nabla u_k)$ converges to $\W(\nabla u)$. Then $\partial\W(\nabla u_k)$ converges to $\partial\W(\nabla u)$ weakly in $L^{p'}(U;\Mmn)$.
\end{lemma}

\begin{proof}
It is enough to prove that
\begin{equation}\label{401}
\langle \partial\W(\nabla u),\Psi\rangle
\le \liminf_{k\to\infty} \,
\langle \partial\W(\nabla u_k),\Psi\rangle
\end{equation}
for every $\Psi\in L^p(U;\Mmn)$. Let $\eta_i$ be a sequence of positive numbers converging to $0$. If we apply the lower semicontinuity theorem for $GSBV$ (Theorem~\ref{semigsbv}) to the function $W(x,\xi+\eta_i\Psi(x))$, for every $i$ we obtain
\begin{equation}\label{401.5}
\frac{\W(\nabla u+ \eta_i \Psi)-\W(\nabla u)}{\eta_i}\le
\liminf_{k\to\infty}\, \frac{\W(\nabla u_k+ \eta_i \Psi)-\W(\nabla u_k)}{\eta_i}\,.
\end{equation}
Therefore there exists an increasing sequence of integers $k_i$ such that
\begin{equation}\label{402}
\frac{\W(\nabla u+ \eta_i \Psi)-\W(\nabla u)}{\eta_i} - \frac{1}{i} \le
\frac{\W(\nabla u_k+ \eta_i \Psi)-\W(\nabla u_k)}{\eta_i}
\end{equation}
for every $k\ge k_i$.
Defining $\e_k:=\eta_i$ for $k_i\le k< k_{i+1}$, from (\ref{402}) we obtain
\begin{equation}\label{403}
\liminf_{k\to\infty}\,\frac{\W(\nabla u+ \e_k\Psi)-\W(\nabla u)}{\e_k} \le
\liminf_{k\to\infty}\, \frac{\W(\nabla u_k+ \e_k \Psi)-\W(\nabla u_k)}{\e_k}\,.
\end{equation}
Since $\W$ is of class $C^1$ on $L^p(U;\Mmn)$, we have
\begin{equation}\label{404}
\langle \partial\W(\nabla u),\Psi\rangle=
\lim_{k\to\infty} \,\frac{\W(\nabla u+ \e_k\Psi)-\W(\nabla u)}{\e_k}
\end{equation}
and
\begin{equation}\label{405}
\frac{\W(\nabla u_k+ \e_k \Psi)-\W(\nabla u_k)}{\e_k}=
\langle \partial\W(\nabla u_k+\tau_k\Psi),\Psi\rangle
\end{equation}
for suitable constants $\tau_k\in[0,\e_k]$. By (\ref{W2}) and by Lemma~\ref{perturbation} we have
\begin{equation}\label{406}
\liminf_{k\to\infty}\,
\langle \partial\W(\nabla u_k+\tau_k\Psi),\Psi\rangle =
\liminf_{k\to\infty}\, \langle \partial\W(\nabla u_k),\Psi\rangle\,.
\end{equation}
Inequality (\ref{401}) follows now from (\ref{403})--(\ref{406}).
\end{proof}

If $W$ is strictly convex, using
\cite[Theorem 2]{Vis} one can prove that
$\nabla u_k$ converges in measure to $\nabla u$ on~$U$. 
By (\ref{W2}) this implies that
$\partial\W(\nabla u_k)$ converges to $\partial\W(\nabla u)$ weakly in $L^{p'}(U;\Mmn)$. We refer to \cite{Bre2} for similar results with $p=1$.

\subsection{Approximation by Riemann sums}
\label{aRiem}
We now prove a lemma concerning the approximation of Lebesgue
integrals by Riemann sums. The convergence result (\ref{riemsum2}) is well-known (see \cite{Hahn}). For the application we have in mind we need the stronger result~(\ref{riemsum*}), that is related to the
Saks-Henstock lemma (see \cite{Saks} and \cite{Hen}) used in the theory of Henstock-Kurzweil integral (see, e.g., \cite{Maw}). We prefer to present an independent proof, based on \cite[page~63]{Doo}, which only uses  Fubini's theorem.

\begin{lemma}\label{riemann1}
Let $[a,b]$ be a closed bounded  interval, let $X$ be a Banach space,
and let $f\colon[a,b]\to X$ be a Bochner integrable
function. Then there exists a sequence of subdivisions $(\tki)_{0\le i\le i_k}$
of the interval $[a,b]$, with
\begin{eqnarray}
& a=t_k^0<t_k^1<\cdots<t_k^{i_k-1}<t_k^{i_k}=b\,,
\label{subdiv}\\
&\displaystyle
\lim_{k\to\infty}\,
\max_{1\le i\le i_k} (\tki-\tkim)= 0\,,
\label{fine}
\end{eqnarray}
such that
\begin{equation}\label{riemsum}
\lim_{k\to\infty}\,\sum_{i=1}^{i_k}\int_{\tkim}^{\tki} \| f(\tki) -
 f(t)\|\,dt \,= 0\,.
\end{equation}
 In particular we have
\begin{eqnarray}
&\displaystyle \sum_{i=1}^{i_k} \Big \| (\tki-\tkim) f(\tki) -
\int_{\tkim}^{\tki} f(t)\,dt \,\Big\|\ \longrightarrow \ 0\,,\label{riemsum*}\\
&\displaystyle \sum_{i=1}^{i_k} (\tki-\tkim) f(\tki)\ \longrightarrow \ \int_a^b
f(t)\,dt \qquad\hbox{strongly in }\, X \label{riemsum2}
\end{eqnarray}
as $k\to\infty$.
\end{lemma}

\begin{proof}
We extend $f$ to $0$ outside $[a,b]$.
Set, for every $m\ge 1$ 
and  $i\in {\mathbb Z}$,  $\taumi:=i/m$.
For every $s\in [0,1]$ we have
\begin{eqnarray*}
&\displaystyle
\sum_{i\in {\mathbb Z}}
\int_{s+\taumim}^{s+\taumi}\|f(s+\taumi)-f(t)\|\,dt=\\
&\displaystyle
=\sum_{i\in {\mathbb Z}}
\int_0^{\frac{1}{m}} \|f(s+\taumi)-f(s+\taumi-\tau)\|\,d\tau\,.
\end{eqnarray*}
Note that there are  at most $m(b-a+1)+2$ non-zero elements in the above sums, namely those with $i \in I_m:=\{i\in{\mathbb Z}:m(a-1)\le i\le mb+1\}$. 
Integrating with respect to $s$ we obtain
\begin{eqnarray}
&\displaystyle
\int_{0}^{1}\Big[\sum_{i \in  {\mathbb Z}}\int_{s+\taumim}^{s+\taumi}  \|  f(s+\taumi) -
 f(t)\|\,dt \,\Big]\,ds\le\nonumber\\
&\displaystyle
\le\sum_{i \in I_m}
\int_0^{\frac{1}{m}}
\Big[\int_{-\infty}^{+\infty}\|f(s+\taumi)-f(s+\taumi-\tau)\|\,ds\Big]\,d\tau=\label{3500}\\
&\displaystyle
=\sum_{i \in I_m}
\int_0^{\frac{1}{m}}
\Big[\int_{-\infty}^{+\infty}\|f(s)-f(s-\tau)\|\,ds\Big]\,d\tau \,.\nonumber
\end{eqnarray}
By  continuity of the translations, for every $\e>0$ there exists
$\delta>0$ such that
\begin{equation}\label{3501}
\int_{-\infty}^{+\infty}\|f(s)-f(s-\tau)\|\,ds<\e
\end{equation}
for $0<\tau<\delta$.
Thus,  from~(\ref{3500})
and~(\ref{3501}) we obtain
$$
\lim_{m\to\infty}\int_{0}^{1}\Big[\sum_{i\in {\mathbb Z}}
\int_{s+\taumim}^{s+\taumi}\|f(s+\taumi)-f(t)\|\,dt\Big]\,ds=0\,.
$$
Therefore there exists a sequence $m_k\to\infty$ such that
\begin{equation}\label{3502}
\lim_{k\to\infty}\sum_{i\in {\mathbb Z}}
\int_{s+\tau_{m_k}^{i-1}}^{s+\tau_{m_k}^i}\| f(s+\tau_{m_k}^i) -f(t)\|\,dt=0
\end{equation}
for a.e.\ $s\in[0,1]$.
Let us fix $s\in[0,1]$ such that (\ref{3502}) holds. Let $\rho_k$ be
the largest integer $i$ such that $s+\tau_{m_k}^i\le a$, and let
$\sigma_k$ be the smallest integer $i$ such that $s+\tau_{m_k}^i\ge b$,
and let $i_k:=\sigma_k-\rho_k$. For $i=1,\ldots,i_k-1$ we define
$\tki:=s+\tau_{m_k}^{\rho_k+i}$ and we set $t_k^0:=a$ and
$t_k^{i_k}:=b$. Then (\ref{subdiv}) and (\ref{fine}) are satisfied.
Moreover
\begin{eqnarray}
&\displaystyle
\sum_{i=1}^{i_k}\int_{\tkim}^{\tki} \| f(\tki) -
 f(t)\|\,dt \,=\nonumber\\
&\displaystyle
=\sum_{i=\rho_k+2}^{\sigma_k-1}\int_{s+\tau_{m_k}^{i-1}}^{s+\tau_{m_k}^i} \|  f(s+\tau_{m_k}^i) -
 f(t)\|\,dt \,+{}\label{3503}\\
&\displaystyle{}
+\int_a^{a_k}\|f(a_k)-
f(t)\|\,dt+
\int_{b_k}^b\|f(b)-
f(t)\|\,dt\,,\nonumber
\end{eqnarray}
where $a_k:=s+\tau_{m_k}^{\rho_k+1}$ and $b_k:=s+\tau_{m_k}^{\sigma_k-1}$.
Since all integers between $\rho_k+2$ and $\sigma_k-1$ belong to
$I_{m_k}$, the first term in the right hand side of (\ref{3503})
tends to $0$ by (\ref{3502}). The second term is estimated by
$$
\int_a^{a_k}\|f(a_k)-f(t)\|\,dt\le
 \int_{s+\tau_{m_k}^{\rho_k}}^{s+\tau_{m_k}^{\rho_k+1}}
\|f(s+\tau_{m_k}^{\rho_k+1})-f(t)\|\,dt
$$
which also tends to $0$ by~(\ref{3502}). The third
term tends to $0$ by the absolute continuity of the integral, since
$b-b_k$ tends to~$0$ by the choice of $\sigma_k$. This concludes the
proof of (\ref{riemsum}).
\end{proof}

\begin{remark}\label{riemann2}
If $X_j$ is a sequence of Banach spaces and $f_j\colon[a,b]\to X_j$
is a sequence of Bochner integrable function, then there exists a
sequence of subdivisions $(\tki)_{0\le i\le i_k}$, independent of $j$
and satisfying (\ref{subdiv}) and (\ref{fine}), such that
(\ref{riemsum}) is satisfied simultaneously for each function $f_j$.
Indeed, we can consider the
Banach space $X$ of all sequences $x:=(x_j)$ such that $x_j\in X_j$
for every $j$ and $\sum_j
\|x_j\|_{X_j}<+\infty$, endowed with the norm $\|x\|_X:= \sum_j
\|x_j\|_{X_j}$. To obtain the result it is enough to
apply Lemma~\ref{riemann1} to the function $g\colon [a,b]\to X$ whose
components $g_j$ are given by $g_j(t):=2^{-j}f_j(t)/\|f_j\|_1$, where
$\|f_j\|_1:=\int_a^b\|f_j(t)\|_{X_j}\,dt$.
\end{remark}

\end{section}

\begin{section}{PRELIMINARY RESULTS}\label{preliminary}

We now return to the framework described in Section~\ref{formulation} 
and adapt to it the tools developed in Section~\ref{tools}. Moreover, 
we extend the jump transfer results of \cite{Fra-Lar} to the space 
$GSBV(\Om;\Rm)$, and use them to
prove the stability of minimum energy configurations (see
Subsection~\ref{minen}).

\subsection{Jump transfer}\label{jtr}
To deal with the interaction
between cracks and boundary deformations it is convenient to extend 
all deformations to a bounded open set $\Om_0$, containing
$\overline\Om$ and with Lipschitz boundary. When we speak of $\sigmap$-convergence we always refer to $\sigmap$-convergence in~$\Om_0$.

For every rectifiable set $\Ga\subset\Rn$ 
we define 
\begin{equation}\label{gammanb}
\Ga^N:=\Ga\cup\partial_N\Om\,.
\end{equation}
{}From (\ref{calK}) we obtain that
\begin{equation}\label{eqganga}
\K(\Ga^N)=\K(\Ga)
\end{equation}
for every rectifiable set $\Ga\subset\Rn$. 
Theorem~\ref{semigamma} implies that
\begin{equation}\label{semiK}
\K(\Ga\cup\Ga')\le \liminf_{k\to\infty} \K(\Ga_k\cup\Ga')
\end{equation}
whenever $\Ga_k$, $\Ga$, and $\Ga'$ are rectifiable sets in $\overline\Om$, $\Ga_k$ $\sigmap$-converges to $\Ga$, and $\hn(\Ga')<+\infty$.

{}From \cite[Theorem~2.1]{Fra-Lar} we shall obtain the following result.

\begin{theorem}[Jump transfer in $SBV$]\label{fl}
Assume that $\Ga_k\in\Rp(\overline\Om_B)$ and that $\Ga_k^N$
$\sigmap$-con\-verges to $\Ga$. Then for every function
$v\in SBV^p(\Om_0;\Rm)$ there exists a sequence $v_k\in
SBV^p(\Om_0;\Rm)$ such that
\begin{itemize}
\item[{\rm (a)}] $v_k=v$ a.e.\ in $\Om_0\setmeno\Om_B$,
\item[{\rm (b)}] $v_k\to v$ strongly in $L^1(\Om_0;\Rm)$,
\item[{\rm (c)}] $\nabla v_k\to\nabla v$ strongly in $L^p(\Om_0;\Mmn)$,
\item[{\rm (d)}] $\hn((S(v_k)\setmeno\Ga_k^N) \setmeno(S(v)\setmeno\Ga^N))\to0$.
\end{itemize}
If, in addition, $v\in L^\infty(\Om_0;\Rm)$, then we may assume that
$v_k$ is bounded in $L^\infty(\Om_0;\Rm)$.
\end{theorem}

To prove Theorem~\ref{fl} we need the following lemma.

\begin{lemma}\label{dnom}
Assume that $\Ga_k\in\Rp(\overline\Om_B)$ and that $\Ga_k^N$
$\sigmap$-con\-verges to $\Ga$. Then
there exist a function
$w\in SBV^p(\Om_0)$ and a sequence $w_k$ converging to $w$ weakly in
$SBV^p(\Om_0)$ such that $S(w)\eqhn{\Ga\setmeno\partial_N\Om}$ and
$S(w_k)\subsethn\Ga_k$ for every $k$.
\end{lemma}

\begin{proof}
Let $\varphi_i\in C^\infty(\Rn)$ with $\varphi_i(x)=1$ if ${\rm
dist}(x,\partial_N\Om)>1/i$, and $\varphi_i(x)=0$ in a neighbourhood
of $\partial_N\Om$. From condition (b) in Definition~\ref{setsbv}
there exist a function $v\in SBV^p(\Om_0)$ and a sequence $v_k$
converging to $v$ weakly in $SBV^p(\Om_0)$ such that
$S(v)\eqhn\Ga$ and
$S(v_k)\subsethn\Ga_k^N$ for every $k$. Let
$v^i:=\varphi_iv$ and $v^i_k:=\varphi_iv_k$. Then $v^i_k$ converges to
$v^i$ weakly in $SBV^p(\Om_0)$, $S(v^i_k)\subsethn\Ga_k$, and
${\Ga\setmeno\partial_N\Om}\eqhn\bigcup_iS(v^i)$. As $\hn(\Ga)<+\infty$ by  
Remark~\ref{bgk}-(a), the conclusion follows from
Lemma~\ref{cond2}.
\end{proof}

\begin{proof}[Proof of Theorem~\ref{fl}] By assumption we have
$\Ga_k^N\subsethn\overline\Om_B\cup\partial_N\Om$ for every $k$.
As $\overline\Om_B\cup\partial_N\Om$ is closed, 
we deduce that $\Ga\subsethn\overline\Om_B\cup\partial_N\Om$.
By Lemma~\ref{dnom} there exist a function $w\in SBV^p(\Om_0)$
and a sequence $w_k$ converging to $w$ weakly in
$SBV^p(\Om_0)$ such that
$S(w)\eqhn{\Ga\setmeno\partial_N\Om} \subsethn\overline\Om_B$ and
$S(w_k)\subsethn\Ga_k \subsethn\overline\Om_B$ for every~$k$. 
Let us fix
$v=(v^1,\ldots,v^m)\in SBV^p(\Om_0;\Rm)$. If we apply \cite[Theorem
2.1]{Fra-Lar} to each component $v^i$ of $v$, with $\Om=\Om_B$,
$\Om'=\Om_0$, $u_k=w_k$, and $u=w$, we construct a sequence
$v_k=(v_k^1,\ldots,v_k^m)\in SBV^p(\Om_0;\Rm)$ which satisfies (a),
(b), (c), and
$$
\hn((S(v^i_k)\setmeno S(w_k))\setmeno(S(v^i)\setmeno S(w))) \to0
$$
for $i=1,\ldots,m$.
By monotonicity this implies
$$
\hn((S(v^i_k)\setmeno \Ga_k)\setmeno(S(v^i)\setmeno 
(\Ga\setmeno \partial_N\Om)) \to0\,.
$$
As $(S(v^i_k)\setmeno \Ga_k^N)\setmeno(S(v^i)\setmeno
\Ga^N)\subsethn (S(v^i_k)\setmeno \Ga_k)\setmeno(S(v^i)\setmeno
(\Ga\setmeno\partial_N\Om))$, we obtain
\begin{equation}\label{di}
\hn((S(v^i_k)\setmeno \Ga_k^N)\setmeno(S(v^i)\setmeno
\Ga^N))\to0\,.
\end{equation}
Since
$S(v_k)\eqhn\bigcup_i S(v^i_k)$
and $S(v)\eqhn\bigcup_i S(v^i)$,
we have
$$
(S(v_k)\setmeno\Ga_k^N)\setmeno(S(v)\setmeno\Ga^N)\subsethn
\bigcup_{i=1}^m
((S(v^i_k)\setmeno\Ga_k^N)\setmeno(S(v^i)\setmeno\Ga^N))\,,
$$
so that (d) follows from (\ref{di}).

If, in addition, $v\in L^\infty(\Om_0;\Rm)$, we can replace
$v_k$ by $\varphi(v_k)$, where $\varphi\in C^1_0(\Rm;\Rm)$ satisfies
$\varphi(z)=z$ for $|z|\le\|v\|_\infty$. The new sequence $\varphi(v_k)$
is bounded in $L^\infty(\Om_0;\Rm)$ and continues to satisfy (a)--(d).
\end{proof}

The following theorem extends the result of Theorem~\ref{fl} to the
case of $GSBV$ functions.

\begin{theorem}[Jump transfer in $GSBV$]\label{fl2}
Assume that $\Ga_k\in\Rp(\overline\Om_B)$ and that $\Ga_k^N$ 
$\sigmap$-con\-verges to $\Ga$. Then for every function
$v\in GSBV^p_{q}(\Om_0;\Rm)$ there exists a sequence $v_k\in
GSBV^p_{q}(\Om_0;\Rm)$ such that
\begin{itemize}
\item[{\rm (a)}] $v_k=v$ a.e.\ in $\Om_0\setmeno\Om_B$,
\item[{\rm (b)}] $v_k\to v$ strongly in $L^{q}(\Om_0;\Rm)$,
\item[{\rm (c)}] $\nabla v_k\to\nabla v$ strongly in $L^p(\Om_0;\Mmn)$,
\item[{\rm (d)}] $\hn((S(v_k)\setmeno\Ga_k^N) 
\setmeno(S(v)\setmeno\Ga^N)) \to0$.
\end{itemize}
\end{theorem}

\begin{proof} Since $\Om_0\setmeno \overline\Om_B$ has a
Lipschitz boundary, by Proposition~\ref{traces} for
$\hn$-a.e.\ $x\in \partial\Om_B$ there exists $\tilde v_B(x)\in\Rm$ such that
\begin{equation}\label{vbx}
\tilde v_B(x):=\aplim\limits_{{\scriptstyle y\to x},\ {\scriptstyle
y\notin\overline\Om_B}}v(y)\,.
\end{equation}
For every integer $i\ge1$ let
\begin{equation}\label{Ri}
S_i:=\{x\in \partial\Om_B: |\tilde v_B(x)|\ge i\}\,.
\end{equation}
Then
\begin{equation}\label{limSi}
\hn(S_i)\to 0\,.
\end{equation}

Let $\varphi\in C^1_c(\Rm;\Rm)$ be a function such that $\varphi(z)=z$
for $|z|\le 1$, and let $\varphi_i(z):=i\varphi(z/i)$. Then
$\varphi_i\in C^1_c(\Rm;\Rm)$,
$\varphi_i(z)=z$ for $|z|\le i$, and
$|\nabla \varphi_i|\le C$ for some constant independent of $i$. Since $v\in GSBV^p(\Om_0;\Rm)$, the functions $v^i:=\varphi_i(v)$ belong to
$SBV^p(\Om_0;\Rm)\cap L^\infty(\Om_0;\Rm)$.

By Theorem~\ref{fl} for every $i$ there exists a sequence
$v_k^i\in SBV^p(\Om_0;\Rm)$ such that $v_k^i=v^i$
a.e.\ in $\Om_0\setmeno\Om_B$,
$v_k^i\to v^i$ strongly
in $L^{q}(\Om_0;\Rm)$,
$\nabla v_k^i\to\nabla v^i$
strongly in $L^p(\Om_0;\Mmn)$,
and
$\hn((S(v_k^i)\setmeno\Ga_k^N)
\setmeno(S(v^i)\setmeno\Ga^N))\to0$.
Therefore there exists an increasing sequence of integers $k_i$ such
that
\begin{eqnarray}
&\|v_k^i- v^i\|_{q}<1/i \label{145}\\
&\|\nabla v_k^i- \nabla v^i\|_{p}<1/i \label{146}\\
&\hn((S(v_k^i)\setmeno\Ga_k^N)
\setmeno(S(v^i)\setmeno\Ga^N))<1/i \label{147}
\end{eqnarray}
for every $k\ge k_i$.

For $k_i\le k<k_{i+1}$ define $v_k:=v_k^i$ a.e.\ on $\Om_B$,
and $v_k:=v$
a.e.\ on $\Om\setmeno\Om_B$. Then $v_k\in
GSBV^p_{q}(\Om_0;\Rm)$ and condition (a) is
satisfied.

As $|\varphi_i(z)|\le C|z|$ for every $i$ and $z$, the sequence $v^i$
converges to $v$ in $L^{q}(\Om_0;\Rm)$. Since $\nabla
v^i=\nabla\varphi_i(v)\nabla v$, the sequence $\nabla v^i$
converges to $\nabla v$
strongly in $L^p(\Om_0;\Mmn)$. Therefore (b) and (c) follow from
(\ref{145}) and (\ref{146}).

For $k_i\le k<k_{i+1}$ let $x\in \partial\Om_B\setmeno (S(v_k^i)\cup 
S_i)$ such that (\ref{vbx}) holds.
As $x\in \partial\Om_B\setmeno S_i$, it is easy to deduce from
(\ref{vbx}) that
\begin{equation}\label{vbxi}
\tilde v_B(x)=\aplim\limits_{\scriptstyle y\to x,\ \scriptstyle
y\notin\overline\Om_B}v^i(y)=
\aplim\limits_{\scriptstyle y\to x,\ \scriptstyle
y\notin\overline\Om_B}v_k^i(y)\,.
\end{equation}
Since $x\notin S(v_k^i)$, the approximate limit of $v_k^i$ at $x$ 
exists, so that we must have
\begin{equation}\label{vbxik}
\tilde v_B(x):=
\aplim\limits_{\scriptstyle y\to x,\ \scriptstyle
y\in\overline\Om_B}v_k^i(y)\,.
\end{equation}
By the definition of $v_k$, from (\ref{vbx}) and (\ref{vbxik}) we
deduce that
$$
\tilde v_B(x):=
\aplim\limits_{\scriptstyle y\to x}v_k(y)\,,
$$
hence $x\notin S(v_k)$.

This shows that $\partial\Om_B \setmeno
(S(v_k^i)\cup S_i)\subsethn \partial\Om_B \setmeno S(v_k)$, which implies
$S(v_k)\cap \partial\Om_B \subsethn
(S(v_k^i)\cap\partial\Om_B)\cup S_i$.
Since
$S(v_k)\cap\Om_B\eqhn S(v_k^i)\cap\Om_B$, we conclude that
$S(v_k)\cap\overline\Om_B\subsethn (S(v_k^i)\cap\overline\Om_B)\cup S_i$.
Recalling that $S(v^i)\subsethn S(v)$, we
obtain
\begin{equation}\label{setmeno1}
[(S(v_k)\setmeno\Ga_k^N)
\setmeno(S(v)\setmeno\Ga^N)]\cap \overline\Om_B\subsethn
[((S(v_k^i)\setmeno\Ga_k^N)
\setmeno(S(v^i)\setmeno\Ga^N))\cup S_i]\cap \overline\Om_B\,.
\end{equation}

On the other hand,
$S(v_k)\setmeno\overline\Om_B\subsethn
S(v)\setmeno\overline\Om_B$ and
$\Ga_k^N\setmeno\overline\Om_B\eqhn 
\partial_N\Om\setmeno\overline\Om_B\eqhn \Ga^N\setmeno\overline\Om_B$ by 
Remark~\ref{bgk}-(b),
so that
\begin{equation}\label{setmeno2}
[(S(v_k)\setmeno\Ga_k^N)
\setmeno(S(v)\setmeno\Ga^N)]\setmeno \overline\Om_B\eqhn\emptyset\,.
\end{equation}
{}From (\ref{setmeno1}) and (\ref{setmeno2})
we deduce that
$$
(S(v_k)\setmeno\Ga_k^N)
\setmeno(S(v)\setmeno\Ga^N)\subsethn
((S(v_k^i)\setmeno\Ga_k^N)
\setmeno(S(v^i)\setmeno\Ga^N))\cup S_i
$$
for $k_i\le k<k_{i+1}$. Therefore (d) follows from (\ref{limSi}) and
(\ref{147}).
\end{proof}


\begin{remark}\label{rem15}
Condition (d) of Theorem~\ref{fl2},
together with (\ref{calK}) and (\ref{K2}), implies that
$$
\lim_{k\to\infty}\K((S(v_k)\setmeno\Ga_k^N)\setmeno
(S(v)\setmeno\Ga^N))=0\,,
$$
hence
$$
\limsup_{k\to\infty}\K(S(v_k)\setmeno\Ga_k^N)\le
\K(S(v)\setmeno\Ga^N)\,.
$$
\end{remark}

\subsection{Convergence of minima}\label{cm}
We begin by proving the existence Theorems~\ref{esistmin00} and~\ref{esistmin}.

\begin{proof}[Proof of Theorem~\ref{esistmin00}] Let us fix 
$t\in[0,T]$, let us extend $\psi(t)$ to a function $\psi_0(t)\in 
W^{1,p}(\Om_0;\Rm)\cap L^{q}(\Om_0;\Rm)$, and let $u_k$ be a 
minimizing sequence of problem (\ref{min00}). We extend also $u_k$ to 
$\Om_0$ by setting $u_k:=\psi_0(t)$ a.e.\ on $\Om_0\setmeno\Om$. The 
extended functions belong to $GSBV^p_q(\Om_0;\Rm)$ and satisfy $S(u_k)\subsethn \Ga(t)^N$ by 
Proposition~\ref{traces} and by the definition of $S(u_k)$.
Since $\Ec(t)(\psi(t))<+\infty$ by (\ref{bound1}), the infimum in 
(\ref{min00}) is finite.
By (\ref{coerc1}) there exists a constant $C\ge 0$ such that
\begin{equation}\label{bbound000}
\|\nabla u_k\|_{p,\Om_0}^p +
\|u_k\|_{q,\Om_0}^{q}\le C
\end{equation}
for every $k$. By Theorem~\ref{compgsbv}
there exists a subsequence,
still denoted $u_k$, which converges
weakly in $GSBV^p(\Om_0;\Rm)$ to a function $u$. Since
$u =\psi_0(t)$ a.e.\ on $\Om_0\setmeno\Om$ and $S(u)\subsethn\Ga(t)^N$ 
by Remark~\ref{Su},
we conclude that
$u\in AD(\psi(t),\Ga(t))$.

By (\ref{W}) in Theorem~\ref{semigsbv} we have
\begin{equation}\label{530000}
\W(\nabla u)\le\liminf_{k\to\infty}\,\W(\nabla u_k)\,.
\end{equation}

By Lemma~\ref{rem3} the functions $u_k$ and $u$ belong to
$W^{1,p}(\Om_S;\Rm)\cap L^{q}(\Om_S;\Rm)$. By (\ref{pq}) and
(\ref{bbound000}) the sequence $u_k$ is bounded in $W^{1,p}(\Om_S;\Rm)$,
so it converges to $u$ weakly in $W^{1,p}(\Om_S;\Rn)$.
{}From the compactness of the trace operator we deduce that $u_k$ 
converges to $u$ strongly in $L^{r}(\partial_S\Om;\Rm)$, and by the continuity of $\G(t)$ we get 
\begin{equation}\label{530200}
\G(t)(u)=\lim_{k\to\infty}\,\G(t)(u_k)\,.
\end{equation}

From (\ref{(2)}), (\ref{530000}), and (\ref{530200}) we obtain
$$
\Ec(t)(u)\le\liminf_{k\to\infty}\Ec(t)(u_k)\,.
$$
Since $u\in AD(\psi(t),\Ga(t))$ and $u_k$ is a minimizing sequence, we 
conclude that $u$ is a minimum point of~(\ref{min00}).
\end{proof}

\begin{proof}[Proof of Theorem~\ref{esistmin}]
Let us fix $t\in[0,T]$ and $\Ga_0\in\Rp(\overline\Om_B)$, and let 
$(u_k,\Ga_k)$ be a minimizing sequence of problem (\ref{min}) with $\Ga_0\subsethn\Ga_k\subsethn\overline\Om_B$ for every~$k$.
 We extend $u_k$ to $\Om_0$ by setting $u_k:=\psi_0(t)$ a.e.\ on 
$\Om_0\setmeno\Om$, where $\psi_0(t)$ is the function introduced 
at the beginning of the proof of Theorem~\ref{esistmin00}. 
Note that these extensions belong to $GSBV^p_q(\Om_0;\Rm)$ and satisfy $S(u_k)\subsethn 
\Ga_k^N$ by Proposition~\ref{traces} and by the definition of 
$S(u_k)$. Since $\E(t)(\psi(t),\Ga_0)<+\infty$ by (\ref{bound1}), 
the infimum in (\ref{min}) is finite.
By (\ref{K9}) and (\ref{coerc1}) there exists a constant $C\ge 0$ such that
\begin{equation}\label{bbound0}
\|\nabla u_k\|_{p,\Om_0}^p +
\|u_k\|_{q,\Om_0}^{q}+\hn(\Ga_k^N)\le C
\end{equation}
for every $k$. By Theorem~\ref{compgsbv}
there exists a subsequence,
still denoted by $u_k$, which converges
weakly in $GSBV^p(\Om_0;\Rm)$ to a function $u$ which satisfies
$u =\psi_0(t)$ a.e.\ on $\Om_0\setmeno\Om$.

By Theorem~\ref{compattezza} there exists a subsequence, still denoted 
$\Ga_k$, such that $\Ga_k^N$ $\sigmap$-converges to a set $\Ga^*$.
Since $\Ga_k^N\subsethn\overline\Om_B\cup\partial_N\Om$ and $\partial_N\Om$ is closed, we deduce, thanks to Remark~\ref{bgk}-(b), that $\Ga^*\subsethn\overline\Om_B\cup\partial_N\Om$.
By Proposition~\ref{gsbvrp} we have
$S(u)\subsethn\Ga^*$. Since $u=\psi_0(t)$ a.e.\ on $\Om_0\setmeno\Om$ we deduce that the traces of $u$ and $\psi(t)$ coincide $\hn$-a.e.\ on $\partial\Om\setmeno\Ga^*$. Let $\Ga:=\Ga^*\setmeno\partial_N\Om$. Then $\Ga\in\Rp(\overline\Om_B)$ and  
$u\in AD(\psi(t ),\Ga)$, because $\partial_D\Om\setmeno\Ga=\partial_D\Om\setmeno\Ga^*$.

Arguing as in the proof of Theorem~\ref{esistmin00} we obtain 
(\ref{530000}) and (\ref{530200}).
By (\ref{eqganga}) and (\ref{semiK}) we have also
\begin{equation}\label{5303}
\K(\Ga\cup\Ga_0)=\K(\Ga^*\cup\Ga_0)\le
\liminf_{k\to\infty}\,\K(\Ga_k^N\cup\Ga_0)=
\liminf_{k\to\infty}\,\K(\Ga_k)\,.
\end{equation}
By (\ref{(2)}), (\ref{530000}), (\ref{530200}), and (\ref{5303}) we have
$$
\E(t)(u,\Ga\cup\Ga_0)\le\liminf_{k\to\infty}\E(t)(u_k,\Ga_k)\,.
$$
Since $\Ga\cup\Ga_0\in\Rp(\overline \Om_B)$ and $u\in AD(\psi(t),\Ga\cup\Ga_0)$, our assumption on $(u_k,\Ga_k)$ implies that $(u,{\Ga\cup\Ga_0})$ is a minimum point 
of~(\ref{min}).
\end{proof}

We now prove the stability of minimizers with respect to the 
$\sigmap$-convergence.

\begin {theorem}\label{convmin} Let $t_k\in [0,T]$ and 
$\Ga_k\in\Rp(\overline\Om_B)$.
Assume that $t_k\to t_\infty$ and 
$\Ga_k^N$ $\sigmap$-con\-verges to $\Ga^*_\infty$, and define $\Ga_\infty:=\Ga^*_\infty\setmeno\partial_N\Om$. For every $k$ let $u_k\in
AD(\psi(t_k),\Ga_k)$ be a function such that
\begin{equation}\label{mink}
\E(t_k)(u_k,\Ga_k)\le \E(t_k)(v,\Ga)
\end{equation}
for every 
$\Ga\in\Rp(\overline\Om_B)$, with $\Ga_k\subsethn\Ga$,
and every $v\in AD(\psi(t_k),\Ga)$.
Then $\Ga_\infty\in \Rp(\overline\Om_B)$ and there exist a subsequence of $u_k$,
not relabelled, and a function
$u_\infty \in AD(\psi(t_\infty ),\Ga_\infty)$ such that
\begin{eqnarray}
& u_k\wto u_\infty \, \hbox{ weakly in }\, GSBV^p(\Om;\Rm)\,,
\label{123}\\
& u_k\wto u_\infty \, \hbox{ weakly in }\, L^{q}(\Om;\Rm)\,,
\label{123.5}\\
& u_k\to u_\infty \, \hbox{ strongly in }\, L^{\dot q}(\Om;\Rm)\,,
\label{124}\\
& u_k\to u_\infty \, \hbox{ strongly in }\, L^{r}(\partial_S\Om;\Rm)\,,
\label{125}\\
& \nabla u_k\wto \nabla u_\infty \, \hbox{ weakly in }\, L^p(\Om;\Mmn)\,.
\label{126}
\end{eqnarray}
Moreover
\begin{equation}\label{min*}
\E(t_\infty )(u_\infty ,\Ga_\infty )\le \E(t_\infty )(v,\Ga)
\end{equation}
for every
$\Ga\in\Rp(\overline\Om_B)$, with $\Ga_\infty
\subsethn\Ga$, and every $ v\in AD(\psi(t_\infty),\Ga)$.
Finally
\begin{eqnarray}
& \W(\nabla u_k)\to \W(\nabla u_\infty )\,,
\label{127}\\
& \F(t_k)(u_k)\to \F(t_\infty )(u_\infty )\,,
\label{128}\\
& \G(t_k)(u_k)\to \G(t_\infty )(u_\infty )\,.
\label{129}
\end{eqnarray}
\end{theorem}

\begin{proof}
Taking $\Ga:=\Ga_k$ and $v=\psi(t_k)$ in (\ref{mink}) we obtain
$$
\E(t_k)(u_k,\Ga_k)\le \E(t_k)(\psi(t_k),\Ga_k)\,.
$$
By (\ref{K9}) and (\ref{bound1}) we have
$$
\E(t_k)(\psi(t_k),\Ga_k)\le
\alpha_1^\E(\|\nabla \psi(t_k)\|_p^{p}+ \|\psi(t_k)\|_{q}^{q}+ \|\psi(t_k)\|_{r,\partial_S\Om}^{r})+
\kappa_2\hn(\Ga_k)
+\beta_1^\E \,,
$$
so that $\E(t_k)(u_k,\Ga_k)$ is bounded uniformly with respect to
$k$ (recall that $\hn(\Ga_k^N)$ is bounded by the definition of the $\sigmap$-convergence). By (\ref{K9}) and (\ref{coerc1}) there exists a constant $C\ge0$ such that
\begin{equation}\label{bbound}
\|\nabla u_k\|_p^p +
\|u_k\|_{q}^{q}+\hn(\Ga_k)\le C
\end{equation}
for every $k$.

Using an extension operator it is possible to construct an absolutely 
continuous function $t\mapsto\psi_0(t)$ from $[0,T]$ into 
$W^{1,p}(\Om_0;\Rm)\cap L^{q}(\Om_0;\Rm)$ such that 
$\psi_0(t)=\psi(t)$ a.e.\ in $\Om$ for every~$t$. Then we extend 
$u_k$ to $\Om_0$ by setting
$u_k:=\psi_0(t_k)$ a.e.\ on $\Om_0\setmeno\Om$.
The extended functions belong to $GSBV_q^p(\Om_0;\Rm)$ and satisfy $S(u_k)\subsethn \Ga_k^N$ by 
Proposition~\ref{traces} and by the definition of $S(u_k)$.
By (\ref{bbound}) and by Theorem~\ref{compgsbv}
there exists a subsequence,
still denoted by $u_k$, which converges to a function $u_\infty$
weakly in $GSBV^p(\Om_0;\Rm)$ and weakly in $L^{q}(\Om_0;\Rm)$. Since
$\psi_0(t_k)$
converges to $\psi_0(t_\infty )$ in $W^{1,p}(\Om_0;\Rm)$ we conclude that
$u_\infty =\psi_0(t_\infty )$ a.e.\ on $\Om_0\setmeno\Om$.
By Proposition~\ref{gsbvrp} we have
$S(u_\infty )\subsethn\Ga^*_\infty$,
hence the traces of  $u_\infty$ and $\psi(t_\infty)$ agree $\hn$-a.e.\ on $\partial\Om\setmeno\Ga^*_\infty$.
 
Since $\Ga_k^N\subsethn\overline\Om_B\cup\partial_N\Om$ and $\partial_N\Om$ is closed, we deduce that $\Ga^*_\infty\subsethn\overline\Om_B\cup\partial_N\Om$ and by Remark~\ref{bgk}-(a) we have $\hn(\Ga^*_\infty)<+\infty$, therefore $\Ga_\infty\in\Rp(\overline\Om_B)$ and 
 $u_\infty \in AD(\psi(t_\infty ),\Ga_\infty)$ because $\partial_D\Om\setmeno\Ga_\infty=\partial_D\Om\setmeno\Ga^*_\infty$.

Properties (\ref{123}) and (\ref{123.5}) have already been proved, 
and (\ref{126}) is a
consequence of (\ref{123}). Since $\dot q<q$, property (\ref{124}) follows from the fact that $u_k$ converges to $u_\infty$ pointwise a.e.\ on $\Om$, and that
$u_k$ is bounded in $L^{q}(\Om;\Rm)$ by (\ref{123.5}).

By Lemma~\ref{rem3} the functions $u_k$ and $u_\infty$ belong to
$W^{1,p}(\Om_S;\Rm)\cap L^{q}(\Om_S;\Rm)$. By (\ref{pq}) and
(\ref{bbound}) the sequence $u_k$ is bounded in $W^{1,p}(\Om_S;\Rm)$,
so it converges to $u_\infty$ weakly in $W^{1,p}(\Om_S;\Rm)$.
Property (\ref{125}) follows from the fact that the trace operator is
compact from $W^{1,p}(\Om_S;\Rm)$ into $L^{r}(\partial_S\Om;\Rm)$.

To prove (\ref{min*}), we fix
$\Ga\in\Rp(\overline\Om_B)$ with $\Ga_\infty
\subsethn\Ga$. Given $v\in AD(\psi(t_\infty),\Ga)$, we extend $v$ by setting $v:=\psi_0(t_\infty)$ a.e.\ on $\Om_0\setmeno\Om$, and define
$w:=v-\psi_0(t_\infty)$ a.e.\ on $\Om_0$.

Since $(\Ga^*_\infty)^N=\Ga_\infty^N$, by the Jump Transfer Theorem~\ref{fl2} there exists a sequence
$w_k\in GSBV^p_{q}(\Om_0;\Rm)$ such that
$w_k=w$ a.e.\ in $\Om_0\setmeno\Om_B$,
$w_k\to w$ strongly
in $L^{q}(\Om_0;\Rm)$,
$\nabla w_k\to\nabla w$
strongly in $L^p(\Om_0;\Mmn)$,
$\hn((S(w_k)\setmeno\Ga_k^N)
\setmeno(S(w)\setmeno\Ga_\infty^N))\to0$, and $S(w_k)\subsethn
\overline\Om_B\cup\partial_N\Om$.

Let $v_k:=w_k+\psi_0(t_k)$. Then $v_k\in GSBV^p(\Om_0;\Rm)$,
$v_k=\psi_0(t_k)$ a.e.\ in $\Om_0\setmeno\Om$,
$v_k\to v$ strongly
in $L^{q}(\Om_0;\Rm)$,
$\nabla v_k\to\nabla v$
strongly in $L^p(\Om_0;\Mmn)$,
$\hn(({S(v_k)\setmeno\Ga_k^N})
\setmeno({S(v)\setmeno\Ga_\infty^N}))\to0$,
and $S(v_k)\subsethn
\overline\Om_B\cup\partial_N\Om$. By Lemma~\ref{rem3} the functions 
$v_k$ and $v$ belong to $W^{1,p}(\Om_S;\Rm)$ and by (\ref{pq}) $v_k$ 
converges to $v$ strongly in $W^{1,p}(\Om_S;\Rm)$. We conclude by 
(\ref{trace}) that $v_k$ converges to $v$ strongly in 
$L^{r}(\partial_S\Om;\Rm)$.

Let $\Ga'_k:=\Ga_k\cup
(S(v_k)\setmeno \partial_N\Om)$. As $S(v_k)\setmeno\partial_N\Om\subsethn
\overline\Om_B\setmeno \partial_N\Om$, we have
$\Ga'_k\in \Rp(\overline\Om_B)$. 
Since $S(v_k)\cap\Om\subsethn\Ga'_k$, $S(v_k)\cap\partial_D\Om\subsethn\Ga'_k$, and $v_k=\psi_0(t_k)$ a.e.\ on $\Om_0\setmeno\Om$, by Proposition~\ref{traces} and by the definition of $S(v_k)$ we have $v_k\in AD(\psi(t_k),\Ga'_k)$.
By the 
minimality condition
(\ref{mink}) we have
\begin{eqnarray*}
&\W(\nabla u_k)+\K(\Ga_k)-\F(t_k)(u_k)- \G(t_k)(u_k) \le
\\
&\le \W(\nabla v_k)+\K(\Ga'_k)-\F(t_k)(v_k)- \G(t_k)(v_k)\,,
\end{eqnarray*}
which implies
\begin{eqnarray}
&\W(\nabla u_k)- \F(t_k)(u_k)- \G(t_k)(u_k) \le \label{150}\\
&\le \W(\nabla v_k)+\K(S(v_k)\setmeno\Ga_k^N)
-\F(t_k)(v_k)- \G(t_k)(v_k)\,.\nonumber
\end{eqnarray}

By (\ref{W}) in Theorem~\ref{semigsbv} we have
\begin{equation}\label{semiW}
\W(\nabla u_\infty)\le\liminf_{k\to\infty} \, \W(\nabla u_k)\,.
\end{equation}
Since $\nabla v_k$ converges to $\nabla v$ strongly in 
$L^p(\Om_0;\Mmn)$, by the continuity of $\W$ we have
\begin{equation}\label{limWv}
\W(\nabla v) = \lim_{k\to\infty} \, \W(\nabla v_k)\,.
\end{equation}

By (\ref{F8}) and (\ref{F11}) we have
\begin{eqnarray}
&\displaystyle
|\F(t_k)(u_k)- \F(t_\infty)(u_k)|\le \int_{t_\infty}^{t_k} |\dot 
\F(s)(u_k)|\,ds \le \int_{t_\infty}^{t_k} 
(\alpha_3^\F(s)\|u_k\|_{\dot q}^{\dot q}
+ \beta_3^\F(s))\,ds\,,\ \ \label{97}\\
&\displaystyle
|\F(t_k)(v_k)- \F(t_\infty)(v_k)|\le \int_{t_\infty}^{t_k} |\dot 
\F(s)(v_k)|\,ds \le \int_{t_\infty}^{t_k} 
(\alpha_3^\F(s)\|v_k\|_{\dot q}^{\dot q}
+ \beta_3^\F(s))\,ds\,.\ \ \label{98}
\end{eqnarray}
As $u_k$ converges to $u_\infty$ strongly in $L^{\dot q}(\Om;\Rm)$ by 
(\ref{124}), using (\ref{(2)}) and (\ref{97}) we obtain
\begin{equation}\label{limF}
\F(t_\infty)(u_\infty)\ge \limsup_{k\to\infty} \, \F(t_\infty)(u_k) =
\limsup_{k\to\infty} \, \F(t_k)(u_k) \,.
\end{equation}
As $\F(t_\infty)$ is continuous in $L^{q}(\Om;\Rm)$ and
$v_k$ converges to $v_\infty$ strongly in $L^{q}(\Om;\Rm)$, 
using (\ref{98}) we obtain
\begin{equation}\label{limFv}
\F(t_\infty)(v)= \lim_{k\to\infty} \, \F(t_\infty)(v_k) =
\lim_{k\to\infty} \, \F(t_k)(v_k) \,.
\end{equation}

In the same way we prove that
\begin{eqnarray}
&\displaystyle
\G(t_\infty)(u_\infty)\ge \limsup_{k\to\infty} \, \G(t_\infty)(u_k) =
\limsup_{k\to\infty} \, \G(t_k)(u_k) \,,\label{limG}\\
&\displaystyle
\G(t_\infty)(v)= \lim_{k\to\infty} \, \G(t_\infty)(v_k) =
\lim_{k\to\infty} \, \G(t_k)(v_k) \,.\label{limGv}
\end{eqnarray}

{}From (\ref{semiW}), (\ref{limF}), and (\ref{limG}) we obtain
\begin{eqnarray}
&\W(\nabla u_\infty)- \F(t_\infty)(u_\infty)- \G(t_\infty)(u_\infty) 
\le \label{151}\\
&\displaystyle
\le\liminf_{k\to\infty}\,\big[ \W(\nabla u_k)-\F(t_k)(u_k)- 
\G(t_k)(u_k) \big]\,.\nonumber
\end{eqnarray}
{}From Remark~\ref{rem15}, and from (\ref{limWv}), (\ref{limFv}), and 
(\ref{limGv}) we obtain
\begin{eqnarray}
&\displaystyle
\limsup_{k\to\infty}\,\big[ \W(\nabla v_k)+\K(S(v_k)\setmeno\Ga_k^N)
- \F(t_k)(v_k)- \G(t_k)(v_k) \big] \le \label{152}\\
&\le \W(\nabla v)+\K(S(v)\setmeno\Ga_\infty^N)
- \F(t_\infty)(v)- \G(t_\infty)(v)\,.\nonumber
\end{eqnarray}
By (\ref{150}), (\ref{151}), and (\ref{152}) we have
\begin{eqnarray}
&\W(\nabla u_\infty)- \F(t_\infty)(u_\infty)- \G(t_\infty)(u_\infty) 
\le \label{153}\\
&\le \W(\nabla v)+\K(S(v)\setmeno\Ga_\infty^N)
- \F(t_\infty)(v)- \G(t_\infty)(v)\,.\nonumber
\end{eqnarray}
As $S(v)\subsethn \Ga^N$, we have $\K(S(v)\setmeno\Ga_\infty^N)\leq
\K(\Ga^N\setmeno\Ga_\infty^N)=\K(\Ga\setmeno\Ga_\infty)$, so that inequality (\ref{153}) implies~(\ref{min*}).

Taking $v:=u_\infty$, from (\ref{150}) and (\ref{152}) we obtain
\begin{eqnarray}
&\displaystyle
\limsup_{k\to\infty}\,\big[ \W(\nabla u_k)- \F(t_k)(u_k)- 
\G(t_k)(u_k) \big] \le \nonumber\\
\nonumber
&\le \W(\nabla u_\infty)- \F(t_\infty)(u_\infty)- \G(t_\infty)(u_\infty)\,.
\end{eqnarray}
This inequality, together with (\ref{semiW}), (\ref{limF}), and 
(\ref{limG}), gives (\ref{127}),
(\ref{128}), and (\ref{129}). 
\end{proof}

\subsection{Convergence of Riemann sums}\label{cRiem}
To prove Theorem~\ref{main} 
we need to approximate the time integrals of $\dot\F(t)(u)$ and 
$\dot\G(t)(u)$ by Riemann sums, uniformly when $u$ varies in a 
compact set. As usual, if $f$ is a measurable function defined a.e.\ on $[0,T]$, the same symbol $f$ denotes its extension by $0$ to $[0,T]$.
We begin by the following remark.

\begin{remark}\label{riemann3}
 Let ${\mathcal V}$ be a countable dense subset of $W^{1,p}(\Om;\Rm)\cap L^{\dot q}(\Om;\Rm)$. 
By Remark~\ref{riemann2} for every $t\in(0,T]$ there exists a sequence of subdivisions
$(\ski)_{0\le i\le i_k}$, with
\begin{eqnarray}
&0=s_k^0<s_k^1<\cdots<s_k^{i_k-1}<s_k^{i_k}=t\,,
\label{subdivs}\\
&\displaystyle
\lim_{k\to\infty}\, \max_{1\le i\le i_k} (\ski-\skim)= 0\,,
\label{fines}
\end{eqnarray}
such that
\begin{eqnarray}
&\displaystyle
\lim_{k\to\infty}\, \sum_{i=1}^{i_k}\Big | (\ski-\skim)
\dot \F(\ski)(u) -
\int_{\skim}^{\ski} \dot \F(s)(u)\,ds \,\Big| = 0
\,,\label{riemsums}\\
&\displaystyle
\lim_{k\to\infty}\, \sum_{i=1}^{i_k} \Big | (\ski-
\skim)
\dot \G(\ski)(u) -
\int_{\skim}^{\ski} \dot \G(s)(u)\,ds \,\Big| = 0
\label{riemsumsg}
\end{eqnarray}
for every $u\in{\mathcal V}$.
\end{remark}

\begin{lemma}\label{riemann4}
Let $t\in (0,T]$, let ${\mathcal V}$ be a countable dense subset of $W^{1,p}(\Om;\Rm)\cap L^{\dot q}(\Om;\Rm)$, and
let $(\ski)_{0\le i\le i_k}$ be a sequence of subdivisions satisfying
(\ref{subdivs})--(\ref{riemsumsg})
for every $u\in{\mathcal V}$.
Assume that
\begin{equation} \label{riemalpha}
\lim_{k\to\infty}\, \sum_{i=1}^{i_k} \Big | (\ski-\skim)
\varphi(\ski) -
\int_{\skim}^{\ski}\varphi(s)\,ds \,\Big| =0 
\end{equation}
whenever $\varphi$ is any one of the four functions $\alpha_4^\F$, $\beta_4^\F$,
$\alpha_4^\G$, and $\beta_4^\G$ which appear
in (\ref{F11.5}) and (\ref{G11.5}). Then
\begin{equation}\label{riemsumk}
\lim_{k\to\infty}\, \sum_{i=1}^{i_k}\, \sup_{u\in{\mathcal H} } \,
\Big | (\ski-\skim) \dot \F(\ski)(u) -
\int_{\skim}^{\ski} \dot \F(s)(u)\,ds \,\Big| = 0
\end{equation}
for every compact subset ${\mathcal H} $ of $L^{\dot q}(\Om;\Rm)$, and
\begin{equation}
\lim_{k\to\infty}\, \sum_{i=1}^{i_k}\, \sup_{u\in{\mathcal H} } \,
\Big | (\ski-\skim) \dot \G(\ski)(u) -
\int_{\skim}^{\ski} \dot \G(s)(u)\,ds \,\Big| = 0 \label{riemsumkg}
\end{equation}
for every compact subset ${\mathcal H} $ of $L^{r}(\partial_S\Om;\Rm)$.
\end{lemma}

\begin{proof}
Let us fix
a compact subset ${\mathcal H} $ of $L^{\dot q}(\Om;\Rm)$ and let
$$
M:=\max\{\|u\|_{\dot q}: u\in {\mathcal H} \}+1\,.
$$
For every $\e\in(0,1)$ there exists a finite number
$u_1,\,\ldots,\,u_h\in{\mathcal V}$, with
$\|u_j\|_{\dot q}< M$, such that for every $u\in {\mathcal H} $
there exists $j$ with $\|u-u_j\|_{\dot q}<\e$. Therefore
for every $u\in{\mathcal H}$ and every $s\in[0,t]$
there exists $v_j(s)\in L^{\dot q}(\Om;\Rm)$, with $\|v_j(s)\|_{\dot q}< M$,
such that
$$
| \dot \F(s)(u)-\dot \F(s)(u_j)| \leq
|\langle \partial \dot\F(s)(v_j(s)),u_j-u\rangle|\,.
$$
By (\ref{F11.5}) this implies that
\begin{equation}\label{700}
| \dot \F(s)(u)-\dot \F(s)(u_j)|\le \e\,(\alpha_4^\F(s) M^{\dot q-1} +
\beta_4^\F(s))
\end{equation}
for every $s\in[0,t]$.
Consequently
\begin{eqnarray*}
&\displaystyle
\Big | (\ski-\skim)
\dot \F(\ski)(u) -
\int_{\skim}^{\ski} \dot \F(s)(u)\,ds \,\Big| \le\\
&\displaystyle
\le \Big | (\ski-\skim)
\dot \F(\ski)(u_j) -
\int_{\skim}^{\ski} \dot \F(s)(u_j)\,ds \,\Big| +{} \\
&\displaystyle
{}+ \e\, (\ski-\skim) \,(\alpha_4^\F(\ski) M^{\dot q-1} +
\beta_4^\F(\ski)) +{}\\
&\displaystyle
{}+\e \int_{\skim}^{\ski} (\alpha_4^\F(s) M^{\dot q-1} +
\beta_4^\F(s))\,ds\,.
\end{eqnarray*}
This yields
\begin{eqnarray*}
&\displaystyle
\sup_{u\in{\mathcal H}}\,
\Big | (\ski-\skim)
\dot \F(\ski)(u) -
\int_{\skim}^{\ski} \dot \F(s)(u)\,ds \,\Big| \le \\
&\displaystyle
\le \sum_{j=1}^h
\Big | (\ski-\skim)
\dot \F(\ski)(u_j) -
\int_{\skim}^{\ski} \dot \F(s)(u_j)\,ds \,\Big| +{} \\
&\displaystyle
{}+ \e\, (\ski-\skim) \,(\alpha_4^\F(\ski) M^{\dot q-1} +
\beta_4^\F(\ski)) +{} \\
&\displaystyle
{}+ \e \int_{\skim}^{\ski} (\alpha_4^\F(s) M^{\dot q-1} +
\beta_4^\F(s))\,ds\,.
\end{eqnarray*}
Therefore
\begin{eqnarray*}
&\displaystyle
\sum_{i=1}^{i_k} \,\sup_{u\in{\mathcal H}}\,
\Big | (\ski-\skim)
\dot \F(\ski)(u) -
\int_{\skim}^{\ski} \dot \F(s)(u)\,ds \,\Big| \le \\
&\displaystyle
\le \sum_{j=1}^h \sum_{i=1}^{i_k}
\Big | (\ski-\skim)
\dot \F(\ski)(u_j) -
\int_{\skim}^{\ski} \dot \F(s)(u_j)\,ds \,\Big| +{} \\
&\displaystyle
{}+ \e\, \sum_{i=1}^{i_k} (\ski-\skim) \,(\alpha_4^\F(\ski) M^{\dot q-1}
+
\beta_4^\F(\ski)) +{} \\
&\displaystyle
{}+ \e \int_{0}^{t} (\alpha_4^\F(s) M^{\dot q-1} +
\beta_4^\F(s))\,ds\,.
\end{eqnarray*}
By (\ref{riemsums}) and (\ref{riemalpha}) we have
\begin{eqnarray*}
&\displaystyle
\limsup_{k\to\infty} \,
\sum_{i=1}^{i_k} \,\sup_{u\in{\mathcal H}}\,
\Big | (\ski-\skim)
\dot \F(\ski)(u) -
\int_{\skim}^{\ski} \dot \F(s)(u)\,ds \,\Big| \le \\
&\displaystyle
\le 2 \e \int_{0}^{t} (\alpha_4^\F(s) M^{\dot q-1} +
\beta_4^\F(s))\,ds\,,
\end{eqnarray*}
and taking the limit as $\e\to 0$ we obtain~(\ref{riemsumk}).

The proof of (\ref{riemsumkg}) is similar.
\end{proof}

\end{section}

\begin{section}{THE DISCRETE-TIME PROBLEMS}\label{discrete}

Theorem~\ref{main} will be proved by time discretization. 
In this section we study the discrete-time problems and prove a fundamental energy estimate, which will be crucial in the proof of the nondissipativity condition for the solution of the continuous-time quasistatic evolution (condition~(c) in Subsection~\ref{quasistatic}). 

Let us fix a sequence of subdivisions $(\tki)_{0\le i\le k}$
of the interval $[0,T]$, with
\begin{eqnarray}
& 0=t_k^0<t_k^1<\cdots<t_k^{k-1}<t_k^{k}=T\,,
\label{subdiv5}\\
&\displaystyle
\lim_{k\to\infty}\,
\max_{1\le i\le k} (\tki-\tkim)= 0\,.
\label{fine5}
\end{eqnarray}
For $i=1,\ldots,k$ we set
\begin{equation}\label{ki}
\F\ki:=\F(\tki)\,, \qquad \G\ki:=\G(\tki)\,,\qquad 
\psiki:=\psi(\tki)\,, \qquad \E\ki:=\E(\tki)\,.
\end{equation}

Let $(u_0,\Ga_0)$ be an initial configuration that satisfies the 
minimality property (\ref{mininit}). For every $k$ we define $\uki$ 
and $\gki$ by induction.
We set $(u_k^0, \Ga_k^0):=(u_0,\Ga_0)$, and for $i=1,\ldots,k$ we 
define $(\uki,\gki)$ as a solution of the problem
\begin{equation}\label{2110}
\min\, \{\E\ki(u,\Ga): \Ga\in \Rp(\overline\Om_B),\ \gkim\subsethn\Ga,\ u\in AD(\psiki,\Ga)\}\,.
\end{equation}
The existence of a solution to this problem is proved in 
Theorem~\ref{esistmin}.
Note that
\begin{equation}\label{6000}
\E\ki(\uki,\gki)\le\E\ki(u,\Ga)
\end{equation}
for every $\Ga\in \Rp(\overline\Om_B)$, with $\gki\subsethn\Ga$ and every $u\in AD(\psiki,\Ga)$.

Since $\psiki\in AD(\psiki,\gkim)$, by (\ref{2110}) we have
$\Ec(\tki)(\uki)+\K(\gki)\le \Ec(\tki)(\psiki)+\K(\gkim)$, 
which gives
$\Ec(\tki)(\uki)\le \Ec(\tki)(\psiki)$.
By (\ref{bound1}) this implies
$$
\Ec(\tki)(\uki)\le 
\alpha_1^\E(\|\nabla\psiki\|_p^{p}+\|\psiki\|_{q}^{q}+\|\psiki\|_{r,\partial_S\Om}^r)+\beta_1^\E\,,
$$
so that by (\ref{coerc1}) there 
exists a constant $C>0$ such that
\begin{equation}\label{6001}
\|\nabla\uki\|_p\le C\,, \quad \|\uki\|_{\dot q}\le C\,,\quad
 \|\uki\|_{q}\le C
\end{equation}
for every $k$ and $i$.
By Lemma~\ref{rem3} the functions $\uki$ belong to 
$W^{1,p}(\Om_S;\Rm)\cap L^{q}(\Om_S;\Rm)$. Therefore 
(\ref{tracepq}) implies that, if we change the constant, we may assume also that
\begin{equation}\label{6002}
\|\uki\|_{r,\partial_S\Om }\le C
\end{equation}
for every $k$ and $i$.

For every $t\in[0,T]$ we define
\begin{equation}\label{ukt}
\begin{array}{ccc}
\vphantom{\displaystyle\min_A }
\tau_k(t):=\tki\,,
&\uk(t):=\uki\,,
&\Ga_k(t):=\gki\,,
\\
\F_k(t):=\F\ki:=\F(\tki)\,,
&\G_k(t):=\G\ki:=\G(\tki)\,,
&\E_k(t):=\E\ki:=\E(\tki)\,,
\end{array}
\end{equation}
where $i$ is the largest integer such that $\tki\le t$. Note that $\uk(t)=\uk(\tau_k(t))$
and $\Ga_k(t)=\Ga_k(\tau_k(t))$.
By (\ref{6001}) and (\ref{6002}) we have
\begin{equation}\label{2170}
\|\nabla\uk(t)\|_p\le C\,, \quad \|\uk(t)\|_{\dot q}\le C\,,\quad
\|\uk(t)\|_{q}\le C\,,\quad \|\uk(t)\|_{r,\partial_S\Om }\le C\,,
\end{equation}
for every $k$ and every $t\in[0,T]$.

We introduce now a sequence of functions which play an important role in our estimates. For a.e.\ $t\in[0,T]$ we set
\begin{equation}\label{thetak}
\begin{array}{c} \vphantom{\displaystyle \min_A}
\theta_k(t):= \langle \partial\W(\nabla\uk(t)), \nabla \dot\psi(t) \rangle
- \langle \partial\F_k(t)(\uk(t)), \dot\psi(t) \rangle - {}\\
{} - \dot \F(t)(\uk(t))- \langle \partial\G_k(t)(\uk(t)), 
\dot\psi(t) \rangle -
\dot \G(t)(\uk(t)) \,.
\end{array}
\end{equation}

The main result of this section is the energy estimate given by the following lemma.
\begin{lemma}\label{energy estimate}
There exists a sequence $R_k\to 0$ such that
\begin{equation}\label{1999'}
\E_k(t)(\uk(t),\Ga_k(t))\le \E(0)(u_0,\Ga_0)+\int_0^{\tau_k(t)} \theta_k(s)\,ds + R_k
\end{equation}
for every $k$ and every $t\in[0,T]$.
\end{lemma}

\begin{proof} We have to prove that there exists a sequence $R_k\to 0$ 
such that
\begin{equation}\label{1999}
\E\ki(\uki,\gki)\le \E(0)(u_0,\Ga_0)+\int_0^{\tki} \theta_k(s)\,ds + R_k
\end{equation}
for every $k$ and for every $i=1,\ldots,k$.

Let us fix $j$ and $k$ with $1\le j\le k$.
Since $\ukjm+\psikj-\psikjm \in AD(\psikj,\gkjm)$, by the minimality 
condition (\ref{2110}) we have
\begin{equation}\label{2001}
\E\kj(\ukj,\gkj)\le\E\kj(\ukjm+\psikj-\psikjm,\gkjm)
\end{equation}
We now estimate $\E\kj(\ukjm+\psikj-\psikjm,\gkjm)$ in terms of 
$\E\kjm(\ukjm,\gkjm)$.

Let us consider first $\W(\ukjm+\psikj-\psikjm)$. There exists a 
constant $\rhokj\in[0,1]$ such that
\begin{equation}\label{2002}
\begin{array}{c}
 \W(\nabla\ukjm+\nabla\psikj-\nabla\psikjm)- \W(\nabla\ukjm)=\\
= \langle \partial\W(\nabla\ukjm+ 
\rhokj(\nabla\psikj-\nabla\psikjm)),\nabla \psikj-\nabla\psikjm 
\rangle \,,
\end{array}
\end{equation}
where $\langle\cdot,\cdot\rangle$ denotes the duality pairing between
$L^{p'}(\Om;\Mmn)$ and $L^{p}(\Om;\Mmn)$.
Let us consider the piecewise constant function $\Psik\in 
L^\infty([0,T];L^p(\Om;\Mmn))$ defined by
\begin{equation}\label{2003}
\Psik(s):=\rhokj(\nabla\psikj-\nabla\psikjm)=
\rhokj\int_{\tkjm}^{\tkj}\nabla \dot\psi(\tau)\,d\tau\qquad \hbox{for }\, 
\tkjm\le s < \tkj\,.
\end{equation}
Since $s\mapsto \nabla \dot\psi(s)$ belongs to 
$L^1([0,T];L^p(\Om;\Mmn))$, by the absolute continuity of the 
integral we have that
\begin{equation}\label{2003.5}
\|\Psik(s)\|_p\to 0\,\hbox{ uniformly with respect to }\, s\in[0,T]\,.
\end{equation}
Therefore by (\ref{2170}) we may assume that 
$\|\nabla u_k(s)+\Psik(s)\|_p\le C+1$ for every $s\in[0,T]$.
{}From (\ref{2002}) and (\ref{2003}) we obtain
\begin{equation}\label{2004}
\begin{array}{c}
\vphantom{ \int_{\tkjm}}
\W(\nabla\ukjm+\nabla\psikj-\nabla\psikjm)- \W(\nabla\ukjm)=\\
\displaystyle = \int_{\tkjm}^{\tkj} \langle \partial\W(\nabla\uk(s)+ 
\Psik(s)),\nabla \dot\psi(s) \rangle \,ds\,. 
\end{array}
\end{equation}

Let us consider now $\F\kj(\ukjm+\psikj-\psikjm)$. There exists a 
constant $\sigmakj\in[0,1]$ such that
\begin{equation}\label{2012}
\F\kj(\ukjm+\psikj-\psikjm)-\F\kj(\ukjm)=\langle \partial\F\kj(\ukjm+ 
\sigmakj(\psikj-\psikjm)), \psikj-\psikjm \rangle \,,
\end{equation}
where $\langle\cdot,\cdot\rangle$ denotes now the duality pairing between
$L^{q'}(\Om;\Rm)$ and $L^{q}(\Om;\Rm)$.
Let us consider the piecewise constant function $\phik\in 
L^\infty([0,T];L^{q}(\Om;\Rm))$ defined by
\begin{equation}\label{2013}
\phik(s):=\sigmakj(\psikj-\psikjm)=
\sigmakj\int_{\tkjm}^{\tkj}\dot\psi(\tau)\,d\tau\qquad \hbox{for }\, 
\tkjm\le s < \tkj\,.
\end{equation}
Since $s\mapsto \dot\psi(s)$ belongs to 
$L^1([0,T];L^{q}(\Om;\Rm))$, by the absolute continuity of the 
integral we have that
\begin{equation}\label{2013.5}
\|\phik(s)\|_{q}\to 0\,\hbox{ uniformly with respect to }\, s\in[0,T]\,.
\end{equation}
Therefore by (\ref{2170}) we may assume that 
$\|\uk(s)+\phik(s)\|_{\dot q}\le C+1$ for every $s\in[0,T]$.

{}From (\ref{2012}) and (\ref{2013}) we obtain
\begin{equation}\label{2014.0}
\begin{array}{c}
\vphantom{ \int_{\tkjm}}
\F\kj(\ukjm+\psikj-\psikjm)-\F\kj(\ukjm)=\\
\displaystyle = \int_{\tkjm}^{\tkj} \langle \partial\F\kj(\uk(s)+ \phik(s)),
\dot\psi(s) \rangle \,ds\,.
\end{array}
\end{equation}
By (\ref{F8.5}) for every $s\in[\tkjm,\tkj)$ we have
\begin{eqnarray*}
& \langle \partial\F\kj(\uk(s)+ \phik(s)), \dot\psi(s) \rangle -
\langle \partial\F(s)(\uk(s)+ \phik(s)), \dot\psi(s) \rangle=\\
&\displaystyle = \int_{s}^{\tkj} \langle \partial\dot \F(\sigma)(\uk(s)+ 
\phik(s)), \dot\psi(s) \rangle\,d\sigma\,,\\
&\displaystyle \langle \partial\F(s)(\uk(s)), \dot\psi(s) \rangle-
\langle \partial\F\kjm(\uk(s)), \dot\psi(s) \rangle=
\int_{\tkjm}^{s} \langle \partial\dot \F(\sigma)(\uk(s)), \dot\psi(s)\rangle\,d\sigma\,,
\end{eqnarray*}
By (\ref{F11.5}) this implies that for every $s\in[\tkjm,\tkj)$
\begin{eqnarray}
& | \langle \partial\F\kj(\uk(s)+ \phik(s)), \dot\psi(s) \rangle -
\langle \partial\F(s)(\uk(s)+ \phik(s)), \dot\psi(s) \rangle| \le
\gamma^\F_k \|\dot\psi(s)\|_{\dot q}\,, \label{2014.5}\\
& | \langle \partial\F(s)(\uk(s)), \dot\psi(s) \rangle -
\langle \partial\F_k(s)(\uk(s)), \dot\psi(s) \rangle| \le
\gamma^\F_k \|\dot\psi(s)\|_{\dot q}\,, \label{2014.6}
\end{eqnarray}
where
$$ 
\gamma^\F_k:=\max_{1\le i\le k} \int_{\tkim}^{\tki} (\alpha_4^\F(s) 
(C+1)^{\dot q-1} +
\beta_4^\F(s))\,ds\,.
$$ 
By the absolute continuity of the integral we have
\begin{equation}\label{2014.8}
\gamma^\F_k\to 0\,,
\end{equation}
and by (\ref{2014.0}) and (\ref{2014.5}) we have
\begin{eqnarray}
&\F\kj(\ukjm+\psikj-\psikjm)- \F\kj(\ukjm)\ge\nonumber\\
&\displaystyle \ge \int_{\tkjm}^{\tkj} \langle \partial\F(s)(\uk(s)+ \phik(s)),
\dot\psi(s) \rangle \,ds - \gamma^\F_k 
\int_{\tkjm}^{\tkj}\|\dot\psi(s)\|_{\dot q} \,ds\,.\label{2014} 
\end{eqnarray}
Since by (\ref{F8})
\begin{equation}\label{2015}
\F\kj(\ukjm)- \F\kjm(\ukjm)=
\int_{\tkjm}^{\tkj} \dot \F(s)(\ukjm) \,ds=
\int_{\tkjm}^{\tkj} \dot \F(s)(\uk(s)) \,ds\,,
\end{equation}
{}from (\ref{2014}) and (\ref{2015}) we obtain
\begin{eqnarray}
&\displaystyle \vphantom{\int_{\tkjm}^{\tkj}}
\F\kj(\ukjm+\psikj-\psikjm)- \F\kjm(\ukjm)\ge \nonumber\\
&\displaystyle \ge \int_{\tkjm}^{\tkj} \langle \partial\F(s)(\uk(s)+ 
\phik(s)),
\dot\psi(s) \rangle \,ds + \label{2016}\\
&\displaystyle +\int_{\tkjm}^{\tkj} \dot \F(s)(\uk(s)) \,ds
- \gamma^\F_k \int_{\tkjm}^{\tkj} \|\dot\psi(s)\|_{\dot q} \,ds\,. \nonumber
\end{eqnarray}

Finally, let us consider $\G\kj(\ukjm+\psikj-\psikjm)$. The same arguments used for 
$\F\kj$ show that there exist a sequence $\gamma^\G_k$ of real numbers, with
\begin{equation}\label{2024.8}
\gamma^\G_k\to 0\,,
\end{equation}
and a sequence $\zetak\in 
L^\infty([0,T];L^{r}(\partial_S\Om;\Rm))$, with
\begin{equation}\label{2023.5}
\|\zetak(s)\|_{r,\partial_S\Om}\to 0\,\hbox{ uniformly with respect 
to }\, s\in[0,T]\,,
\end{equation}
such that
\begin{equation} \label{2024.6}
\displaystyle \vphantom{\int_{\tkjm}^{\tkj}}
 | \langle \partial\G(s)(\uk(s)), \dot\psi(s) \rangle -
\langle \partial\G_k(s)(\uk(s)), \dot\psi(s) \rangle| \le
\gamma^\G_k \|\dot\psi(s)\|_{r,\partial_S\Om}
\end{equation}
and
\begin{eqnarray}
&\displaystyle \vphantom{\int_{\tkjm}^{\tkj}}
\G\kj(\ukjm+\psikj-\psikjm)- \G\kjm(\ukjm)\ge \nonumber \\
&\displaystyle \ge \int_{\tkjm}^{\tkj} \langle \partial\G(s)(\uk(s)+ 
\zetak(s)),
\dot\psi(s) \rangle \,ds +\label{2026}\\
&\displaystyle +\int_{\tkjm}^{\tkj} \dot \G(s)(\uk(s)) \,ds
-\gamma^\G_k
\int_{\tkjm}^{\tkj} \|\dot\psi(s)\|_{r,\partial_S\Om} \, ds\,. \nonumber
\end{eqnarray}
By (\ref{2001}), (\ref{2004}), (\ref{2016}), and (\ref{2026}) we have
\begin{eqnarray}
&\displaystyle \vphantom{\int_{\tkjm}^{\tkj}}
\E\kj(\ukj,\gkj) - \E\kjm(\ukjm,\gkjm) \le \nonumber\\
&\displaystyle
\le \int_{\tkjm}^{\tkj} \langle \partial\W(\nabla\uk(s)+\Psik(s)),
\nabla \dot\psi(s) \rangle \,ds - \nonumber\\
&\displaystyle - \int_{\tkjm}^{\tkj} \langle \partial\F(s)(\uk(s)+\phik(s)),
\dot\psi(s) \rangle \,ds - \int_{\tkjm}^{\tkj} \dot \F(s)(\uk(s)) 
\,ds- \label{2030}\\
&\displaystyle -\int_{\tkjm}^{\tkj} \langle \partial\G(s)(\uk(s)+ \zetak(s)),
\dot\psi(s) \rangle \,ds- \int_{\tkjm}^{\tkj} \dot \G(s)(\uk(s)) \,ds +
\nonumber\\
&\displaystyle +\gamma^\F_k
\int_{\tkjm}^{\tkj} \|\dot\psi(s)\|_{\dot q} \, ds+
\gamma^\G_k
\int_{\tkjm}^{\tkj} \|\dot\psi(s)\|_{r,\partial_S\Om} \, ds\,.
\nonumber
\end{eqnarray}

Let us fix now $i$ with $1\le i\le k$. By summing for $j=1,\ldots,i$ we obtain
\begin{eqnarray}
&\displaystyle \vphantom{\int_{\tkim}^{\tki}}
\E\ki(\uki,\gki) - \E(0)(u_0,\Ga_0) \le \nonumber\\
&\displaystyle
\le \int_{0}^{\tki} \langle \partial\W(\nabla\uk(s)+\Psik(s)),
\nabla \dot\psi(s) \rangle \,ds - \nonumber\\
&\displaystyle - \int_{0}^{\tki} \langle \partial\F(s)(\uk(s)+\phik(s)),
\dot\psi(s) \rangle \,ds - \int_{0}^{\tki} \dot \F(s)(\uk(s)) \,ds- 
\label{2040}\\
&\displaystyle - \int_{0}^{\tki} \langle \partial\G(s)(\uk(s)+ \zetak(s)),
\dot\psi(s) \rangle \,ds - \int_{0}^{\tki} \dot \G(s)(\uk(s)) \,ds +
\nonumber\\
&\displaystyle{} +\gamma^\F_k
\int_{0}^{\tki} \|\dot\psi(s)\|_{\dot q} \, ds+
\gamma^\G_k
 \int_{0}^{\tki} \|\dot\psi(s)\|_{r,\partial_S\Om} \, ds\,.
\nonumber
\end{eqnarray}

Using Lemma~\ref{perturbation}, from (\ref{alpha2W}), (\ref{W7.5}), (\ref{2170}), and (\ref{2003.5}) 
for a.e.\ $s\in[0,T]$ we obtain
$$
\langle \partial\W(\nabla\uk(s)+\Psik(s)), \nabla \dot\psi(s) \rangle -
\langle \partial\W(\nabla\uk(s)), \nabla \dot\psi(s) \rangle \to 0\,,
$$
and by (\ref{W10}) and (\ref{2170}) we have
\begin{eqnarray*}
&|\langle \partial\W(\nabla\uk(s)+\Psik(s)), \nabla \dot\psi(s) \rangle -
\langle \partial\W(\nabla\uk(s)), \nabla \dot\psi(s) \rangle| \le
\\
&\le 2 (\alpha_2^\W (C+1)^{p-1} + \beta_2^\W)\, \|\nabla \dot\psi(s)\|_p\,.
\end{eqnarray*}
As $s\mapsto \nabla \dot\psi(s)$ belongs to 
$L^1([0,T];L^p(\Om;\Mmn))$, we 
deduce that
\begin{equation}\label{2202}
\int_0^T |\langle \partial\W(\nabla\uk(s)+\Psik(s)), \nabla 
\dot\psi(s) \rangle -
\langle \partial\W(\nabla\uk(s)), \nabla \dot\psi(s) \rangle|\, ds \to 0 \,.
\end{equation}

By (\ref{F7.5}) and Remark~\ref{nonat} we can apply Lemma~\ref{perturbation}. Therefore from (\ref{2170}) and (\ref{2013.5}) for a.e.\ $s\in[0,T]$ we obtain
$$
\langle \partial\F(s)(\uk(s)+\phik(s)), \dot\psi(s) \rangle -
\langle \partial\F(s)(\uk(s)), \dot\psi(s) \rangle \to 0\,,
$$
and by (\ref{F10}) and (\ref{2170}) we have
\begin{eqnarray*}
&|\langle \partial\F(s)(\uk(s)+\phik(s)), \dot\psi(s) \rangle -
\langle \partial\F(s)(\uk(s)), \dot\psi(s) \rangle| \le
\\
&\le 2 (\alpha_2^\F (C+1)^{q-1} + \beta_2^\F)\,\| \dot\psi(s)\|_{q}\,.
\end{eqnarray*}
As $s\mapsto \dot\psi(s)$ belongs to $L^1([0,T];L^{q}(\Om;\Rm))$,
we deduce that
\begin{equation}\label{2212}
\int_0^T |\langle \partial\F(s)(\uk(s)+\phik(s)), \dot\psi(s) \rangle -
\langle \partial\F(s)(\uk(s)), \dot\psi(s) \rangle|\, ds \to 0 \,.
\end{equation}

Finally, by (\ref{G7.5}) and Remark~\ref{nonat} we can apply Lemma~\ref{perturbation} again. Therefore from (\ref{2170}) and 
 (\ref{2023.5}) for a.e.\ $s\in[0,T]$ we obtain
$$
\langle \partial\G(s)(\uk(s)+\zetak(s)), \dot\psi(s) \rangle -
\langle \partial\G(s)(\uk(s)), \dot\psi(s) \rangle \to 0\,,
$$
and by (\ref{G10}) and (\ref{2170}) we have
\begin{eqnarray*}
&|\langle \partial\G(s)(\uk(s)+\zetak(s)), \dot\psi(s) \rangle -
\langle \partial\G(s)(\uk(s)), \dot\psi(s) \rangle| \le
\\
&\le 2 (\alpha_2^\G (C+1)^{r-1} + \beta_2^\G)\,
\| \dot\psi(s)\|_{r,\partial_S\Om}\,.
\end{eqnarray*}
As $s\mapsto \dot\psi(s)$ belongs to 
$L^1([0,T];L^{r}(\partial_S\Om;\Rm))$, we deduce that
\begin{equation}\label{2222}
\int_0^T |\langle \partial\G(s)(\uk(s)+\zetak(s)), \dot\psi(s) \rangle -
\langle \partial\G(s)(\uk(s)), \dot\psi(s) \rangle|\, ds \to 0 \,.
\end{equation}

Let
\begin{eqnarray}
&\displaystyle R_k:=
\int_0^T |\langle \partial\W(\nabla\uk(s)+\Psik(s)), \nabla 
\dot\psi(s) \rangle -
\langle \partial\W(\nabla\uk(s)), \nabla \dot\psi(s) \rangle|\, ds +
\nonumber\\
&\displaystyle
+\int_0^T |\langle \partial\F(s)(\uk(s)+\phik(s)), \dot\psi(s) \rangle -
\langle \partial\F(s)(\uk(s)), \dot\psi(s) \rangle|\, ds +
\nonumber\\
&\displaystyle
+\int_0^T |\langle \partial\F(s)(\uk(s)), \dot\psi(s) \rangle -
\langle \partial\F_k(s)(\uk(s)), \dot\psi(s) \rangle|\, ds +
\label{2230}\\
&\displaystyle
+\int_0^T |\langle \partial\G(s)(\uk(s)+\zetak(s)), \dot\psi(s) \rangle -
\langle \partial\G(s)(\uk(s)), \dot\psi(s) \rangle|\, ds +
\nonumber\\
&\displaystyle
+\int_0^T |\langle \partial\G(s)(\uk(s)), \dot\psi(s) \rangle -
\langle \partial\G_k(s)(\uk(s)), \dot\psi(s) \rangle|\, ds +
\nonumber\\
&\displaystyle
{}+\gamma^\F_k
\int_{0}^{T} \|\dot\psi(s)\|_{\dot q} \, ds+
\gamma^\G_k
\int_{0}^{T} \|\dot\psi(s)\|_{r,\partial_S\Om} \, ds\,.
\nonumber
\end{eqnarray}
Then $R_k\to 0$ by (\ref{2014.6}), (\ref{2014.8}), (\ref{2024.8}), (\ref{2024.6}), (\ref{2202}), 
(\ref{2212}), and (\ref{2222}). 
Inequality (\ref{1999}) follows from (\ref{2040}) and~(\ref{2230}).
\end{proof}

\begin{remark}\label{6005}
By (\ref{W10}), (\ref{F10}), (\ref{F11}), (\ref{G10}), (\ref{G11}), and (\ref{2170}), the 
integral term in (\ref{1999}) is bounded uniformly with respect to 
$k$ and $i$, therefore the same property holds for 
$\E\ki(\uki,\gki)$. By (\ref{K9}), (\ref{E}), and (\ref{coerc1}), 
this implies that there exists a constant $M>0$ such 
that
\begin{equation}\label{boundGamma}
\hn(\Ga_k(t))\le M
\end{equation}
for every $t\in[0,T]$ and every $k$, where $\Ga_k(t)$ is defined in~(\ref{ukt}).
\end{remark}

\begin{remark}\label{6075}
We notice that the integral which appears in (\ref{1999}) can be written as a sum which involves $\F(t)$, $\G(t)$, $\psi(t)$ only at the discrete times $t=t\kj$. Indeed we have
\begin{eqnarray}
&\displaystyle \int_0^{\tki} \theta_k(s)\,ds=
\sum_{j=1}^i\, \langle \partial\W(\nabla u\kjm),
\nabla \psi\kj- \nabla \psi\kjm\rangle-{}\nonumber\\
&\displaystyle {}-\sum_{j=1}^i\, \langle \partial\F\kjm( u\kjm),
\psi\kj- \psi\kjm\rangle-
\sum_{j=1}^i\,[\F\kj(u\kjm)-\F\kjm(u\kjm)]-
{}\nonumber\\
&\displaystyle {}-\sum_{j=1}^i\, \langle \partial\G\kjm( u\kjm),
\psi\kj- \psi\kjm\rangle-
\sum_{j=1}^i\,[\G\kj(u\kjm)-\G\kjm(u\kjm)]
{}\nonumber
\end{eqnarray}
for every $k$ and for every $i=1,\ldots,k$.
\end{remark}
\end{section}

\begin{section}{PROOF OF THE MAIN RESULT}\label{proof}

We are now in a position to prove the main result of the paper.

\begin{proof}[Proof of Theorem \ref{main}]
Let us fix a sequence of subdivisions $(\tki)_{0\le i\le k}$ of the interval $[0,T]$ satisfying (\ref{subdiv5}) and (\ref{fine5}), and
let $(u_0,\Ga_0)$ be an initial configuration satisfying the 
minimality property (\ref{mininit}).
For every $k$ let $(\uki,\gki)$, $i=1,\ldots,k$ be defined inductively as solutions of the discrete problems (\ref{2110}), with $(u_k^0,\Ga_k^0):=(u_0,\Ga_0)$, and let 
$\tau_k(t)$, $\uk(t)$, $\Ga_k(t)$, $\F_k(t)$, $\G_k(t)$, $\E_k(t)$ be defined by (\ref{ukt}).
By (\ref{boundGamma}) and Theorem~\ref{helly} there exist a 
subsequence, still denoted $\Ga_k$, and an increasing set function 
$t\mapsto\Ga^*(t)$, such that
\begin{equation}\label{convGamma}
\Ga_k^N(t)\ \,\sigmap\hbox{-converges to } \Ga^*(t)
\end{equation}
for every $t\in[0,T]$, where $\Ga_k^N(t):=\Ga_k(t)\cup \partial_N\Om$, according to~(\ref{gammanb}). Since $\Ga_k^N(t)\subsethn\overline\Om_B\cup\partial_N\Om$ and $\partial_N\Om$ is closed, we deduce, thanks to Remark~\ref{bgk}-(b),  that $\Ga^*(t)\subsethn\overline\Om_B\cup\partial_N\Om$ for every $t\in[0,T]$. Let $\Ga\colon[0,T]\to \Rp(\overline\Om_B)$ be the increasing set function defined by
\begin{equation}\label{Gammat}
\Ga(t):=\Ga^*(t)\setmeno\partial_N\Om
\end{equation}
By (\ref{eqganga}) and (\ref{semiK}) for every $t\in[0,T]$ we have
\begin{equation}\label{semiKK}
\K(\Ga(t))= \K(\Ga^*(t))\le \liminf_{k\to\infty} \K(\Ga_k^N(t))=
\liminf_{k\to\infty} \K(\Ga_k(t))\,.
\end{equation}

For a.e.\ $t\in[0,T]$ we set
\begin{equation}\label{theta}
\theta_\infty(t):=\limsup_{k\to\infty}\,\theta_k(t)\,,
\end{equation}
where $\theta_k$ is defined by (\ref{thetak}).
By (\ref{W10}), (\ref{F10}), (\ref{F11}), (\ref{G10}), (\ref{G11}), and (\ref{2170}) we have
\begin{equation}\label{dominated}
\begin{array}{c}
\vphantom{\displaystyle\min_A}
|\theta_k(t)| \le (\alpha_2^\W C^{p-1} + \beta_2^\W)\|\nabla \dot\psi(t)\|_p
+ (\alpha_2^\F C^{q-1} + \beta_2^\F)\,\| \dot\psi(t)\|_{q}+{}
\\ {}\qquad\qquad+ \alpha_3^\F(t) C^{\dot q} + \beta_3^\F(t)+ (\alpha_2^\G C^{r-1} + \beta_2^\G)\,
\| \dot\psi(t)\|_{r,\partial_S\Om}+\alpha_3^\G(t) C^{r} + \beta_3^\G(t)\,.
\end{array}
\end{equation}
As the right hand side of this inequality belongs to $L^1([0,T])$,
the function $\theta_\infty$ belongs to $L^1([0,T])$ and, by the Fatou lemma,
\begin{equation}\label{Fatou}
\limsup_{k\to\infty}\int_0^{\tau_k(t)} \theta_k(s)\,ds\le
\int_0^t \theta_\infty(s) \,ds\,.
\end{equation}

For every $t\in[0,T]$ can extract a subsequence $\theta_{k_j}$ of 
$\theta_k$, depending on $t$, such that
\begin{equation}\label{thetalimit}
\theta_\infty(t):=\lim_{j\to\infty}\,\theta_{k_j}(t)\,.
\end{equation}

By (\ref{6000}) and (\ref{convGamma}) we can 
apply Theorem~\ref{convmin} to $\tau_{k_j}(t)$, $\Ga_{k_j}(t)$, and 
$u_{k_j}(t)$. Therefore there exist a further subsequence, still 
denoted $u_{k_j}$, and a function $u(t)$ in $AD(\psi(t),\Ga(t))$ such 
that
\begin{eqnarray}
& u_{k_j}(t)\wto u(t) \, \hbox{ weakly in }\, GSBV^p(\Om;\Rm)\,,
\label{3123}\\
& u_{k_j}(t)\wto u(t) \, \hbox{ weakly in }\, L^{q}(\Om;\Rm)\,,
\label{3123.5}\\
& u_{k_j}(t)\to u(t) \, \hbox{ strongly in }\, L^{\dot q}(\Om;\Rm)\,,
\label{3124}\\
& u_{k_j}(t)\to u(t) \, \hbox{ strongly in }\, L^{r}(\partial_S\Om;\Rm)\,,
\label{3125}\\
& \nabla u_{k_j}(t)\wto \nabla u(t) \, \hbox{ weakly in }\, L^p(\Om;\Mmn)\,.
\label{3126}
\end{eqnarray}
Moreover
\begin{equation}\label{3min*}
\E(t)(u(t) ,\Ga(t))\le \E(t)(v,\Ga)
\end{equation}
for every $\Ga\in\Rp(\overline\Om_B)$, with $\Ga(t) \subsethn\Ga$, and every 
$v \in AD(\psi(t), \Ga)$.
This proves that $t\mapsto(u(t),\Ga(t))$ satisfies properties 
(a) and (b) of the definition of quasistatic evolution (Subsection~\ref{quasistatic}).

By (\ref{3123.5})--(\ref{3126}) there exists a constant $C>0$ such that
\begin{equation}\label{3170}
\|\nabla u(t)\|_p\le C\,, \quad \| u(t)\|_{\dot q}\le C\,,\quad
\| u(t)\|_{q}\le C\,,\quad \| u(t)\|_{r,\partial_S\Om }\le C
\end{equation}
for every $t\in[0,T]$.

Finally Theorem~\ref{convmin} implies that
\begin{eqnarray}
& \W(\nabla u_{k_j}(t))\to \W(\nabla u(t) )\,,
\label{3127}\\
& \F_{k_j}(t)(u_{k_j}(t))\to \F(t)(u(t))\,,\label{3128}\\
& \G_{k_j}(t)(u_{k_j}(t))\to \G(t)(u(t))\,.
\label{3129}
\end{eqnarray}
Using Lemma~\ref{stressconv}, from
(\ref{3123}) and (\ref{3127}), we obtain
\begin{equation}\label{3200}
\langle \partial\W(\nabla u_{k_j}(t)), \nabla \dot\psi(t) \rangle \to
\langle \partial\W(\nabla u(t)), \nabla \dot\psi(t) \rangle\,.
\end{equation}
By (\ref{3124}) the sequence $\partial_zF(t,x,u_{k_j}(t))$ converges to $\partial_zF(t,x,u(t))$ in measure on~$\Om$. By (\ref{F10}) and (\ref{3123.5}) it is bounded in $L^{q'}(\Om;\Rm)$. Therefore  $\partial_zF(t,x,u_{k_j}(t))$ converges to $\partial_zF(t,x,u(t))$ weakly in $L^{q'}(\Om;\Rm)$, and consequently by (\ref{F7.5})
\begin{equation}\label{3200.5}
\langle \partial\F(t)(u_{k_j}(t)), \dot\psi(t) \rangle
\to \langle \partial\F(t)(u(t)), \dot\psi(t) \rangle \,,
\end{equation}
so that by (\ref{2014.6}) and (\ref{2014.8}) we obtain
\begin{equation}
\langle \partial\F_{k_j}(t)(u_{k_j}(t)), \dot\psi(t) \rangle
\to \langle \partial\F(t)(u(t)), \dot\psi(t) \rangle \,.\label{3201}
\end{equation}
Observing that $\dot \F(t)$ is continuous on $L^{\dot q}(\Om;\Rm)$, while
$\partial\G(t)$ and $\dot \G(t)$ are continuous on 
$L^{r}(\partial_S\Om;\Rm)$,
by (\ref{3124}) and (\ref{3125}) we have
\begin{eqnarray}
& \dot \F(t)(u_{k_j}(t)) \to \dot \F(t)(u(t)) \,, \label{3202}\\
& \langle \partial\G(t)(u_{k_j}(t)), \dot\psi(t) \rangle \to
\langle \partial\G(t)(u(t)), \dot\psi(t) \rangle \,, \label{3202.5}\\
& \dot \G(t)(u_{k_j}(t)) \to \dot \G(t)(u(t)) \,, \label{3204}
\end{eqnarray}
so that by (\ref{2024.8}) and (\ref{2024.6}) we obtain
\begin{equation}
\langle \partial\G_{k_j}(t)(u_{k_j}(t)), \dot\psi(t) \rangle
\to \langle \partial\G(t)(u(t)), \dot\psi(t) \rangle \,.\label{3203}
\end{equation}

For a.e.\ $t\in[0,T]$ we set
\begin{equation} \label{3205}
\begin{array}{c} \vphantom{\displaystyle \min_A}
\theta(t):= \langle \partial\W(\nabla u(t)), \nabla \dot\psi(t) \rangle
- \langle \partial\F(t)(u(t)), \dot\psi(t) \rangle - {}\\
 \vphantom{\displaystyle \min_A}
{} - \dot \F(t)(u(t))- \langle \partial\G(t)(u(t)), 
\dot\psi(t) \rangle -
\dot \G(t)(u(t)) =\\
= \langle g(t), \dot\psi(t) \rangle-\dot \F(t)(u(t)) - \dot \G(t)(u(t))\,,
\end{array}
\end{equation}
where $g(t)$ is defined by (\ref{g(t)}).
From (\ref{thetak}), (\ref{thetalimit}), (\ref{3200}), 
(\ref{3201}), (\ref{3202}), (\ref{3204}), and (\ref{3203}) we obtain
\begin{equation} \label{3205.5}
\theta_\infty(t)=\theta(t)
\end{equation}
for a.e.\ $t\in[0,T]$.

By (\ref{semiKK}) and 
(\ref{3127})--(\ref{3129}) we have
\begin{equation}\label{3206}
\E(t)(u(t),\Ga(t))\le \liminf_{j\to\infty}\, \E_{k_j}(t)(u_{k_j}(t), 
\Ga_{k_j}(t))\le \limsup_{k\to\infty} \E_{k}(t)(u_{k}(t), 
\Ga_{k}(t))\,.
\end{equation}
{}From (\ref{1999'}), (\ref{Fatou}), and (\ref{3205.5}) we obtain
\begin{equation}\label{3207}
\limsup_{k\to\infty} \E_{k}(t)(u_{k}(t), \Ga_{k}(t)) \le 
\E(0)(u_0,\Ga_0) + \int_0^{t} \theta(s)\,ds\,.
\end{equation}
By (\ref{3206}) and (\ref{3207}) we have
\begin{equation}\label{3208}
\E(t)(u(t),\Ga(t))\le \E(0)(u_0,\Ga_0) + \int_0^{t} \theta(s)\,ds\,.
\end{equation}

To conclude the proof of the theorem it is enough to show that
\begin{equation}\label{7000}
\E(t)(u(t),\Ga(t))\ge \E(0)(u_0,\Ga_0) + \int_0^{t} \theta(s)\,ds\,.
\end{equation}
for every $t\in[0,T]$. Indeed (\ref{3208}) and (\ref{7000}) imply that
$t\mapsto E(t):=\E(t)(u(t),\Ga(t))$ is absolutely continuous in 
$[0,T]$ and that $\dot E(t)=\theta(t)$ for a.e.\ $t\in[0,T]$. By (\ref{3205}) this 
yields condition~(c) 
of the definition of quasistatic evolution (Subsection~\ref{quasistatic}).
\end{proof}

In order to prove (\ref{7000}) we need the following lemma.

\begin{lemma}\label{estimate below}
Assume that $t\mapsto(u(t),\Ga(t))$ satisfies conditions (a) and (b) in the definition of quasistatic evolution (Subsection~\ref{quasistatic}).
Let $t\in[0,T]$ and let $\ski$ be a subdivision of $[0,t]$ 
satisfying (\ref{subdivs}) and (\ref{fines}). Then
there exists a sequence $S_k(t)\to 0$ such that
\begin{eqnarray}
&\displaystyle
\E(t)(u(t),\Ga(t)) \ge \E(0)(u_0,\Ga_0)
+ \sum_{i=1}^{i_k} \int_{\skim}^{\ski} \langle \partial\W(\nabla u(\ski)),
\nabla \dot\psi(s) \rangle \,ds - \nonumber\\
&\displaystyle - \sum_{i=1}^{i_k} \int_{\skim}^{\ski} \langle 
\partial\F(\ski)(u(\ski)),
\dot\psi(s) \rangle \,ds - \sum_{i=1}^{i_k} \int_{\skim}^{\ski} 
\dot \F(s)(u(\ski)) \,ds- \label{2999}\\
&\displaystyle -\sum_{i=1}^{i_k} \int_{\skim}^{\ski} \langle 
\partial\G(\ski)(u(\ski)),
\dot\psi(s) \rangle \,ds - \sum_{i=1}^{i_k} \int_{\skim}^{\ski} 
\dot \G(s)(u(\ski)) \,ds - S_k(t)\,.
\nonumber
\end{eqnarray}
\end{lemma}

\begin{proof} The minimality property in condition (a) gives $\E(t)(u(t),\Ga(t))\le \E(t)(\psi(t),\Ga(t))$, hence $\Ec(t)(u(t))\le \Ec(t)(\psi(t))$.
By (\ref{coerc1}) and (\ref{bound1}), this implies that 
 there exists a constant $C>0$ such that
\begin{equation}\label {3170'}
\|\nabla u(t)\|_p\le C\,, \quad \|u(t)\|_{\dot q}\le C\,,\quad
\|u(t)\|_{q}\le C
\end{equation}
for every $t\in[0,T]$.
By Lemma~\ref{rem3} the functions $u(t)$ belong to 
$W^{1,p}(\Om_S;\Rm)\cap L^{q}(\Om_S;\Rm)$. Therefore 
(\ref{tracepq}) implies that, if we change the constant, we may assume also that
\begin{equation}\label{3170"}
\|u(t)\|_{r,\partial_S\Om }\le C
\end{equation}
for every $t\in[0,T]$.

For every $k$ let
$\sk\colon(0,t]\to(0,t]$ be the piecewise constant function defined by
$$
\sk(s)=\ski\qquad\hbox{for }\skim<s\le\ski\,,
$$
and let
$\vk\colon (0,t]\to GSBV^p_{q}(\Om;\Rm)$ be the piecewise constant 
function defined by
\begin{equation}\label{502.5}
\vk(s):= u(\ski)=u(\sk(s)) \qquad \hbox{for }\, \skim< s \le \ski\,.
\end{equation}

For every $i=1,\ldots,i_k$ we have $u(\ski)-\psi(\ski)+\psi(\skim)\in 
AD(\psi(\skim), \Ga(\ski))$ and $\Ga(\skim)\subsethn\Ga(\ski)$. By 
the minimality property in condition (a) we have
\begin{equation}\label{501}
\E(\skim)(u(\skim),\Ga(\skim))\le
\E(\skim)(u(\ski)-\psi(\ski)+\psi(\skim),\Ga(\ski))\,.
\end{equation}

Arguing as in the proof of (\ref{2040}) we obtain that there exist two sequences $\delta^\F_k$ and $\delta^\G_k$ of real numbers, with 
\begin{equation}\label{514.8}
\delta^\F_k\to 0\,, \qquad \delta^\G_k\to 0\,,
\end{equation}
and three sequences $\Xk\in L^\infty([0,T];L^p(\Om;\Mmn))$, $\chik\in 
L^\infty([0,T];L^{q}(\Om;\Rm))$, and $\etak\in 
L^\infty([0,T];L^{r}(\partial_S\Om;\Rm))$, with
\begin{equation}\label{503.5}
\|\Xk(s)\|_p+\|\chik(s)\|_{q}+\|\etak(s)\|_{r,\partial_S\Om}\to 0\,\hbox{ uniformly with respect to }\, s\in[0,T]\,, 
\end{equation}
such that
\begin{eqnarray}
&\displaystyle \vphantom{\int_{\tkim}^{\tki}}
\E(t)(u(t),\Ga(t)) - \E(0)(u_0,\Ga_0) \ge \nonumber\\
&\displaystyle
\ge \int_{0}^t \langle \partial\W(\nabla\vk(s)+\Xk(s)),
\nabla \dot\psi(s) \rangle \,ds - \nonumber\\
&\displaystyle - \int_{0}^{t} \langle \partial\F(\sk(s))(\vk(s)+\chik(s)),
\dot\psi(s) \rangle \,ds - \int_{0}^{t} \dot \F(s)(\vk(s)) \,ds- 
\label{540}\\
&\displaystyle - \int_{0}^{t} \langle \partial\G(\sk(s))(\vk(s)+ \etak(s)),
\dot\psi(s) \rangle \,ds - \int_{0}^{t} \dot \G(s)(\vk(s)) \,ds -
\nonumber\\
&\displaystyle -\delta^\F_k
\int_{0}^{t} \|\dot\psi(s)\|_{\dot q} \, ds-
\delta^\G_k
\int_{0}^{t} \|\dot\psi(s)\|_{r,\partial_S\Om} \, ds\,.
\nonumber
\end{eqnarray}
Let
\begin{eqnarray}
&\displaystyle S_k(t):=
\int_0^T |\langle \partial\W(\nabla\vk(s)+\Xk(s)), \nabla \dot\psi(s) \rangle -
\langle \partial\W(\nabla\vk(s)), \nabla \dot\psi(s) \rangle|\, ds +
\nonumber\\
&\displaystyle
+\int_0^T |\langle \partial\F(\sk(s))(\vk(s)+\chik(s)), \dot\psi(s) \rangle -
\langle \partial\F(\sk(s))(\vk(s)), \dot\psi(s) \rangle|\, ds +
\label{630}\\
&\displaystyle
+\int_0^T |\langle \partial\G(\sk(s))(\vk(s)+\etak(s)), \dot\psi(s) \rangle -
\langle \partial\G(\sk(s))(\vk(s)), \dot\psi(s) \rangle|\, ds +
\nonumber\\
&\displaystyle
{}+\delta^\F_k
\int_{0}^{T} \|\dot\psi(s)\|_{\dot q} \, ds+
\delta^\G_k
\int_{0}^{T} \|\dot\psi(s)\|_{r,\partial_S\Om} \, ds\,.
\nonumber
\end{eqnarray}
Using Lemma~\ref{perturbation} as in the last part of the proof of Lemma~\ref{energy estimate}, from (\ref{514.8}) and (\ref{503.5})
 we obtain that $S_k(t)\to 0$. Inequality (\ref{2999}) follows from 
(\ref{540}) and~(\ref{630}).
\end{proof}

\begin{proof}[Proof of Theorem~\ref{main} continued]
Let us fix $t\in(0,T]$.
By Lemmas~\ref{riemann1} and~\ref{riemann4} and 
Remarks~\ref{riemann2} and~\ref{riemann3} there exists a subdivision 
$\ski$ of $[0,t]$, satisfying (\ref{subdivs}), (\ref{fines}), 
(\ref{riemsumk}), and (\ref{riemsumkg}), such that
\begin{eqnarray}
&\displaystyle
\lim_{k\to\infty} \,\sum_{i=1}^{i_k} \Big | (\ski-\skim) \theta(\ski) -
\int_{\skim}^{\ski} \theta (s)\,ds \,\Big| = 0\,,
\label{riemsumtheta}\\
&\displaystyle
\lim_{k\to\infty} \,\sum_{i=1}^{i_k} \Big \| (\ski-\skim)\nabla 
\dot\psi(\ski) -
\int_{\skim}^{\ski} \nabla \dot \psi(s)\,ds \,\Big\|_p = 0\,,
\label{riemsumnablapsi}\\
&\displaystyle
\lim_{k\to\infty} \,\sum_{i=1}^{i_k} \Big \| (\ski-\skim) \dot \psi(\ski) -
\int_{\skim}^{\ski} \dot \psi(s)\,ds \,\Big\|_{q} = 0\,,
\label{riemsumpsi}\\
&\displaystyle
\lim_{k\to\infty} \,\sum_{i=1}^{i_k} \Big \| (\ski-\skim) \dot \psi(\ski) -
\int_{\skim}^{\ski} \dot \psi(s)\,ds \,\Big\|_{r,\partial_S\Om} = 0\,,
\label{riemsumpsiS}
\end{eqnarray}

Let
$$
{\mathcal H}:=\{u\in GSBV^p_{q}(\Om;\Rm):S(u)\subsethn\overline\Om_B,\ 
\hn(S(u))\le M,\ \|\nabla u\|_p\le C,\ \|u\|_{q}\le C\}\,,
$$
where $M$ and $C$ are the constants which appear in (\ref{boundGamma}) and (\ref{3170'}). 
Arguing as in Theorem~\ref{convmin} it is easy to prove that ${\mathcal H}$ is 
compact in $L^{\dot q}(\Om;\Rm)$ and in $L^{r}(\partial_S\Om;\Rm)$.

Let
\begin{eqnarray}
&\displaystyle
\om_k(t):=
\sum_{i=1}^{i_k} \Big\| \partial\W(\nabla u(\ski))\Big\|_{p'} \Big\| 
(\ski-\skim)\nabla \dot\psi(\ski) - \int_{\skim}^{\ski}
\nabla \dot\psi(s) \,ds \Big\|_p + {}\nonumber\\
&\displaystyle {}+\sum_{i=1}^{i_k}
\Big\| \partial\F(\ski)(u(\ski))\Big\|_{q'}
\Big\| (\ski-\skim) \dot\psi(\ski) - \int_{\skim}^{\ski}
\dot\psi(s) \,ds\Big\|_{q} +{}\nonumber\\
&\displaystyle {}+
\sum_{i=1}^{i_k}\, \sup_{u\in{\mathcal H} } \,
\Big | (\ski-\skim) \dot \F(\ski)(u) -
\int_{\skim}^{\ski} \dot \F(s)(u)\,ds \,\Big| +
{}\label{4000}\\
&\displaystyle{} +\sum_{i=1}^{i_k}
\Big\| \partial\G(\ski)(u(\ski))\Big\|_{r',\partial_S\Om}
\Big\| (\ski-\skim) \dot\psi(\ski) - \int_{\skim}^{\ski}
\dot\psi(s) \,ds\Big\|_{r,\partial_S\Om}+{}\nonumber\\
&\displaystyle {}+\sum_{i=1}^{i_k}\, \sup_{u\in{\mathcal H} } \,
\Big | (\ski-\skim) \dot \G(\ski)(u) -
\int_{\skim}^{\ski} \dot \G(s)(u)\,ds \,\Big|+ S_k(t)\,,
\nonumber
\end{eqnarray}
where $S_k(t)$ is given by Lemma~\ref{estimate below}. 
By (\ref{W10}), (\ref{F10}), (\ref{G10}), (\ref{3170'}), and (\ref{3170"}) the terms
$\|\partial\W(\nabla u(\ski))\|_{p'}$, 
$\|\partial\F(\ski)(u(\ski))\|_{q'}$, and
$\| \partial\G(\ski)(u(\ski))\|_{r',\partial_S\Om}$ are bounded 
uniformly with respect to $k$ and $i$. Therefore, using 
(\ref{riemsumk}), (\ref{riemsumkg}), 
(\ref{riemsumnablapsi})--(\ref{riemsumpsiS}), and 
Lemma~\ref{estimate below} we obtain
\begin{equation}\label{4001}
\lim_{k\to\infty}\,\om_k(t)=0\,.
\end{equation}

By (\ref{boundGamma}) and (\ref{3170'}) we have
$u(\ski)\in{\mathcal H}$ for every $k$ and $i$.  Taking (\ref{3205}) and (\ref{4000}) into account, the estimate from below in
Lemma~\ref{estimate below} gives
\begin{equation}\label{800}
\E(t)(u(t),\Ga(t)) \ge \E(0)(u_0,\Ga_0) + \sum_{i=1}^{i_k}
(\ski-\skim) \theta(\ski)- \om_k(t)\,.
\end{equation}
Passing to the limit as $k\to\infty$,
by (\ref{riemsumtheta}) we obtain (\ref{7000}),
which, together with (\ref{3205}) and (\ref{3208}), yields
condition (c) in the definition of quasistatic evolution.
\end{proof}

\end{section}

\begin{section}{CONVERGENCE OF THE DISCRETE-TIME PROBLEMS}\label{convergence}
 
In this section we show that for every $t\in [0,T]$ the elastic energies and the crack energies of the solutions to the discrete-time problems converge to the corresponding energies for the continuous-time problem. Note that this result can be proved even if the minimum energy deformations corresponding to a crack $\Ga(t)$ are not unique, but that it only holds
for the discretizations that produce a given crack $\Ga(t)$. 

Let $\tki$ be as in Section~\ref{discrete}, let $(u_0,\Ga_0)$ be an 
initial configuration that satisfies the minimality property 
(\ref{mininit}), and let $(\uki,\gki)$ be the solutions of
the minimum problems (\ref{2110}), with 
$(u_k^0,\Ga_k^0):=(u_0,\Ga_0)$. Let $\tau_k$, $\uk$,  $\Ga_k$, $\F_k$, $\G_k$, and $\E_k$ be the 
piecewise constant functions introduced in (\ref{ukt}), and let 
 $\E_{k}^{el}\colon[0,T]\to\R$ be the piecewise constant function defined 
by
\begin{equation}\label{eckt}
\E_{k}^{el}(t):=\Ec(\tki)=\Ec(\tau_k(t))\,,
\end{equation}
where $i$ is the largest integer such that $\tki\le t$.

\begin{theorem}Let $t\mapsto(v(t),\Ga(t))$ be a quasistatic evolution, let $\theta_k$ be defined by (\ref{thetak}), and let 
\begin{equation} \label{3205+}
\begin{array}{c} \vphantom{\displaystyle \min_A}
\theta(t):= \langle \partial\W(\nabla v(t)), \nabla \dot\psi(t) \rangle
- \langle \partial\F(t)(v(t)), \dot\psi(t) \rangle - {}\\
 \vphantom{\displaystyle \min_A}
{} - \dot \F(t)(v(t))- \langle \partial\G(t)(v(t)), 
\dot\psi(t) \rangle -
\dot \G(t)(v(t))\,.
\end{array}
\end{equation}
Assume that $\Ga_k(t)$ and $\Ga(t)$ satisfy (\ref{convGamma}) and  (\ref{Gammat})
for every $t\in[0,T]$. Then
\begin{eqnarray}
&\displaystyle\vphantom{L^1}
\Ec(t)(v(t))= \lim_{k\to\infty} \, \E_{k}^{el}(t)(\uk(t))
\,,
\label{lime0}\\
&\displaystyle \vphantom{L^1}
\K(\Ga(t)) =
\lim_{k\to\infty} \, \K(\Ga_k(t))
\label{limkappa}
\end{eqnarray}
for every $t\in[0,T]$. Moreover
\begin{equation}\label{limth}
 \theta_k\to \theta\quad\hbox{in }L^1([0,T])\,,
\end{equation}
so that there exists a subsequence of $\theta_k$ which converges to 
$\theta$ a.e.\ in $[0,T]$.
\end{theorem}
 
\begin{proof} 
For a.e.\ $t\in[0,T]$ let $\theta_\infty(t)$ be defined by (\ref{theta}). In the proof of Theorem~\ref{main} for every $t\in[0,T]$ we constructed a function $u(t)\in AD(\psi(t),\Ga(t))$ such that $t\mapsto(u(t),\Ga(t))$ is a quasistatic evolution and 
\begin{equation}\label{eu}
\E(t)(u(t),\Ga(t))=\E(0)(u_0,\Ga_0)+\int_0^t\theta_\infty(s)\,ds
\end{equation}
for every $t\in[0,T]$ (see (\ref{3205.5}), (\ref{3208}), and (\ref{7000})).
Since $(u(t),\Ga(t))$ and $(v(t),\Ga(t))$ satisfy the minimality condition (a) in the definition of quasistatic evolution (see Subsection~\ref{quasistatic}), we have 
\begin{equation}\label{euev}
\E(t)(u(t),\Ga(t))=\E(t)(v(t),\Ga(t))
\end{equation}
for every $t\in[0,T]$. By condition (c) for $t\mapsto (v(t),\Ga(t))$ we have
\begin{equation}\label{ev}
\E(t)(v(t),\Ga(t))=\E(0)(u_0,\Ga_0)+\int_0^t\theta(s)\,ds
\end{equation}
{}From (\ref{eu})--(\ref{ev}) we deduce that 
\begin{equation}\label{thetainfty}
\theta(t)=\theta_\infty(t)
\end{equation}
for a.e.\ $t\in[0,T]$.

{}By Lemma~\ref{energy estimate} for every $t\in[0,T]$ we have
\begin{eqnarray}
&\displaystyle \liminf_{k\to\infty} \, \E_k(t)(\uk(t),\Ga_k(t))\le 
\E(0)(u_0,\Ga_0) +\liminf_{k\to\infty} \int_0^{\tau_k(t)} 
\theta_k(s)\,ds\,,
\label{liminf0}\\
&\displaystyle
\limsup_{k\to\infty} \, \E_k(t)(\uk(t),\Ga_k(t))\le \E(0)(u_0,\Ga_0) 
+\limsup_{k\to\infty} \int_0^{\tau_k(t)} \theta_k(s)\,ds\,.
\label{limsup}
\end{eqnarray}
Let us fix $t\in[0,T]$ and let $u_{k_j}(t)$ be a subsequence of 
$\uk(t)$ such that
\begin{equation}\label{1101}
\lim_{j\to\infty} \, \E_{k_j}^{el}(t)(u_{k_j}(t))=
\liminf_{k\to\infty} \, \E_{k}^{el}(t)(\uk(t))\,.
\end{equation}
Since $\Ga_{k_j}^N(t)$ $\sigmap$-converges to $\Ga^*(t)$, using (\ref{6000}) we can 
apply Theorem~\ref{convmin} to $\tau_{k_j}(t)$, $\Ga_{k_j}(t)$, and 
$u_{k_j}(t)$. Therefore there exist a further subsequence, still 
denoted $u_{k_j}(t)$, and a function $u^*(t)\in AD(\psi(t),\Ga(t))$ 
such that
\begin{eqnarray}
& \W(\nabla u_{k_j}(t))\to \W(\nabla u^*(t) )\,,
\label{1127}\\
&\F_{k_j}(t)(u_{k_j}(t))\to \F(t)(u^*(t))\,,
\label{1128}\\
& \G_{k_j}(t)(u_{k_j}(t))\to \G(t)(u^*(t))\,.
\label{1129}
\end{eqnarray}
Moreover 
$$
\E(t)(u^*(t) ,\Ga(t))\le \E(t )(v,\Ga)
$$
for every 
$\Ga\in \Rp(\overline\Om_B)$, with $\Ga(t) \subsethn\Ga$, and for every 
$v \in AD(\psi(t),\Ga)$.
Since $(v(t) ,\Ga(t))$ satisfies the same minimality property by 
condition (a) in Subsection~\ref{quasistatic}, we have
\begin{equation}\label{1130}
\Ec(t)(v(t))=\Ec(t)(u^*(t))\,.
\end{equation}
{}From (\ref{1127})--(\ref{1129}) we obtain
$$
\Ec(t)(u^*(t))= \lim_{j\to\infty} \, \E_{k_j}^{el}(t)(u_{k_j}(t))\,,
$$
which, together with (\ref{1101}) and (\ref{1130}) gives
\begin{equation}\label{liminfe0}
\Ec(t)(v(t))=
\liminf_{k\to\infty} \, \E_{k}^{el}(t)(\uk(t))\,.
\end{equation}

By (\ref{eqganga}) and (\ref{semiK}) we have
\begin{equation}\label{liminfkappa}
\K(\Ga(t))=\K(\Ga^*(t))\le \liminf_{k\to\infty} \, \K(\Ga_k^N(t))
=\liminf_{k\to\infty} \, \K(\Ga_k(t))\,,
\end{equation}
so that (\ref{liminfe0}) and (\ref{liminfkappa}) yield
\begin{equation}\label{liminf}
\E(t)(v(t),\Ga(t))\le
\liminf_{k\to\infty} \, \E_k(t)(\uk(t),\Ga_k(t))\,.
\end{equation}
From (\ref{Fatou}), (\ref{eu}), (\ref{euev}), (\ref{thetainfty}), (\ref{liminf0}), (\ref{limsup}), and (\ref{liminf}) we obtain
\begin{eqnarray}
&\displaystyle \E(t)(v(t),\Ga(t))= \lim_{k\to\infty} \, \E_k(t)(\uk(t),\Ga_k(t))\label{lim}\\
&\displaystyle \int_0^t \theta(s)\,ds=\lim_{k\to\infty} 
\int_0^{\tau_k(t)} \theta_k(s)\,ds
\label{limtheta}
\end{eqnarray}
for every $t\in[0,T]$.
Equalities (\ref{lime0}) and~(\ref{limkappa}) follow easily from (\ref{liminfe0}), (\ref{liminfkappa}), and (\ref{lim}). 

By (\ref{limtheta}) we have
\begin{equation}\label{1140}
\int_0^T \theta(t)\,dt=\lim_{k\to\infty} \int_0^T \theta_k(t)\,dt \,.
\end{equation}
By (\ref{theta}) and (\ref{thetainfty}) $\theta_k\lor\theta$ converges to $\theta$ pointwise 
on $[0,T]$, so that by (\ref{dominated})
\begin{equation}\label{1141}
\theta_k\lor\theta\to \theta\quad\hbox{in }L^1([0,T])\,.
\end{equation}
Since $\theta_k+\theta= (\theta_k\lor\theta) + (\theta_k\land \theta)$, from
(\ref{1140}) and (\ref{1141}) we obtain
$$
\int_0^T \theta(t)\,dt=\lim_{k\to\infty} \int_0^T (\theta_k\land 
\theta)(t)\,dt\,.
$$
As $\theta_k\land \theta\le \theta$, this implies that $ 
\theta_k\land \theta$ converges to $\theta$ in $L^1([0,T])$, which, 
together with (\ref{1141}), gives~(\ref{limth}).
\end{proof}
\end{section}

\begin{section}{APPENDIX}\label{ipgen}

In Remarks~\ref{0109} and~\ref{01091} we introduced  elementary conditions on $F$ and $G$ that ensure  that $\F$, $\dot\F$, $\G$, and $\dot\G$ satisfy all properties required in 
Subsections \ref{body} and~\ref{surfacef}.
It is easy to see that (\ref{F2}) and (\ref{G2}) are much stronger than what we need for this purpose. Indeed, it is enough to assume that the partial derivatives $\partial_zF$, $\partial_t F$, $\partial_z\partial_t F$, $\partial_zG$, $\partial_t G$, and $\partial_z\partial_t G$ exist, are measurable with respect to $x$, and continuous with respect to $(t,z)$. In the rest of this section we show that the same results can be obtained under a less regular dependence on time, in the spirit of the hypotheses considered for the boundary deformations in Subsection~\ref{bounddef}.

\subsection{Weaker hypotheses on the body forces}\label{newbody}
We are interested in particular in the case of {\em dead loads\/}, in which the density of the body force $f\colon[0,T]{\times}\Om\to\Rm$ per unit volume in the reference configuration does not depend on the deformation. In this case the simplest choice for the potential is $F(t,x,z):=f(t,x)z$. But this linear dependence on $z$ can not be accepted, because it violates the first inequality in (\ref{F9}), where $\alpha_0^\F>0$. We may consider a slight variant, namely 
\begin{equation}\label{dlf}
F(t,x,z):=f(t,x)z+F_0(x,z)\,,
\end{equation}
where $F_0\colon\Om{\times}\Rm\to\R$ is a Carath\'eodory function. We assume that for every $x\in\Om$ the function $z\mapsto F_0(x,z)$ belongs to $C^1(\Rm)$ and that for every $(x,z)\in\Om{\times}\Rm$
\begin{eqnarray}
& a_0 |z|^{q} - b_0(x)\le -F_0(x,z) \le
a_1 |z|^{q} + b_1(x)\,,\label{alphaF01}
\\
& |\partial_z F_0(x,z)| \le a_2 |z|^{q-1} + b_2(x)\,,\label{alphaF02}
\end{eqnarray}
where $q>1$, $a_0,\,a_1,\,a_2$ are positive constants, 
$b_0,\,b_1\in L^1(\Om)$, $b_2\in L^{q'}(\Om)$, and $q':=q/(q-1)$.

When the body is subject to a deformation $u$, the body force acting 
at time $t$ has a density per unit volume in the reference 
configuration given by $f(t,x)+\partial_zF_0(x,u(x))$. The term 
$f(t,x)$ is a time dependent dead load (that we may think of as 
determined by an experimental device), while $\partial_zF_0(x,u(x))$ 
can be interpreted as a background time independent body force.
As a consequence of our hypotheses on $F_0$, this force will prevent  broken parts of the body from finding   infinity as only equilibrium configuration.

Let $\dot q$ be a constant in $(1,q)$ and let
$\dot q'=\dot q/(\dot q-1)$.
If $t\mapsto f(t,\cdot)$ is absolutely continuous from $[0,T]$ into 
$L^{\dot q'}(\Om;\Rm)$ and $t\mapsto\dot f(t,\cdot)$ is its time 
derivative, which belongs to $L^1([0,T];L^{\dot q'}(\Om;\Rm))$, we can consider the functional $\F(t)\colon L^q(\Om;\Rm)\to\R$ defined  for every $t\in[0,T]$ by (\ref{F7}), with $F$ given by (\ref{dlf}), and the functional 
$\dot\F(t)\colon L^{\dot q}(\Om;\Rm)\to\R$ defined for a.e.\ $t\in[0,T]$ by
$$
\dot\F(t)(u):=\int_\Om \dot f(t,x)\,u(x)\,dx\,.
$$
We can easily check that in this case $\F$ and $\dot\F$ satisfy all properties required in Subsection~\ref{body}.

Note that these hypotheses do not guarantee the existence of the partial derivative $\partial_t F(t,x,z)$ for a.e.\ $t\in[0,T]$ and for every $(x,z)\in\Om{\times}\Rm$, so that Remark~\ref{0109} can not be applied.

We now present  a more general set of hypotheses, which includes this case as well as those considered in Remark~\ref{0109}. More precisely, we assume that
$F\colon[0,T]{\times}\Om{\times}\Rm\to\R$ satisfies the following conditions
\begin{eqnarray}
&&\hbox{\hspace{-1.3cm}for every $z\in\Rm$ the function $(t,x)\mapsto F(t,x,z)$ is 
${\mathcal L}^1{\times} \Ln$-measurable on 
$[0,T]{\times}\Om$,}\label{FF1}\\
&&\hbox{\hspace{-1.3cm}for every $(t,x)\in[0,T]{\times}\Om$ the function $z\mapsto F(t,x,z)$
belongs to $C^1(\Rm)$.}\label{FF2}
\end{eqnarray}
Moreover, we assume that that there exist four
constants $q> 1$, $a_0^F>0$, $a_1^F>0$, $a_2^F>0$ and
three nonnegative functions
$b_0^F$, $b_1^F\in C^0([0,T];L^1(\Om))$,
$b_2^F\in C^0([0,T];L^{q'}(\Om))$,
such that
\begin{eqnarray}
&a_0^F |z|^{q{}} - b_0^F(t,x)\le - F(t,x,z) \le
a_1^F  |z|^{q{}} + b_1^F(t,x)\,,\label{alphaF*}
\\
& |\partial_z F(t,x,z)| \le a_2^F  |z|^{q{}-1} + b_2^F(t,x)\,,
\label{alphaF2*}
\end{eqnarray}
for every $(t,x,z)\in [0,T]{\times}\Om{\times}\Rm$.

To deal with the dependence of $F$ on $t$, we assume that there exists  a function \break
$\dot F\colon[0,T]{\times}\Om{\times}\Rm\to \R$ 
such that for every $t\in [0,T]$ and for every $z\in\Rm$
\begin{eqnarray}
&\displaystyle
F(t,x,z)=F(0,x,z)+\int_0^t\dot F(s,x,z)\,ds
\quad\hbox{for a.e.\ $x\in\Om$,}\label{F5}\\
&\displaystyle
\partial_z F(t,x,z)=\partial_z F(0,x,z)+\int_0^t \partial_z \dot 
F(s,x,z)\,ds \quad\hbox{for a.e.\ $x\in\Om$.}\label{F5.5}
\end{eqnarray}
The integrals in (\ref{F5}) and (\ref{F5.5}) are well defined since 
we assume also that
\begin{eqnarray}
&&\hbox{\hspace{-1.5cm}for every $z\in\Rm$ the function 
$(t,x)\mapsto \dot F(t,x,z)$ is 
${\mathcal L}^1{\times} \Ln$-measurable on 
$[0,T]{\times}\Om$,}\label{F3}\\
&&\hbox{\hspace{-1.5cm}for every $(t,x)\in[0,T]{\times}\Om$ the function $z\mapsto 
\dot F(t,x,z)$
belongs to $C^1(\Rm)$,}\label{F4}
\end{eqnarray}
and that there exist a constant $\dot q\in [1, q)$ and
four nonnegative functions
$a_3^F,a_4^F\!\in\! L^1([0,T])$,
$b_3^F\!\in\! L^1([0,T]; L^1(\Om))$, and $b_4^F\in L^1([0,T];
L^{\dot q'}(\Om))$ such that
\begin{eqnarray}
&|\dot F(t,x,z)|\le a_3^F(t)|z|^{\dot q}+b_3^F(t,x)\,,\label{F6*}\\
&|\partial_z \dot F(t,x,z)|\le a_4^F(t)|z|^{\dot q-1}+b_4^F(t,x)
\label{F6.5*}
\end{eqnarray}
for every $(t,x,z)\in [0,T]{\times}\Om{\times}\Rm$.

The functionals $\F(t)\colon L^q(\Om;\Rm)\to\R$ and 
$\dot \F(t)\colon 
L^{\dot q}(\Om;\Rm)\to\R$ are defined by (\ref{F7}) and by
\begin{equation}\label{star}
\dot\F(t)(u):=\int_\Om \dot F(t,x,u(x))\,dx\,.
\end{equation}
Using (\ref{FF1})--(\ref{alphaF2*}) it is easy to see that $\F(t)$ is of class $C^1$ on $L^q(\Om;\Rm)$ and that (\ref{F7.5}) and (\ref{(2)}) hold.
By (\ref{F3})--(\ref{star}) for a.e.\ $t\in[0,T]$ the functional $\dot\F(t)$ is of class $C^1$ on 
$L^{\dot q}(\Om;\Rm)$ and
\begin{equation}\label{star1}
\langle\partial \dot\F(t)(u), v\rangle=
\int_\Om \partial_z \dot F(t,x,u(x))\,v(x)\,dx
\end{equation}
for every $u,\,v\in L^{\dot q}(\Om;\Rm)$, so that the functions $t\mapsto \dot \F(t)(u)$ and $t\mapsto \langle \partial \dot \F(t)(u),v\rangle$
are measurable on $[0,T]$ for every $u,\,v\in L^q(\Om;\Rm)$.
{}From (\ref{F5}) and (\ref{F5.5}) we obtain (\ref{F8}) and (\ref{F8.5})
for every pair of simple functions $u$ and $v$ from $\Om$ into $\Rm$.
An easy approximation argument shows that (\ref{F8}) and (\ref{F8.5}) hold for every $u,\, v\in L^{q}(\Om;\Rm)$. Inequalities 
(\ref{F9})--(\ref{F11.5}) follow immediately from (\ref{F7}), (\ref{F7.5}), (\ref{alphaF*}), (\ref{alphaF2*}), and (\ref{F6*})--(\ref{star1}).

\begin{remark}\label{dl}
Let us check that, under the hypotheses on $F_0$ considered above, if the function $t\mapsto f(t,\cdot)$ is absolutely continuous from $[0,T]$ into $L^{\dot q'}(\Om;\Rm)$ and $1<\dot q<1$, then the function $F$ defined by (\ref{dlf}) satisfies (\ref{FF1})--(\ref{F6.5*}) with $\dot F(t,x,z):=\dot f(t,x)z$. Properties (\ref{FF2}) and (\ref{F4}) are trivial; (\ref{FF1}) and (\ref{F5})--(\ref{F3})  follow from well-known properties of absolutely continuous functions with values in reflexive Banach spaces (see, e.g., \cite[Appendix]{Bre}).
By the Cauchy inequality we have
$$
-F(t,x,z)\ge 
\frac{a_0}{q'}|z|^{q}-\frac{1}{q'}a_0^{\frac{1}{1-q}}|f(t,x)|^{q'}-b_0(x)\,,
$$
so that the first inequality in (\ref{alphaF*}) is satisfied with 
$a_0^F:=a_0/q'$ and 
$b_0^F(t,x):=b_0(x)+(a_0^{1/(1-q)}/q')|f(t,x)|^{q'}$. The 
second inequality in (\ref{alphaF*}) can be obtained in a similar way, 
while (\ref{alphaF02}) yields (\ref{alphaF2*}) with $a_2^F:=a_2$ and 
$b_2^F(t,x):=b_2(x)+|f(t,x)|$.

To prove (\ref{F6*}) we observe that by the Cauchy inequality we have
$$
|\dot F(t,x,z)|=|\dot f(t,x)z|\le
\frac{1}{\dot q'}
\frac{|\dot f(t,x)|^{\dot q'}}{\|\dot f(t,\cdot)\|_{\dot q'}^{\dot q'-1}}+
\frac{1}{\dot q}\|\dot f(t,\cdot)\|_{\dot q'}\, |z|^{\dot q} \,,
$$
so that inequality (\ref{F6*}) is satisfied with $b_3^F(t,x):=(1/\dot q')
|\dot f(t,x)|^{\dot q'}\|\dot f(t,\cdot)\|_{\dot q'}^{1-\dot q'}$ and 
$a_3^F(t)\allowbreak:=(1/\dot q)\|\dot f(t,\cdot)\|_{\dot q'}$. Finally, 
(\ref{F6.5*}) holds with $a_4^F(t,x):=0$ and 
$b_4^F(t,x):=|\dot f(t,x)|$.
\end{remark}

\subsection{Weaker hypotheses on the surface forces}\label{newsurface}
In the case of a time dependent dead load, the density $g\colon[0,T]{\times}\partial_S\Om\to \Rm$ of
the applied surface force per unit area in the reference
configuration is independent of the deformation $u$.
Then,  the simplest choice for the
potential is $G(t,x,z):=g(t,x)z$. Let $r$ and $r'$ be as in Subsection~\ref{surfacef}.
If $t\mapsto g(t,\cdot)$ is absolutely
continuous from $[0,T]$ into $L^{r'}(\partial_S\Om;\Rm)$, and 
$t\mapsto\dot g(t,\cdot)$ is its time derivative, which belongs to
$L^1([0,T];L^{r'}(\partial_S\Om;\Rm))$, we can consider
 the functional $\G(t)\colon L^r(\partial_S\Om;\Rm)\to\R$ defined  for every $t\in[0,T]$ by 
$$
\G(t)(u):=\int_{\partial_S\Om} g(t,x)\,u(x)\,d\hn(x)\,,
$$
and the functional 
$\dot\G(t)\colon L^r(\Om;\Rm)\to\R$ defined for a.e.\ $t\in[0,T]$ by
$$
\dot\G(t)(u):=\int_{\partial_S\Om} \dot g(t,x)\,u(x)\,d\hn(x)\,.
$$
We can easily check that in this case $\G$ and $\dot\G$ satisfy all properties required in Subsection~\ref{surfacef}.

Note that these hypotheses do not guarantee the existence of the partial derivative $\partial_t G(t,x,z)$ for a.e.\ $t\in[0,T]$ and for every $(x,z)\in\partial_S\Om{\times}\Rm$, so that Remark~\ref{01091} can not be applied.

We present now a more general set of hypotheses which includes this case as well as those considered in Remark~\ref{01091}. More precisely, we assume that $G\colon[0,T]{\times}\partial_S\Om{\times}\Rm\to\R$ satisfies the following conditions:
\begin{eqnarray}
&&\hbox{\hspace{-1.3cm}for every $z\in\Rm$ the function
$(t,x)\mapsto G(t,x,z)$ is ${\mathcal L}^1{\times}
\hn$-measurable,}\label{GG1}\\
&&\hbox{\hspace{-1.3cm}for every
$(t,x)\in[0,T]{\times}\partial_S\Om$ the function $z\mapsto
G(t,x,z)$
belongs to $C^1(\Rm)$.}\label{GG2}
\end{eqnarray}
Moreover, we assume that there exist two constants $a_1^G\ge 0$, $a_2^G\ge0$ and four nonnegative functions 
$a_0^G\in L^\infty([0,T];L^{r'}(\partial_S\Om))$,
$b_0^G$, $b_1^G\in C^0([0,T];L^1(\partial_S\Om))$,
and $b_2^G\in C^0([0,T];L^{r'}(\partial_S\Om))$
such that
\begin{eqnarray}
& -a_0^G(t,x) |z| - b_0^G(t,x)\le - G(t,x,z) \le
a_1^G  |z|^{r} + b_1^G(t,x)\,,\label{alphaG*}
\\
& |\partial_z G(t,x,z)| \le a_2^G  |z|^{r-1} +
b_2^G(t,x) \label{alphaG2*}
\end{eqnarray}
for every $(t,x,z)\in [0,T]{\times}\partial_S\Om{\times}\Rm$.

To deal with the dependence of $G$ on $t$, we assume that there exists a function \break$\dot G\colon[0,T]{\times}\partial_S\Om{\times}\Rm\to \R$ 
such that for every $t\in [0,T]$ and for every $z\in\Rm$
\begin{eqnarray}
&\displaystyle
G(t,x,z)=G(0,x,z)+\int_0^t\dot G(s,x,z)\,ds
\quad\hbox{for $\hn$-a.e.\ $x\in\partial_S\Om$,}\label{G5}\\
&\displaystyle
\partial_z G(t,x,z)=\partial_z G(0,x,z)+\int_0^t \partial_z \dot 
G(s,x,z)\,ds \quad\hbox{for $\hn$-a.e.\ $x\in\partial_S\Om$.}\label{G5.5}
\end{eqnarray}
The integrals in (\ref{G5}) and (\ref{G5.5}) are well defined since we
assume also that
\begin{eqnarray}
&&\hbox{\hspace{-1.3cm}for every $z\in\Rm$ the function 
$(t,x)\mapsto \dot G(t,x,z)$ is ${\mathcal L}^1{\times} \hn$-measurable,}\label{G3}\\
&&\hbox{\hspace{-1.3cm}for every $(t,x)\in[0,T]{\times}\partial_S\Om$ the function $z\mapsto \dot
G(t,x,z)$
belongs to $C^1(\Rm)$,}\label{G4}
\end{eqnarray}
and that there exist four nonnegative functions
$a_3^G,a_4^G\!\in\! L^1([0,T])$,
$b_3^G\!\in\!L^1([0,T];L^1(\partial_S\Om))$, and
$b_4^G\in L^1([0,T], L^{r'}(\Om))$ such that
\begin{eqnarray}
&|\dot G(t,x,z)|\le
a_3^G(t)|z|^{r}+b_3^G(t,x)\,,\label{G6*}\\
&|\partial_z \dot G(t,x,z)|\le a_4^G(t)|z|^{r-1}+b_4^G(t,x) \label{G6.5*}
\end{eqnarray}
for every $(t,x,z)\in [0,T]{\times}\partial_S\Om{\times}\Rm$.

The functionals $\G(t)$ and $\dot \G(t)\colon L^{r}(\partial_S\Om;\Rm)\to\R$
are now defined by (\ref{G7}) and by
$$
\dot\G(t)(u):=\int_{\partial_S\Om} \dot G(t,x,u(x))\,d\hn(x)\,.
$$
Arguing as in Subsection~\ref{newbody} it is easy to prove that $\G$ and $\dot\G$ satisfy all properties required in Subsection~\ref{surfacef}.

\begin{remark}\label{dls} As in Remark~\ref{dl} we can prove that, if the function $t\mapsto g(t,\cdot)$ is absolutely
continuous from $[0,T]$ into $L^{r'}(\partial_S\Om;\Rm)$, then 
the function $G(t,x,z):=g(t,x)z$
satisfies (\ref{GG1})--(\ref{G6.5*}) with 
$\dot G(t,x,z):=\dot g(t,x)z$.
\end{remark}

\end{section}

\bigskip

\noindent {\bf Acknowledgments.} { This
research began while G.A.\ Francfort
 was visiting
SISSA, whose support is gratefully acknowledged. G.A.\ Francfort also wishes
to warmly acknowledge  C.J.\ Larsen, his co-author for \cite{Fra-Lar}, for  his  amazing insight
and invaluable contribution to the subject. R.\ Toader wishes to thank SISSA for its support during the preparation of the paper. The work of G.\ Dal Maso and 
R.\ Toader is part of the Project ``Calculus of Variations" 2002, supported by the Italian Ministry of Education, University, and Research.}

{
}

\end{document}